\documentclass{amsart}

\usepackage[english]{babel}
\usepackage[latin1]{inputenc}
\usepackage{amscd,amssymb,amsmath,amsfonts,amsthm}

\usepackage{graphicx}

\makeatletter
\newcounter{maxima@count}
\newcommand{\maximacount}{\value{maxima@count}}
\def\maxima@def{mxp\romannumeral\maximacount}   

\newcommand{\maxima@macro}[1]{%
  \gdef\maxima@current{#1}%
  \@ifundefined{#1}%
    {\PackageWarning{maxima}{Undefined result macro `#1'}}%
    {}%
}
\newcommand{\maximacurrent}{\csname\maxima@current\endcsname}

\newenvironment{maximaenv}[1][\maxima@def]{%
  \maxima@macro{#1}%
  \stepcounter{maxima@count}%
  \maximacurrent
}{}

\newcommand{\imaximaenv}[1][\maxima@def]{%
  \maxima@macro{#1}%
  \stepcounter{maxima@count}%
  \maximacurrent
}

\newcommand{\imaximaenvenv}[1][\maxima@def]{%
  \stepcounter{maxima@count}%
}

\makeatother

\newcommand{\CuentasMaxima}[1]{
\textrm{
$
\begin{maximaenv}
\end{maximaenv}
$
}
}

\newcommand{\imt}[1]{\imaximaenv}

\usepackage[active]{srcltx}

\usepackage[small,nohug,heads=vee]{diagrams2024}
\usepackage{hyperref}

\usepackage[x11names]{xcolor}

\usepackage{rotating}

\newtheorem{Theorem}{Theorem}

\makeatletter
\newcommand*\bigcdot{\mathpalette\bigcdot@{.5}}
\newcommand*\bigcdot@[2]{\mathbin{\vcenter{\hbox{\scalebox{#2}{$\m@th#1\bullet$}}}}}
\makeatother

\newcommand{\footnoteocultar}[1]{}

\newcommand{\quitar}[1]{
}

\newcommand{\ojoborrar}[1]{
}

\newcommand{\imparocultar}[1]{
}

\def\cets{CETs}
\def\cet{CET}

\def\O{\mathcal{O}}

\newcommand{\permuclaseinverh}[4]{
\begin{array}{c}
\pi_{#1,#2}\\
\pi_{#1,#3}^{-1}\\
\pi_{#1,#4}^h
\end{array}
}

\newcommand{\permuclaseinver}[3]{
\begin{array}{c}
\pi_{#1,#2}\\
\pi_{#1,#3}^{-1}
\end{array}}

\newcommand{\permuclaseh}[3]{
\begin{array}{c}
\pi_{#1,#2}\\
\pi_{#1,#3}^{h}
\end{array}}

\def\hole{\textrm{hole}}

\def\+{\oplus}
\def\-{\ominus}
\def\X{\otimes}

\newcommand{\antin}{{|\pi|^{-1}(n)}}

\makeatletter
\def\@seccntformat#1{\csname the#1\endcsname.\quad}
\makeatother

\newtheorem{theorem}{Theorem}

\newtheorem*{mteo}{Main Theorem}
\newtheorem{coro}[theorem]{Corollary}
\newtheorem{lemma}[theorem]{Lemma}
\newtheorem{proposition}[theorem]{Proposition}

\newtheorem{propo}[theorem]{Proposition}

\theoremstyle{definition}

\theoremstyle{remark}

\newtheorem{claim}[theorem]{Claim}
\newtheorem{remark}[theorem]{Remark}

\def\diam{\mathop{\mathrm{diam}}}

\def\menos{\backslash}

\def\R{\mathbb{R}}

\def\N{\mathbb{N}}
\def\S{\mathbb{S}}
\def\R{\mathbb{R}}
\def\Z{\mathbb{Z}}

\def\intervalog{[0,\ell]}
\def\intervalogab{(0,\ell)}  

\def\intervalogp{[0,\ell']}   
\def\irreducibleset{S_n^{\sigma,*}}
\def\permuset{S_n^{\sigma}}
\def\permusetplus{S_n^{\sigma,+}}

\def\c1nuno{c_1(n+1)}

\newcommand{\iet}{IET}
\newcommand{\iets}{IETs}
\newcommand{\nkiet}{$(n,k)$-IET}
\newcommand{\nkiets}{$(n,k)$-IETs}

\def\cets{CETs}
\def\cet{CET}

\renewcommand\S{\mathbb{S}}
\def\codes{\mathcal{C}_n}
\def\R{\mathbb{R}}

\def\intervalog{[0,\ell]} 
\def\intervalogab{(0,\ell)} 

\def\O{\mathcal{O}}

\def\Disc{\mathrm{Disc}}

\def\Dom{\mathrm{Dom}}

\title{On the non existence of minimal non-uniquely ergodic  4-\iets\ with flips}

\author{Antonio Linero Bas}

\address{Departamento de Matem\'aticas\\
Universidad de Murcia (Campus de Espinardo)\\
30100-Espinardo-Murcia (Spain)
} 

\email{lineroba@um.es}

\author{
Gabriel Soler L\'{o}pez
}

\address{
Departamento de Matem\'atica Aplicada y Estad\'\i{}stica\\
Universidad Polit\'ecnica de Cartagena (Campus Alfonso XIII)\\
Paseo Alfonso XIII, 52\\
30203-Cartagena (Spain)}

\email{gabriel.soler@upct.es}

\def\negro{\color{black}}

\def\blanco{\color{white}}

\begin{document}

\CuentasMaxima{
clase4norientable: [20,[[[-4,-3,-2,-1]],[[-4,-3,2,-1]],[[-4,-2,-3,-1]],[[-4,-2,3,-1]],[[-4,-3,-2,1],[-4,-3,1,2],[-4,1,2,3],[-3,-4,1,2],[-3,-2,-4,1],
              [-2,-4,-3,1]],[[-4,2,-3,-1]],[[-4,2,3,-1]],
             [[-4,-3,2,1],[-4,2,3,1],[-4,3,1,2],[-3,-2,4,1],[-3,4,1,2],
              [-3,4,2,1],[-2,-4,3,1],[-2,3,4,1]],[[-4,3,-2,-1]],[[-4,3,2,-1]],
             [[-4,3,2,1],[-3,2,4,1],[-2,4,1,3],[-2,4,3,1]],
             [[-4,-1,-2,3],[-4,-1,3,-2],[-4,1,-2,-3],[-4,2,-1,-3],[-4,3,-1,2],
              [-4,3,1,-2],[-3,-4,-1,-2],[-3,1,4,-2],[-3,2,-4,-1],[-3,4,-1,2],
              [-3,4,2,-1],[-2,-4,3,-1],[-2,-3,4,-1],[2,-4,-1,3],[2,-4,3,1],
              [2,-3,-4,-1],[2,4,-3,1],[2,4,1,-3],[3,-4,1,-2],[3,-4,2,-1],
              [3,-2,4,1],[3,1,-4,2],[3,2,-4,1],[4,-3,2,1],[4,-2,1,3],
              [4,1,-3,2],[4,1,3,-2],[4,2,1,-3],[4,3,-2,1]],
             [[-4,-3,-1,-2],[-4,-3,-1,2],[-4,-3,1,-2],[-4,-2,-1,-3],
              [-4,-2,-1,3],[-4,-1,-3,-2],[-4,-1,2,3],[-4,1,-3,-2],[-4,1,2,-3],
              [-3,-4,-2,-1],[-3,-4,-1,2],[-3,-2,-4,-1],[-3,-2,4,-1],
              [-3,1,-4,-2],[-3,4,-2,-1],[-3,4,-1,-2],[-2,-4,-3,-1],
              [-2,3,4,-1],[2,-4,-3,-1],[2,-4,-3,1],[2,-4,-1,-3],[2,-4,1,-3],
              [2,3,-4,-1],[2,3,-4,1],[3,-4,-2,-1],[3,-4,1,2],[3,-2,-4,1],
              [3,1,-4,-2],[3,4,-2,1],[3,4,1,-2],[4,-3,-2,1],[4,-3,1,2],
              [4,-2,1,-3],[4,1,-3,-2],[4,1,2,-3],[4,2,-3,1]],
             [[-4,-1,-2,-3],[-4,-1,3,2],[-4,1,-2,3],[-4,2,1,-3],[-4,3,-2,1],
              [-3,-4,1,-2],[-3,-1,4,-2],[-3,-1,4,2],[-3,2,-4,1],[-2,-4,-1,3],
              [-2,-3,-4,-1],[-2,4,-3,1],[-2,4,-1,3],[-2,4,1,-3],[-2,4,3,-1],
              [2,-3,4,-1],[2,-3,4,1],[3,-4,-1,-2],[3,-4,2,1],[3,-1,-4,2],
              [3,2,-4,-1],[4,-3,2,-1],[4,-2,3,1],[4,-1,-3,2],[4,1,-2,3],
              [4,2,-1,-3],[4,3,1,-2]],
             [[2,3,4,-1],[3,4,-2,-1],[3,4,-1,-2],[4,-3,-2,-1],[4,-2,-1,-3],
              [4,-1,-3,-2]],
             [[2,3,4,1],[3,4,1,2],[3,4,2,1],[4,1,2,3],[4,2,3,1],[4,3,1,2]],
             [[2,4,1,3],[2,4,3,1],[3,1,4,2],[3,2,4,1],[4,1,3,2],[4,2,1,3],
              [4,3,2,1]],[[3,-1,4,2],[4,-1,3,2],[4,2,-1,3],[4,3,2,-1]],
             [[-4,-2,-3,1],[-4,-2,1,-3],[-4,-2,1,3],[-4,-2,3,1],[-4,-1,-3,2],
              [-4,-1,2,-3],[-4,1,-3,2],[-4,1,3,-2],[-4,1,3,2],[-4,2,-3,1],
              [-4,2,-1,3],[-4,2,1,3],[-4,3,-1,-2],[-3,-4,-2,1],[-3,-4,2,-1],
              [-3,-4,2,1],[-3,-1,-4,-2],[-3,-1,-4,2],[-3,1,-4,2],[-3,1,4,2],
              [-3,2,4,-1],[-3,4,-2,1],[-3,4,1,-2],[-2,-4,-1,-3],[-2,-4,1,-3],
              [-2,-4,1,3],[-2,-3,-4,1],[-2,-3,4,1],[-2,3,-4,-1],[-2,3,-4,1],
              [-2,4,-3,-1],[-2,4,-1,-3],[2,-4,1,3],[2,-4,3,-1],[2,-3,-4,1],
              [2,4,-3,-1],[2,4,-1,-3],[2,4,-1,3],[2,4,3,-1],[3,-4,-2,1],
              [3,-4,-1,2],[3,-2,-4,-1],[3,-2,4,-1],[3,-1,-4,-2],[3,-1,4,-2],
              [3,1,4,-2],[3,2,4,-1],[4,-3,-1,-2],[4,-3,-1,2],[4,-3,1,-2],
              [4,-2,-3,-1],[4,-2,-3,1],[4,-2,3,-1],[4,-1,-2,-3],[4,-1,-2,3],
              [4,-1,2,-3],[4,1,-2,-3],[4,2,-3,-1],[4,3,-1,-2]],
             [[3,4,-1,2],[3,4,2,-1],[4,-2,-1,3],[4,-1,2,3],[4,-1,3,-2],
              [4,2,3,-1],[4,3,-2,-1],[4,3,-1,2]]]],
              cl4:clase4norientable[2]
}

\CuentasMaxima{clase11:cl4[11]}
\CuentasMaxima{clase1:cl4[1]}
\CuentasMaxima{clase2:cl4[2]}
\CuentasMaxima{clase3:cl4[3]}
\CuentasMaxima{clase4:cl4[4]}
\CuentasMaxima{clase5:cl4[5]}
\CuentasMaxima{clase6:cl4[6]}
\CuentasMaxima{clase7:cl4[7]}
\CuentasMaxima{clase8:cl4[8]}
\CuentasMaxima{clase9:cl4[9]}
\CuentasMaxima{clase10:cl4[10]}
\CuentasMaxima{clase11:cl4[11]}
\CuentasMaxima{clase12:cl4[12]}
\CuentasMaxima{clase13:cl4[13]}
\CuentasMaxima{clase14:cl4[14]}
\CuentasMaxima{clase15:cl4[15]}
\CuentasMaxima{clase16:cl4[16]}
\CuentasMaxima{clase17:cl4[17]}
\CuentasMaxima{clase18:cl4[18]}
\CuentasMaxima{clase19:cl4[19]}
\CuentasMaxima{clase20:cl4[20]}

%

\def\guion{--}

\newarrow{Corresponds}<--->
\newarrow{Sinflecha}{-}{-}{-}{-}{-}
\newarrow{Tachado}{-}{-}{-}{-}{>}
\newarrow{Noantecesor}--+->


\newcommand{\puntoblanco}{\blanco .... \negro}

\begin{abstract}

The existence of minimal non uniquely ergodic oriented \iets\ of $n$ intervals requires that $n\geq 4$ by \cite[Th. 2.12]{veech}; moreover,  for $n=4$ there exists  examples of  minimal non uniquely ergodic $4$-\iets\ with two independent invariant measures, see \cite{keane1}. For the non oriented case, also called \iets\ with flips, the first examples of minimal non uniquely ergodic \iets\ were built in \cite{NoUnicamenteErgodica10}. However, the examples are $(10,k)$-\iets\ ($k$ is the number of flips, $1\leq k\leq 10$). Since there are no minimal $(n,k)$-\iets\ for $n\leq 3,$ it is interesting to decide if, as in the orientable case, there exist minimal non uniquely ergodic \iets\ with flips and four intervals. 
Our Main Theorem answers this question in the negative.

Documento de trabajo: \jobname 

Fecha de generaci\'on: \today

\end{abstract}

\subjclass[2010]{Primary 37E05. Secondary 37A05, 37B05, 28D05}

\keywords{Interval  exchange transformations, flips, minimality, transitivity, Rauzy induction,  invariant measures, non uniquely ergodic, ergodicity, nonnegative matrices}

\maketitle


\definecolor{RED}{rgb}{1,0,0}

\definecolor{red}{rgb}{1,0,0}

\gdef\mxpi { 
} 
\gdef\mxpii { 
} 
\gdef\mxpiii { 
} 
\gdef\mxpiv { 
} 
\gdef\mxpv { 
} 
\gdef\mxpvi { 
} 
\gdef\mxpvii { 
} 
\gdef\mxpviii { 
} 
\gdef\mxpix { 
} 
\gdef\mxpx { 
} 
\gdef\mxpxi { 
} 
\gdef\mxpxii { 
} 
\gdef\mxpxiii { 
} 
\gdef\mxpxiv { 
} 
\gdef\mxpxv { 
} 
\gdef\mxpxvi { 
} 
\gdef\mxpxvii { 
} 
\gdef\mxpxviii { 
} 
\gdef\mxpxix { 
} 
\gdef\mxpxx { 
} 
\gdef\mxpxxi { 
} 
\gdef\mxpxxii { 
} 
\gdef\mxpxxiii { 
\left( -4 , -3 , -2 , -1 \right) 
} 
\gdef\mxpxxiv { 
\left( -4 , -3 , 2 , -1 \right) 
} 
\gdef\mxpxxv { 
\left( -4 , -2 , -3 , -1 \right) 
} 
\gdef\mxpxxvi { 
\left( -4 , -2 , 3 , -1 \right) 
} 
\gdef\mxpxxvii { 
\left( -4 , 2 , -3 , -1 \right) 
} 
\gdef\mxpxxviii { 
\left( -4 , 2 , 3 , -1 \right) 
} 
\gdef\mxpxxix { 
\left( -4 , 3 , -2 , -1 \right) 
} 
\gdef\mxpxxx { 
\left( -4 , 3 , 2 , -1 \right) 
} 
\gdef\mxpxxxi { 
} 
\gdef\mxpxxxii { 
} 
\gdef\mxpxxxiii { 
\left( 2 , 3 , 1 \right) 
} 
\gdef\mxpxxxiv { 
\left( 3 , 1 , 2 \right) 
} 
\gdef\mxpxxxv { 
\left( 3 , 2 , 1 \right) 
} 
\gdef\mxpxxxvi { 
\left( -3 , -2 , -1 \right) 
} 
\gdef\mxpxxxvii { 
\left( -3 , -2 , 1 \right) 
} 
\gdef\mxpxxxviii { 
\left( -3 , 1 , 2 \right) 
} 
\gdef\mxpxxxix { 
\left( -2 , -3 , 1 \right) 
} 
\gdef\mxpxl { 
\left( -3 , 2 , -1 \right) 
} 
\gdef\mxpxli { 
\left( -3 , 2 , 1 \right) 
} 
\gdef\mxpxlii { 
\left( -2 , 3 , 1 \right) 
} 
\gdef\mxpxliii { 
\left( -3 , -1 , -2 \right) 
} 
\gdef\mxpxliv { 
\left( -3 , -1 , 2 \right) 
} 
\gdef\mxpxlv { 
\left( -3 , 1 , -2 \right) 
} 
\gdef\mxpxlvi { 
\left( -2 , -3 , -1 \right) 
} 
\gdef\mxpxlvii { 
\left( -2 , 3 , -1 \right) 
} 
\gdef\mxpxlviii { 
\left( 2 , -3 , -1 \right) 
} 
\gdef\mxpxlix { 
\left( 2 , -3 , 1 \right) 
} 
\gdef\mxpl { 
\left( 3 , -2 , 1 \right) 
} 
\gdef\mxpli { 
\left( 3 , 1 , -2 \right) 
} 
\gdef\mxplii { 
\left( 2 , 3 , -1 \right) 
} 
\gdef\mxpliii { 
\left( 3 , -2 , -1 \right) 
} 
\gdef\mxpliv { 
\left( 3 , -1 , -2 \right) 
} 
\gdef\mxplv { 
\left( 2 , 3 , 1 \right) 
} 
\gdef\mxplvi { 
\left( 3 , 1 , 2 \right) 
} 
\gdef\mxplvii { 
\left( 3 , 2 , 1 \right) 
} 
\gdef\mxplviii { 
\left( 3 , -1 , 2 \right) 
} 
\gdef\mxplix { 
\left( 3 , 2 , -1 \right) 
} 
\gdef\mxplx { 
\left( 2 , 3 , 4 , 1 \right) 
} 
\gdef\mxplxi { 
\left( 3 , 4 , 1 , 2 \right) 
} 
\gdef\mxplxii { 
\left( 3 , 4 , 2 , 1 \right) 
} 
\gdef\mxplxiii { 
\left( 4 , 1 , 2 , 3 \right) 
} 
\gdef\mxplxiv { 
\left( 4 , 2 , 3 , 1 \right) 
} 
\gdef\mxplxv { 
\left( 4 , 3 , 1 , 2 \right) 
} 
\gdef\mxplxvi { 
\left( 2 , 4 , 1 , 3 \right) 
} 
\gdef\mxplxvii { 
\left( 2 , 4 , 3 , 1 \right) 
} 
\gdef\mxplxviii { 
\left( 3 , 1 , 4 , 2 \right) 
} 
\gdef\mxplxix { 
\left( 3 , 2 , 4 , 1 \right) 
} 
\gdef\mxplxx { 
\left( 4 , 1 , 3 , 2 \right) 
} 
\gdef\mxplxxi { 
\left( 4 , 2 , 1 , 3 \right) 
} 
\gdef\mxplxxii { 
\left( 4 , 3 , 2 , 1 \right) 
} 
\gdef\mxplxxiii { 
\left( -4 , 3 , 2 , 1 \right) 
} 
\gdef\mxplxxiv { 
\left( -3 , 2 , 4 , 1 \right) 
} 
\gdef\mxplxxv { 
\left( -2 , 4 , 1 , 3 \right) 
} 
\gdef\mxplxxvi { 
\left( -2 , 4 , 3 , 1 \right) 
} 
\gdef\mxplxxvii { 
\left( 3 , -1 , 4 , 2 \right) 
} 
\gdef\mxplxxviii { 
\left( 4 , -1 , 3 , 2 \right) 
} 
\gdef\mxplxxix { 
\left( 4 , 2 , -1 , 3 \right) 
} 
\gdef\mxplxxx { 
\left( 4 , 3 , 2 , -1 \right) 
} 
\gdef\mxplxxxi { 
\left( -4 , -3 , 1 , 2 \right) 
} 
\gdef\mxplxxxii { 
\left( -2 , -4 , -3 , 1 \right) 
} 
\gdef\mxplxxxiii { 
\left( -3 , -4 , 1 , 2 \right) 
} 
\gdef\mxplxxxiv { 
\left( -3 , -2 , -4 , 1 \right) 
} 
\gdef\mxplxxxv { 
\left( -4 , -3 , -2 , 1 \right) 
} 
\gdef\mxplxxxvi { 
\left( -4 , 1 , 2 , 3 \right) 
} 
\gdef\mxplxxxvii { 
\left( 4 , -3 , -2 , -1 \right) 
} 
\gdef\mxplxxxviii { 
\left( 3 , 4 , -2 , -1 \right) 
} 
\gdef\mxplxxxix { 
\left( 2 , 3 , 4 , -1 \right) 
} 
\gdef\mxpxc { 
\left( 3 , 4 , -1 , -2 \right) 
} 
\gdef\mxpxci { 
\left( 4 , -1 , -3 , -2 \right) 
} 
\gdef\mxpxcii { 
\left( 4 , -2 , -1 , -3 \right) 
} 
\gdef\mxpxciii { 
\left( -4 , 2 , 3 , 1 \right) 
} 
\gdef\mxpxciv { 
\left( -4 , -3 , 2 , 1 \right) 
} 
\gdef\mxpxcv { 
\left( -3 , -2 , 4 , 1 \right) 
} 
\gdef\mxpxcvi { 
\left( -2 , -4 , 3 , 1 \right) 
} 
\gdef\mxpxcvii { 
\left( -2 , 3 , 4 , 1 \right) 
} 
\gdef\mxpxcviii { 
\left( -3 , 4 , 2 , 1 \right) 
} 
\gdef\mxpxcix { 
\left( -3 , 4 , 1 , 2 \right) 
} 
\gdef\mxpc { 
\left( -4 , 3 , 1 , 2 \right) 
} 
\gdef\mxpci { 
\left( -4 , -1 , -2 , 3 \right) 
} 
\gdef\mxpcii { 
\left( -4 , -1 , 3 , -2 \right) 
} 
\gdef\mxpciii { 
\left( -4 , 1 , -2 , -3 \right) 
} 
\gdef\mxpciv { 
\left( -4 , 2 , -1 , -3 \right) 
} 
\gdef\mxpcv { 
\left( -4 , 3 , -1 , 2 \right) 
} 
\gdef\mxpcvi { 
\left( -4 , 3 , 1 , -2 \right) 
} 
\gdef\mxpcvii { 
\left( -3 , -4 , -1 , -2 \right) 
} 
\gdef\mxpcviii { 
\left( -3 , 1 , 4 , -2 \right) 
} 
\gdef\mxpcix { 
\left( -3 , 2 , -4 , -1 \right) 
} 
\gdef\mxpcx { 
\left( -3 , 4 , -1 , 2 \right) 
} 
\gdef\mxpcxi { 
\left( -3 , 4 , 2 , -1 \right) 
} 
\gdef\mxpcxii { 
\left( -2 , -4 , 3 , -1 \right) 
} 
\gdef\mxpcxiii { 
\left( -2 , -3 , 4 , -1 \right) 
} 
\gdef\mxpcxiv { 
\left( 2 , -4 , -1 , 3 \right) 
} 
\gdef\mxpcxv { 
\left( 2 , -4 , 3 , 1 \right) 
} 
\gdef\mxpcxvi { 
\left( 2 , -3 , -4 , -1 \right) 
} 
\gdef\mxpcxvii { 
\left( 2 , 4 , -3 , 1 \right) 
} 
\gdef\mxpcxviii { 
\left( 2 , 4 , 1 , -3 \right) 
} 
\gdef\mxpcxix { 
\left( 3 , -4 , 1 , -2 \right) 
} 
\gdef\mxpcxx { 
\left( 3 , -4 , 2 , -1 \right) 
} 
\gdef\mxpcxxi { 
\left( 3 , -2 , 4 , 1 \right) 
} 
\gdef\mxpcxxii { 
\left( 3 , 1 , -4 , 2 \right) 
} 
\gdef\mxpcxxiii { 
\left( 3 , 2 , -4 , 1 \right) 
} 
\gdef\mxpcxxiv { 
\left( 4 , -3 , 2 , 1 \right) 
} 
\gdef\mxpcxxv { 
\left( 4 , -2 , 1 , 3 \right) 
} 
\gdef\mxpcxxvi { 
\left( 4 , 1 , -3 , 2 \right) 
} 
\gdef\mxpcxxvii { 
\left( 4 , 1 , 3 , -2 \right) 
} 
\gdef\mxpcxxviii { 
\left( 4 , 2 , 1 , -3 \right) 
} 
\gdef\mxpcxxix { 
\left( 4 , 3 , -2 , 1 \right) 
} 
\gdef\mxpcxxx { 
} 
\gdef\mxpcxxxi { 
} 
\gdef\mxpcxxxii { 
\left \{15 , 17 , 18 , 21 , 22 , 23 , 24 , 25 , 26 , 27 , 28 , 29
  \right \}
} 
\gdef\mxpcxxxiii { 
} 
\gdef\mxpcxxxiv { 
\left \{1 , 2 , 3 , 4 , 5 , 6 , 7 , 8 , 9 , 10 , 11 , 12 , 13 , 14 , 
 15 , 16 , 17 , 18 , 19 , 20 , 21 , 22 , 23 , 24 , 25 , 26 , 27 , 28
  , 29 \right \}
} 
\gdef\mxpcxxxv { 
\left( -4 , -1 , -2 , -3 \right) 
} 
\gdef\mxpcxxxvi { 
\left( -4 , -1 , 3 , 2 \right) 
} 
\gdef\mxpcxxxvii { 
\left( -4 , 1 , -2 , 3 \right) 
} 
\gdef\mxpcxxxviii { 
\left( -4 , 2 , 1 , -3 \right) 
} 
\gdef\mxpcxxxix { 
\left( -4 , 3 , -2 , 1 \right) 
} 
\gdef\mxpcxl { 
\left( -3 , -4 , 1 , -2 \right) 
} 
\gdef\mxpcxli { 
\left( -3 , -1 , 4 , -2 \right) 
} 
\gdef\mxpcxlii { 
\left( -3 , -1 , 4 , 2 \right) 
} 
\gdef\mxpcxliii { 
\left( -3 , 2 , -4 , 1 \right) 
} 
\gdef\mxpcxliv { 
\left( -2 , -4 , -1 , 3 \right) 
} 
\gdef\mxpcxlv { 
\left( -2 , -3 , -4 , -1 \right) 
} 
\gdef\mxpcxlvi { 
\left( -2 , 4 , -3 , 1 \right) 
} 
\gdef\mxpcxlvii { 
\left( -2 , 4 , -1 , 3 \right) 
} 
\gdef\mxpcxlviii { 
\left( -2 , 4 , 1 , -3 \right) 
} 
\gdef\mxpcxlix { 
\left( -2 , 4 , 3 , -1 \right) 
} 
\gdef\mxpcl { 
\left( 2 , -3 , 4 , -1 \right) 
} 
\gdef\mxpcli { 
\left( 2 , -3 , 4 , 1 \right) 
} 
\gdef\mxpclii { 
\left( 3 , -4 , -1 , -2 \right) 
} 
\gdef\mxpcliii { 
\left( 3 , -4 , 2 , 1 \right) 
} 
\gdef\mxpcliv { 
\left( 3 , -1 , -4 , 2 \right) 
} 
\gdef\mxpclv { 
\left( 3 , 2 , -4 , -1 \right) 
} 
\gdef\mxpclvi { 
\left( 4 , -3 , 2 , -1 \right) 
} 
\gdef\mxpclvii { 
\left( 4 , -2 , 3 , 1 \right) 
} 
\gdef\mxpclviii { 
\left( 4 , -1 , -3 , 2 \right) 
} 
\gdef\mxpclix { 
\left( 4 , 1 , -2 , 3 \right) 
} 
\gdef\mxpclx { 
\left( 4 , 2 , -1 , -3 \right) 
} 
\gdef\mxpclxi { 
\left( 4 , 3 , 1 , -2 \right) 
} 
\gdef\mxpclxii { 
} 
\gdef\mxpclxiii { 
\left( -2 , 4 , -3 , 1 \right) 
} 
\gdef\mxpclxiv { 
\left( -2 , 4 , 1 , -3 \right) 
} 
\gdef\mxpclxv { 
\left \{5 , 9 , 12 , 14 , 20 , 22 , 24 , 26 \right \}
} 
\gdef\mxpclxvi { 
\left( 2 , -3 , 4 , 1 \right) 
} 
\gdef\mxpclxvii { 
\left( 4 , -2 , 3 , 1 \right) 
} 
\gdef\mxpclxviii { 
\left( 4 , 3 , 1 , -2 \right) 
} 
\gdef\mxpclxix { 
\left( 4 , 1 , -2 , 3 \right) 
} 
\gdef\mxpclxx { 
} 
\gdef\mxpclxxi { 
\left \{5 , 9 , 12 , 14 , 17 , 19 , 20 , 22 , 23 , 24 , 25 , 26 , 27
  \right \}
} 
\gdef\mxpclxxii { 
} 
\gdef\mxpclxxiii { 
\left \{1 , 2 , 4 , 5 , 6 , 9 , 11 , 12 , 14 , 17 , 18 , 19 , 20 , 21
  , 22 , 23 , 24 , 25 , 26 , 27 \right \}
} 
\gdef\mxpclxxiv { 
} 
\gdef\mxpclxxv { 
\left \{1 , 2 , 3 , 4 , 5 , 6 , 7 , 8 , 9 , 10 , 11 , 12 , 13 , 14 , 
 15 , 16 , 17 , 18 , 19 , 20 , 21 , 22 , 23 , 24 , 25 , 26 , 27
  \right \}
} 
\gdef\mxpclxxvi { 
\left( 3 , 4 , -1 , 2 \right) 
} 
\gdef\mxpclxxvii { 
\left( 3 , 4 , 2 , -1 \right) 
} 
\gdef\mxpclxxviii { 
\left( 4 , -2 , -1 , 3 \right) 
} 
\gdef\mxpclxxix { 
\left( 4 , -1 , 2 , 3 \right) 
} 
\gdef\mxpclxxx { 
\left( 4 , -1 , 3 , -2 \right) 
} 
\gdef\mxpclxxxi { 
\left( 4 , 2 , 3 , -1 \right) 
} 
\gdef\mxpclxxxii { 
\left( 4 , 3 , -2 , -1 \right) 
} 
\gdef\mxpclxxxiii { 
\left( 4 , 3 , -1 , 2 \right) 
} 
\gdef\mxpclxxxiv { 
\left( -3 , 4 , 1 , 2 \right) 
} 
\gdef\mxpclxxxv { 
\left( -4 , 3 , 1 , 2 \right) 
} 
\gdef\mxpclxxxvi { 
\left( -3 , -2 , 4 , 1 \right) 
} 
\gdef\mxpclxxxvii { 
\left( -2 , 3 , 4 , 1 \right) 
} 
\gdef\mxpclxxxviii { 
\left( -3 , 2 , -4 , 1 \right) 
} 
\gdef\mxpclxxxix { 
\left( -4 , 2 , 3 , 1 \right) 
} 
\gdef\mxpcxc { 
\left( -4 , -3 , 2 , 1 \right) 
} 
\gdef\mxpcxci { 
\left( 3 , -4 , 2 , 1 \right) 
} 
\gdef\mxpcxcii { 
\left( -4 , -3 , -1 , -2 \right) 
} 
\gdef\mxpcxciii { 
\left( -4 , -3 , -1 , 2 \right) 
} 
\gdef\mxpcxciv { 
\left( -4 , -3 , 1 , -2 \right) 
} 
\gdef\mxpcxcv { 
\left( -4 , -2 , -1 , -3 \right) 
} 
\gdef\mxpcxcvi { 
\left( -4 , -2 , -1 , 3 \right) 
} 
\gdef\mxpcxcvii { 
\left( -4 , -1 , -3 , -2 \right) 
} 
\gdef\mxpcxcviii { 
\left( -4 , -1 , 2 , 3 \right) 
} 
\gdef\mxpcxcix { 
\left( -4 , 1 , -3 , -2 \right) 
} 
\gdef\mxpcc { 
\left( -4 , 1 , 2 , -3 \right) 
} 
\gdef\mxpcci { 
\left( -3 , -4 , -2 , -1 \right) 
} 
\gdef\mxpccii { 
\left( -3 , -4 , -1 , 2 \right) 
} 
\gdef\mxpcciii { 
\left( -3 , -2 , -4 , -1 \right) 
} 
\gdef\mxpcciv { 
\left( -3 , -2 , 4 , -1 \right) 
} 
\gdef\mxpccv { 
\left( -3 , 1 , -4 , -2 \right) 
} 
\gdef\mxpccvi { 
\left( -3 , 4 , -2 , -1 \right) 
} 
\gdef\mxpccvii { 
\left( -3 , 4 , -1 , -2 \right) 
} 
\gdef\mxpccviii { 
\left( -2 , -4 , -3 , -1 \right) 
} 
\gdef\mxpccix { 
\left( -2 , 3 , 4 , -1 \right) 
} 
\gdef\mxpccx { 
\left( 2 , -4 , -3 , -1 \right) 
} 
\gdef\mxpccxi { 
\left( 2 , -4 , -3 , 1 \right) 
} 
\gdef\mxpccxii { 
\left( 2 , -4 , -1 , -3 \right) 
} 
\gdef\mxpccxiii { 
\left( 2 , -4 , 1 , -3 \right) 
} 
\gdef\mxpccxiv { 
\left( 2 , 3 , -4 , -1 \right) 
} 
\gdef\mxpccxv { 
\left( 2 , 3 , -4 , 1 \right) 
} 
\gdef\mxpccxvi { 
\left( 3 , -4 , -2 , -1 \right) 
} 
\gdef\mxpccxvii { 
\left( 3 , -4 , 1 , 2 \right) 
} 
\gdef\mxpccxviii { 
\left( 3 , -2 , -4 , 1 \right) 
} 
\gdef\mxpccxix { 
\left( 3 , 1 , -4 , -2 \right) 
} 
\gdef\mxpccxx { 
\left( 3 , 4 , -2 , 1 \right) 
} 
\gdef\mxpccxxi { 
\left( 3 , 4 , 1 , -2 \right) 
} 
\gdef\mxpccxxii { 
\left( 4 , -3 , -2 , 1 \right) 
} 
\gdef\mxpccxxiii { 
\left( 4 , -3 , 1 , 2 \right) 
} 
\gdef\mxpccxxiv { 
\left( 4 , -2 , 1 , -3 \right) 
} 
\gdef\mxpccxxv { 
\left( 4 , 1 , -3 , -2 \right) 
} 
\gdef\mxpccxxvi { 
\left( 4 , 1 , 2 , -3 \right) 
} 
\gdef\mxpccxxvii { 
\left( 4 , 2 , -3 , 1 \right) 
} 
\gdef\mxpccxxviii { 
\left( 3 , -4 , -1 , -2 \right) 
} 
\gdef\mxpccxxix { 
\left( -3 , -1 , 4 , -2 \right) 
} 
\gdef\mxpccxxx { 
\left( -3 , -4 , 1 , -2 \right) 
} 
\gdef\mxpccxxxi { 
\left( 4 , -2 , -1 , 3 \right) 
} 
\gdef\mxpccxxxii { 
\left( -2 , -4 , -1 , 3 \right) 
} 
\gdef\mxpccxxxiii { 
\left( -2 , 4 , -1 , 3 \right) 
} 
\gdef\mxpccxxxiv { 
\left( 4 , -1 , 2 , 3 \right) 
} 
\gdef\mxpccxxxv { 
\left( 3 , 4 , -1 , 2 \right) 
} 
\gdef\mxpccxxxvi { 
\left( 4 , -1 , -3 , 2 \right) 
} 
\gdef\mxpccxxxvii { 
\left( -3 , -1 , 4 , 2 \right) 
} 
\gdef\mxpccxxxviii { 
\left( -3 , 4 , 1 , 2 \right) 
} 
\gdef\mxpccxxxix { 
\left( -2 , 4 , -3 , 1 \right) 
} 
\gdef\mxpccxl { 
\left( -3 , -2 , 4 , 1 \right) 
} 
\gdef\mxpccxli { 
\left( -2 , 3 , 4 , 1 \right) 
} 
\gdef\mxpccxlii { 
\left( 4 , -2 , 3 , 1 \right) 
} 
\gdef\mxpccxliii { 
} 
\gdef\mxpccxliv { 
\left \{11 , 14 , 16 , 20 , 21 , 22 , 24 , 26 , 27 , 28 , 30 , 33 , 
 34 , 35 , 36 \right \}
} 
\gdef\mxpccxlv { 
} 
\gdef\mxpccxlvi { 
} 
\gdef\mxpccxlvii { 
\left \{1 , 2 , 3 , 4 , 5 , 6 , 7 , 8 , 9 , 10 , 12 , 13 , 14 , 15 , 
 17 , 18 , 19 , 23 , 25 , 29 , 31 , 32 \right \}
} 
\gdef\mxpccxlviii { 
\left( -4 , -2 , -3 , 1 \right) 
} 
\gdef\mxpccxlix { 
\left( -4 , -2 , 1 , -3 \right) 
} 
\gdef\mxpccl { 
\left( -4 , -2 , 1 , 3 \right) 
} 
\gdef\mxpccli { 
\left( -4 , -2 , 3 , 1 \right) 
} 
\gdef\mxpcclii { 
\left( -4 , -1 , -3 , 2 \right) 
} 
\gdef\mxpccliii { 
\left( -4 , -1 , 2 , -3 \right) 
} 
\gdef\mxpccliv { 
\left( -4 , 1 , -3 , 2 \right) 
} 
\gdef\mxpcclv { 
\left( -4 , 1 , 3 , -2 \right) 
} 
\gdef\mxpcclvi { 
\left( -4 , 1 , 3 , 2 \right) 
} 
\gdef\mxpcclvii { 
\left( -4 , 2 , -3 , 1 \right) 
} 
\gdef\mxpcclviii { 
\left( -4 , 2 , -1 , 3 \right) 
} 
\gdef\mxpcclix { 
\left( -4 , 2 , 1 , 3 \right) 
} 
\gdef\mxpcclx { 
\left( -4 , 3 , -1 , -2 \right) 
} 
\gdef\mxpcclxi { 
\left( -3 , -4 , -2 , 1 \right) 
} 
\gdef\mxpcclxii { 
\left( -3 , -4 , 2 , -1 \right) 
} 
\gdef\mxpcclxiii { 
\left( -3 , -4 , 2 , 1 \right) 
} 
\gdef\mxpcclxiv { 
\left( -3 , -1 , -4 , -2 \right) 
} 
\gdef\mxpcclxv { 
\left( -3 , -1 , -4 , 2 \right) 
} 
\gdef\mxpcclxvi { 
\left( -3 , 1 , -4 , 2 \right) 
} 
\gdef\mxpcclxvii { 
\left( -3 , 1 , 4 , 2 \right) 
} 
\gdef\mxpcclxviii { 
\left( -3 , 2 , 4 , -1 \right) 
} 
\gdef\mxpcclxix { 
\left( -3 , 4 , -2 , 1 \right) 
} 
\gdef\mxpcclxx { 
\left( -3 , 4 , 1 , -2 \right) 
} 
\gdef\mxpcclxxi { 
\left( -2 , -4 , -1 , -3 \right) 
} 
\gdef\mxpcclxxii { 
\left( -2 , -4 , 1 , -3 \right) 
} 
\gdef\mxpcclxxiii { 
\left( -2 , -4 , 1 , 3 \right) 
} 
\gdef\mxpcclxxiv { 
\left( -2 , -3 , -4 , 1 \right) 
} 
\gdef\mxpcclxxv { 
\left( -2 , -3 , 4 , 1 \right) 
} 
\gdef\mxpcclxxvi { 
\left( -2 , 3 , -4 , -1 \right) 
} 
\gdef\mxpcclxxvii { 
\left( -2 , 3 , -4 , 1 \right) 
} 
\gdef\mxpcclxxviii { 
\left( -2 , 4 , -3 , -1 \right) 
} 
\gdef\mxpcclxxix { 
\left( -2 , 4 , -1 , -3 \right) 
} 
\gdef\mxpcclxxx { 
\left( 2 , -4 , 1 , 3 \right) 
} 
\gdef\mxpcclxxxi { 
\left( 2 , -4 , 3 , -1 \right) 
} 
\gdef\mxpcclxxxii { 
\left( 2 , -3 , -4 , 1 \right) 
} 
\gdef\mxpcclxxxiii { 
\left( 2 , 4 , -3 , -1 \right) 
} 
\gdef\mxpcclxxxiv { 
\left( 2 , 4 , -1 , -3 \right) 
} 
\gdef\mxpcclxxxv { 
\left( 2 , 4 , -1 , 3 \right) 
} 
\gdef\mxpcclxxxvi { 
\left( 2 , 4 , 3 , -1 \right) 
} 
\gdef\mxpcclxxxvii { 
\left( 3 , -4 , -2 , 1 \right) 
} 
\gdef\mxpcclxxxviii { 
\left( 3 , -4 , -1 , 2 \right) 
} 
\gdef\mxpcclxxxix { 
\left( 3 , -2 , -4 , -1 \right) 
} 
\gdef\mxpccxc { 
\left( 3 , -2 , 4 , -1 \right) 
} 
\gdef\mxpccxci { 
\left( 3 , -1 , -4 , -2 \right) 
} 
\gdef\mxpccxcii { 
\left( 3 , -1 , 4 , -2 \right) 
} 
\gdef\mxpccxciii { 
\left( 3 , 1 , 4 , -2 \right) 
} 
\gdef\mxpccxciv { 
\left( 3 , 2 , 4 , -1 \right) 
} 
\gdef\mxpccxcv { 
\left( 4 , -3 , -1 , -2 \right) 
} 
\gdef\mxpccxcvi { 
\left( 4 , -3 , -1 , 2 \right) 
} 
\gdef\mxpccxcvii { 
\left( 4 , -3 , 1 , -2 \right) 
} 
\gdef\mxpccxcviii { 
\left( 4 , -2 , -3 , -1 \right) 
} 
\gdef\mxpccxcix { 
\left( 4 , -2 , -3 , 1 \right) 
} 
\gdef\mxpccc { 
\left( 4 , -2 , 3 , -1 \right) 
} 
\gdef\mxpccci { 
\left( 4 , -1 , -2 , -3 \right) 
} 
\gdef\mxpcccii { 
\left( 4 , -1 , -2 , 3 \right) 
} 
\gdef\mxpccciii { 
\left( 4 , -1 , 2 , -3 \right) 
} 
\gdef\mxpccciv { 
\left( 4 , 1 , -2 , -3 \right) 
} 
\gdef\mxpcccv { 
\left( 4 , 2 , -3 , -1 \right) 
} 
\gdef\mxpcccvi { 
\left( 4 , 3 , -1 , -2 \right) 
} 
\gdef\mxpcccvii { 
\left( 4 , -3 , 1 , -2 \right) 
} 
\gdef\mxpcccviii { 
} 
\gdef\mxpcccix { 
\left \{7 , 8 , 11 , 13 , 15 , 19 , 21 , 33 , 34 , 36 , 37 , 46
  \right \}
} 
\gdef\mxpcccx { 
} 
\gdef\mxpcccxi { 
\left \{1 , 3 , 5 , 6 , 12 , 14 , 17 , 18 , 24 , 26 , 27 , 29 , 31 , 
 32 , 43 , 45 , 47 , 48 , 51 , 54 \right \}
} 
\gdef\mxpcccxii { 
} 
\gdef\mxpcccxiii { 
\left \{4 , 9 , 16 , 20 , 28 , 38 , 39 , 53 , 55 , 59 \right \}
} 
\gdef\mxpcccxiv { 
\left( 4 , -1 , -2 , 3 \right) 
} 
\gdef\mxpcccxv { 
\left( 4 , 3 , -1 , -2 \right) 
} 
\gdef\mxpcccxvi { 
\left( 4 , -2 , 3 , -1 \right) 
} 
\gdef\mxpcccxvii { 
\left( -2 , -3 , 4 , 1 \right) 
} 
\gdef\mxpcccxviii { 
\left( -3 , -4 , 2 , 1 \right) 
} 
\gdef\mxpcccxix { 
\left( -4 , -2 , 3 , 1 \right) 
} 
\gdef\mxpcccxx { 
} 
\gdef\mxpcccxxi { 
\left \{2 , 10 , 22 , 23 , 25 , 30 , 35 , 40 , 41 , 42 , 44 , 49 , 50
  , 52 , 56 , 57 , 58 \right \}
} 
\gdef\mxpcccxxii { 
\left( 4 , -3 , 1 , -2 \right) 
} 
\gdef\mxpcccxxiii { 
\left( -2 , 3 , -4 , 1 \right) 
} 
\gdef\mxpcccxxiv { 
\left( 4 , -1 , 2 , -3 \right) 
} 
\gdef\mxpcccxxv { 
\left( -3 , 4 , -2 , 1 \right) 
} 
\gdef\mxpcccxxvi { 
\left( -3 , 4 , 1 , -2 \right) 
} 
\gdef\mxpcccxxvii { 
\left( -4 , 2 , -3 , 1 \right) 
} 
\gdef\mxpcccxxviii { 
\left( -2 , 3 , -4 , 1 \right) 
} 
\gdef\mxpcccxxix { 
\left( 4 , -1 , 2 , -3 \right) 
} 
\gdef\mxpcccxxx { 
\left( 4 , -3 , -1 , 2 \right) 
} 
\gdef\mxpcccxxxi { 
\left( 4 , 2 , -3 , -1 \right) 
} 
\gdef\mxpcccxxxii { 
\left( 3 , -4 , -1 , 2 \right) 
} 
\gdef\mxpcccxxxiii { 
\left( 4 , -3 , 1 , -2 \right) 
} 
\gdef\mxpcccxxxiv { 
\left( 4 , -2 , -3 , 1 \right) 
} 
\gdef\mxpcccxxxv { 
\left( 4 , 1 , -2 , -3 \right) 
} 
\gdef\mxpcccxxxvi { 
\left( 4 , -3 , 1 , -2 \right) 
}

\section{Introduction}

Let $\N$ and $\Z$ denote  the sets of natural \textcolor{black}{(positive integer)} numbers and integers, respectively,  and let $n\in\N$.  An
\emph{n-\iet} is  an injective map $T:D\subset\intervalogab\to \intervalogab$
such that:

\begin{enumerate}
\item $D$ is the union of $n$ pairwise disjoint open intervals, $D=\bigcup_{i=1}^n I_i$, \textcolor{black}{where
$I_i=(a_i,a_{i+1})$, $1\leq i\leq n$, and} $0=a_1<a_2<a_3<\dots<a_{n+1}=\ell$.
\item $T|_{I_i}$ is an affine  map of constant slope of either $1$ or $-1$, for  $i=1,2,\dots, n$.
\end{enumerate}

An \iet\ $T$  is said to have  flips if there exist  subintervals where the slope is $-1$. Otherwise,  $T$  is 
called \emph{an oriented interval exchange transformation}, which can be extended to $[0,\ell)$ right-continously at its  discontinuity points while  preserving the bijectivity. Tipically,  papers
that focus solely on  \iets\ without flips define these maps over  the entire semi-open interval, however we consider that \iets\ are  not defined in the discontinuities.
We denote an \iet\ with $k$-flips and $n-1$ discontinuities as an  \emph{(n,k)-\iet}. If  $T$ has 
no flips it is referred as an  $n$-\iet.  
The set of discontinuity poits of $T$ is denoted by $\Disc(T)=\{a_i\}_{i=2}^n$.  When 
no false discontinuities exist,  
 $T$ is called  a \emph{proper} \iet. A discontinuity $a\in (0,\ell)$ is
 considered \textcolor{black}{false} if $T$ can be continuously extended to  $a$.

The orbit of  $x\in (0,\ell)$ is given by
 $$\O_T(x):=\{T^m(x): m \in \mathbb{Z} \textrm{ and } T^{m}(x) \textrm{ makes sense}\},$$
 where $T^0(x)=x$ for any $x\in[0,\ell]$ and for any $m>0$, $T^{m+1}=T\circ T^{m}$ and
 $T^{-m-1}=T^{-1}\circ T^{-m}$.
 
Furthermore, we define the orbits  $\O_T(0):=\{0\}\cup\O_T(\lim_{x\to 0^+}T(x))$ and  $\O_T(\ell):=\{\ell\}\cup\O_T(\lim_{x\to \ell^-}T(x))$. 
 An \iet\  $T$ is said to be \emph{minimal} if $\O_T(x)$ is dense in $[0,\ell]$
for any  $x\in [0,\ell]$. A weaker condition, 
transitivity, states that  $T$ is 
\emph{transitive} if there exists some $x\in [0,\ell]$ such that
$\O_T(x)$ is dense in $[0,\ell]$. An $n$-\iet\ is said to be \emph{uniquely ergodic} if it
 admits only  the Lebesgue measure (and its multiples) 
as invariant measures. Recall that a finite 
measure $\mu$ is invariant under $T$ if $\mu (A)=\mu (T^{-1}(A))$ for all Borel sets $A$.
 For \iets, unique ergodicity implies  ergodicity with respect to the   Lebesgue measure
 (cf. \cite[Section II.6, Th. 6.1]{MagneLibro}), meaning that if a  Borel set
$E\subset\intervalog$ satisfies {$T^{-1}(E)=E$}, then    $\mu(E)=0$ or
$\mu(E)=1.$

By replacing  $[0,\ell]$ with $\S^1=[0,\ell]/\equiv$, (where $0 \equiv \ell$), we obtain the notion of \emph{circle exchange transformation of $n$-intervals}, denoted \emph{$(n,k)$-\cet}. Given an $(n,k)$-\iet\ $T:D \to (0,\ell)$  as before, we  denote by $T^c$ the $(n,k)$-\cet\ obtained by identifying $0$ and $\ell$. 
The definitions of   \textcolor{black}{false} discontinuities, minimality, 
transitivity, unique ergodicity and properness for \cets\ follow
directly from those for  \iets.

The relationship between minimality and unique ergodicity has been extensively studied
since the early investigations of \iets. Keane conjectured in \cite{keane} 
that minimality would imply unique ergodicity. However  
Keynes-Newton and Keane soon constructed counterexamples of minimal 
non uniquely ergodic oriented \iets\
 in \textcolor{black}{\cite{keane1} and \cite{keynesnewton}, respectively.} 
 Masur \cite{masur} and Veech \cite{veechannals} later independently proved 
 that almost all  \textcolor{black}{oriented} \iets\ associated to 
 an irreducible  permutation are uniquely ergodic 
 (see Sections \ref{sub1} and \ref{sub2} for notation). 
 In the case of \iets\ with flips, Gutierrez, Lloyd, Medvedev, Pires and  Zhuzhoma 
 constructed transitive uniquely ergodic \iets\  with flips and posed
 as an open problem the existence of  minimal non-uniquely ergodic \iets\ with flips, see \cite{gutierrez4b}. \textcolor{black}{In this direction,} the present authors constructed examples of $(10,k)$-\iets\ with $1\leq k\leq 10$ and two independent invariant measures, see \cite{NoUnicamenteErgodica10}.


Returning to the orientable case, Veech showed  in \cite[Th. 2.12]{veech}  
that the maximum number, $\textcolor{black}{N}$, of invariant normalized measures 
that an oriented $n$-$\iet$ can admit satisfies $N\leq \frac{n}{2}$.  In fact, 
he proved that $N\leq \frac{R}{2},$ where $R$ is the rank of the $n\times n$ translation matrix associated to $T$, see \cite[Sec. 2]{veech}. Furthermore, in \cite[S. IV]{YoccozCurso}, 
the reader can follow an analysis in the following terms: $N$ is bounded by $g$,  the genus of the translation surface  obtained by the  suspension of $T$, and this genus satisfies $g=1+\frac{n-m}{2}$, where $m$ is the number of the so called marked points of the surface,
 see \cite[S. 5-6]{espinperaltasoler} for more details on the suspension surface. 
 \textcolor{black}{ In \cite{YoccozCurso}, for any even natural number $n$, Yoccoz even builds $n$-\iets\ with $\frac{n}2$ invariant measures.}
Recently,  J. S. Athreya and J. Chaika  have shown that  the Hausdorff dimension of the set of non uniquely ergodic minimal \iets\ having as associated permutation $\pi=(4,3,2,1)$ is $\frac52$, see \cite{chaikaathreya}. For transitive \nkiets\ with flips, $n> 4$,  A. Skripchenko and S. Troubetzkoy \cite{troubetzkoy2018} have shown that, for any fixed irreducible  permutation $\pi$, the Hausdorff dimension of the set of transitive non uniquely ergodic \iets\ associated to $\pi$ is bigger than or equal to $n-3$.

In the oriented case, we  know the existence of  $4$-\iets\ with two invariant measures, see \cite{keane1,keynesnewton}.
However:

\begin{mteo}\label{maint1}
 Let $T$ be a transitive $(4,k)$-\iet\ with $k\geq 1$. Then, $T$  is uniquely ergodic. 
 Moreover, there are no transitive $(n,k)$-\iets, $1\leq k\leq  n\leq 3$,  and there are no 
 non uniquely ergodic transitive (orientable) $n$-\iets, $1\leq n\leq 3$.
\end{mteo}

%
%


\section{Folklore results}

\subsection{\label{sub1}Rauzy-Veech-Nogueira maps}

A \emph{signed permutation} is an injective map $$\pi:N_n=\{1,2,\ldots,n\}\to
N_n^\sigma=\{-n,-(n-1),\ldots,-1,1,2,\ldots,n\}$$ such that $|\pi|:N_n\to
N_n$ is bijective, that is, a \emph{standard permutation}. Let $S_n^\sigma$ be the set of signed permutations. \textcolor{black}{A \emph{non
standard permutation} will be a signed permutation $\pi$ whenever $\pi(i)<0$ for
some $i\in N_n$.} We represent
$\pi$  by $(\pi(1),\pi(2),\ldots,\pi(n))\in
(N_n^\sigma)^n.$ A permutations is \emph{irreducible} if $|\pi(\{1,2,\ldots t\})|\not=\{1,2,\ldots,t\}$ for any $1\leq t<  n$. The set of \emph{irreducible permutations} is denoted by $S_n^{\sigma,*}$. The subset $S_n^{\sigma,+}$ of  permutations $\pi$ with $|\pi|(n)\not=n$ will play an important role. Observe that $S_n^{\sigma,*}\subsetneq S_n^{\sigma,+}$.

We define the maps $a$ and $b$ on $\permusetplus$ which helps to study the Rauzy induction.
For a given $x\in\R\menos\{0\}$  the sign of $x$ is denoted by $\sigma(x).$
The \emph{generalized Rauzy maps} were introduced by Nogueira in \cite{nogueira2} 
(cf. also \cite{rauzy}) and are   $
a,b:\permusetplus\rightarrow\permuset,$
where $a(\pi)$ and $b(\pi)$ are the permutations defined  by:
\begin{equation}\label{def:a}
a(\pi)(i)=\left\{
\begin{array}{ll}
\pi(i)                          &\textrm{ if }|\pi(i)|\leq |\pi(n)| - \frac{1-\sigma(\pi(n))}{2},\\
\sigma(\pi(n))\sigma(\pi(i))(|\pi(n)|+\frac{1+\sigma(\pi(n))}{2})        &\textrm{ if }|\pi(i)|=n,\\
\sigma(\pi(i))(|\pi(i)|+1)        &\textrm{ otherwise,}\\
\end{array}\right.
\end{equation}
and 
\begin{equation}\label{def:b}
b(\pi)(i)=\left\{
\begin{array}{ll}
\pi(i)                          &\textrm{ if } i \leq |\pi|^{-1}(n) + 
\frac{\sigma(\pi(|\pi|^{-1}(n)))-1}{2},\\
\sigma(\pi(|\pi|^{-1}(n)))\pi(n)        &\textrm{ if } i = |\pi|^{-1}(n) + 
\frac{\sigma(\pi(|\pi|^{-1}(n)))-1}{2} +1, \\
\pi(i-1)        &\textrm{ otherwise.}\\
\end{array}\right.
\end{equation}

Together with these maps, we also define the \emph{generalized
Rauzy matrices} associated to a permutation $\pi\in\permusetplus$, $M_a(\pi)$ and
$M_b(\pi)$. Given  $1\leq i,j\leq n$, $E_{i,j}$ denotes the $n\times n$
matrix having zeros in all the positions except for the position $(i,j)$
which is equal to 1, and $\mathrm{Id}_n$ denotes the $n\times n$ identity matrix. The
definitions of $M_a(\pi)$ and $M_b(\pi)$ are:
\begin{eqnarray}
M_a(\pi)&=&\mathrm{Id}_n+E_{n,|\pi|^{-1}(n)},\nonumber\\
M_b(\pi)&=&\left(\sum\limits_{i=1}^{|\pi|^{-1}(n)}E_{i,i}\right)
+ E_{n, |\pi|^{-1}(n)+ \frac{1+\sigma(\pi(|\pi|^{-1}(n)))}{2}} 
+ \left(\sum\limits_{i=|\pi|^{-1}(n)}^{n-1}E_{i,i+1}\right).
\label{E:matrizAB}
\end{eqnarray}

In the following, a matrix $A$ is said to be \emph{nonnegative} if all the entries are in $[0,+\infty)$, \textcolor{black}{shortly we write $A\geq 0$}; and it is said to be \emph{positive} when all the entries are strictly positive, $A>0$.

\subsection{Coordinates}\label{sub2}

From now on, we write $\R_+=(0,\infty)$,  $\Lambda^n=\R_+^n=\{\sum_{{{i=1}}}^{{n}}\lambda_ie_i:\lambda_i>0\}$, 
 where $e_i$ denotes the $i$-th vector of the canonical basis of $\R^{{n}}$.  We establish a one-to-one map between the set of $n$-\iets\  and $\codes=\Lambda^n\times
S_n^\sigma$.
 Let $T:D=\bigcup_{i=1}^n I_i\subset(0,\ell)\to(0,\ell)$  be an $n$-\iet. Then its associated coordinates in $(\lambda,\pi)\in \codes$ are  defined
by:
\begin{itemize}
\item \label{notacion}$\lambda_i=a_{i+1}-a_i$ for all $i\in N_n$.
\item $\pi(i)$ is positive (resp. negative) if $T|_{I_i}$ has
slope $1$ (resp. $-1$). Moreover $|\pi(i)|$ is the position of the
interval $T(I_i)$ in the set $\{T(I_i)\}_{i=1}^n$ with the usual order in $\R$.
\end{itemize}

These coordinates allow us to identify $T=(\lambda,\pi)$. For a fixed permutation $\pi,$ we can consider
the Lebesgue measure of the cone $\Lambda^n$ on the set of  $n$-\iets\ having associated permutation $\pi$.

\subsection{Rauzy induction}

Roughly speaking, the \emph{generalized Rauzy induction} is an
operator in the set of \iets\  which sends any $T:D\subset \intervalog\to \intervalog$
to its first return map on some subinterval
$\intervalogp\subsetneq\intervalog.$

We are now ready to present formally the \emph{generalized Rauzy 
operator} $R$, we follow \cite{linerosolerergodic} (see also \cite{nogueira2}). Let 
\begin{equation}\label{E:DominioRauzy}
 \mathcal{D}=\{(\lambda,\pi)\in\Lambda^n\times
S_n^\sigma: \lambda_{n}\not=\lambda_{|\pi|^{-1}(n)}\},
\end{equation}
 then
$$
\begin{array}{rcl}
R:\mathcal{D}\subset \Lambda^n\times S_n^\sigma&\longrightarrow&
\Lambda^n\times S_n^\sigma\\
T=(\lambda,\pi)&\mapsto&T'=(\lambda',\pi')
\end{array}$$
is defined by:

$$T'=(\lambda',\pi')=\left\{
\begin{array}{ll}
\left(M_a(\pi)^{-1}\lambda,a(\pi)\right)&\textrm{if }
\lambda_{|\pi|^{-1}(n)}<\lambda_n,\\
\left(M_b(\pi)^{-1}\lambda,b(\pi)\right)&\textrm{if }
\lambda_{|\pi|^{-1}(n)}>\lambda_n.\\
\end{array}\right.$$

If $T'$ is obtained from $T$ by means of the operator $a,$ $T$ is said to be of \emph{type} $a$, otherwise $T$ is of \emph{type} $b$. In any of the two cases, the largest of the two intervals is called the \emph{winner} and the shortest one is called the \emph{loser}  and  $T'$  is the \emph{Poincar\'e first return map} induced by $T$ on $[0,\ell'],$ with $\ell'=\ell-\min\{\lambda_n,\lambda_{\antin}\}$.

It is interesting to investigate how the intervals are exchanged when we apply the Rauzy operator.

\begin{claim}\label{Cl:ExchangeIntervalsBya}\label{Cl:ExchangeIntervalsByb}
Assume that $T=(\lambda,\pi)$ is of type $a$ and we have the length of the domain intervals placed as $(\lambda_1,\lambda_2,\dots ,\lambda_n)$; then, the length of the domain intervals of $R(T)$ will be in the order  $(\lambda_1,\lambda_2,\dots, \lambda_{n-1},\lambda_n')$, where the length of the last interval is $\lambda_n'=\lambda_n-\lambda_{|\pi|^{-1}(n)}$.

Let $T=(\lambda,\pi)$ be of type $b$ and suppose that the lengths of the  domain intervals are placed as $(\lambda_1,\lambda_2,\dots ,\lambda_n)$; then, the length of the  domain intervals of $R(T)$ will be placed as 
$$(\lambda_1,\lambda_2,\dots ,\lambda_{|\pi|^{-1}(n)-1},\lambda_n,\lambda'_{|\pi|^{-1}(n)},\lambda_{|\pi|^-1(n)+1},\dots, \lambda_{n-1})\quad\textrm{ if }\quad \sigma(\pi(|\pi|^{-1}(n)))<0$$ 
or
$$(\lambda_1,\lambda_2,\dots \lambda_{|\pi|^{-1}(n)-1}, \lambda_{|\pi|^{-1}(n)}',\lambda_n,\lambda_{|\pi|^{-1}(n)+1},\dots, \lambda_{n-1})\quad\textrm{ if }\quad\sigma(\pi(|\pi|^{-1}(n)))>0,$$ 
where  $\lambda_{|\pi|^{-1}(n)}'=\lambda_{|\pi|^{-1}(n)}-\lambda_n$. $\hfill\Box$
\end{claim}

\subsection{Rauzy graphs}

The operators $a$ and $b$ induce in the set $S_n^{\sigma,*}$ \textcolor{black}{of irreducible permutations} a directed graph
structure whose vertices are all the points from $S_n^{\sigma,*}$ and the
directed edges are arrows labelled by $a$ and $b$. Given $\pi,\pi'\in
S_n^{\sigma,*},$ there is  an arrow labelled by $a$ (resp. $b$) from $\pi$ to
$\pi'$ if and only if $v(\pi)=\pi'$ ($v\in\{a,b\}$). In this situation, $\pi'$ is said to be a \emph{successor} of $\pi$ and $\pi$ is a \emph{predecessor} of $\pi'$.  This graph is called the \emph{Rauzy graph}.

We remark that we only take into account irreducible permutations because they are the only ones for which the associated
\iets\ can be minimal. Moreover, it is worth noticing that if $\pi$ is an irreducible standard permutation then 
$a(\pi)$ and $b(\pi)$ are irreducible, while it is not always the case for non standard irreducible permutations, 
observe for instance that $a(-4,3,2,-1)=(1,4,3,-2)$.  When $v(\pi)$ is reducible,  the vertex $\pi$ has not an outgoing arrow labelled by $v$; however, for convenience reasons, we can stress this situation in the Rauzy graph by drawing an arrow labelled by $v$, starting in $\pi$ and ending in a so called \emph{hole} (which identifies all the reducible permutations). Also, sometimes it may be convenient to remark when a permutation has no predecessor by $v\in\{a,b\}$ by adding an arrow incoming to $\pi$ but having no vertex in the beginning. Observe the different ways of representing Rauzy subgraphs in Figure~\ref{clase19A} and Figure~\ref{nonrepellor}. The first one has richer information while the second is a proper Rauzy subgraph.


A \emph{path (of operators)} is an element of $\{a,b\}^L$, where $L\in\N$
or $L=\N$. An easy way of constructing Rauzy subgraphs
is to take a vertex $\pi\in S_n^{\sigma,*}$ and to construct recursively other
vertices by applying a path of operators. The \emph{Rauzy subgraph}, $\mathcal{G}^{\pi_1,v}$,
associated to $\pi_1\in S_n^{\sigma,*}$ and $v=(v_1,v_2,\dots)\in \{a,b\}^L$, 
 is the graph with vertices $\{\pi_i\}_{i=1}^L$
satisfying $v_i(\pi_i)=\pi_{i+1}$, $1\leq i < L$, the edges of this graph being arrows
labelled by $v_i$ from $\pi_i$ to $\pi_{i+1}$, \textcolor{black}{that is, $\pi_{i} \rTo^{v_i} \pi_{i+1}$.} Observe that any $n$-\iet,
$T=(\lambda,\pi)\in\mathcal{D}$, defines a Rauzy subgraph in a natural way, the one
associated to $\pi$ and the path of operators $v$ defined by the Rauzy 
induction, that is, $v_i$ is $a$ (resp $b$) if $R^{i-1}(T)$ is of type $a$ (resp. $b$); we denote this subgraph by $\mathcal{G}^T$. We will say that $T$
is \emph{infinitely inducible} if $v$ has infinite dimension, i.e.,
$v\in\{a,b\}^{\N}$. For a finite path of operators, $v\in\{a,b\}^k${, $v=(v_1,\ldots,v_k)$}, $k\in\N$, and a vertex
$\pi\in S_n^{\sigma,*}$, $v(\pi)$ denotes the vertex obtained after applying sequentially, from the left, the operators in $v$; also, 
\begin{equation}M_v(\pi):=
M_{v_1}(\pi)
M_{v_2}(v_1(\pi))
M_{v_3}(v_2(v_1(\pi)))
\dots
M_{v_k}(v_{k-1}(v_{k-2}(\dots v_1(\pi))))\label{E:matrizv}
\end{equation}
\textcolor{black}{is the matrix constructed by the product of subsequent generalized Rauzy matrices.}

\subsection{Rauzy classes}

Operators $a$ and $b$ also induce in  $S_n^{\sigma,*}$  \textcolor{black}{a relation} that divides $S_n^{\sigma,*}$ into  \textcolor{black}{different} classes, so called \emph{Rauzy classes}. Two elements $\pi,\pi'\in S^{\sigma,*}_n$ are said to be related \textcolor{black}{(or \textit{connected})} if and only if one of the following conditions hold:
\begin{itemize}
 \item $\pi=\pi'$.
 \item There exists two \emph{Rauzy subgraphs}, $\mathcal{G}^{\pi_1,v_1}$ and $\mathcal{G}^{\pi_2,v_2}$
in such a way that $\pi$  is a vertex of $\mathcal{G}^{\pi_1,v_1}$,   $\pi'$  is a vertex
of $\mathcal{G}^{\pi_2,v_2}$  and both subgraphs have at least a common vertex.
\end{itemize}
Observe that the elements of any Rauzy class 
\textcolor{black}{
   are the vertices of a maximal connected Rauzy subgraph, where here the connection is understood at least in a direction, but not necessarily in a bidirectional way; for instance, in class $\mathcal{C}_5$, see Figure~\ref{clase5}, the elements $\pi_{5,1}$ and $\pi_{5,3}$ belong to the same maximal class, but it is not possible to access to $\pi_{5,3}$ from $\pi_{5,1}$, although at least we can have a path from $\pi_{5,3}$ to $\pi_{5,1}$.  
}
 This Rauzy subgraph may have no arrows in classes with one element.

The orientable Rauzy-classes have already been  studied in \cite{kontsevichzorich}. However, this has not been the case in the non oriented case, in which a  lot of classes appear; in particular, for four intervals, $n=4$, we find twenty Rauzy classes that will be introduced \textcolor{black}{later},  only two of them  \textcolor{black}{being} orientable. Certainly, some of them have few elements, for example, the following elements form individual classes:
\begin{eqnarray*}
\imt{cl4[1][1]},\,
\imt{cl4[2][1]},\,
\imt{cl4[3][1]},\,
\imt{cl4[4][1]},\,
\imt{cl4[6][1]},\\
\imt{cl4[7][1]},\,
\imt{cl4[9][1]},\,
\imt{cl4[10][1]}.
\end{eqnarray*}
The cardinality of the remaining classes are 
$4,6,7,8,27,29,36$ and $59$. All of them will be introduced later.

\subsection{Inverse and Conjugate}

\newcommand{\Id}{\mathrm{Id}}
Observe that if $T=(\lambda,\pi)$ is a \nkiet\ then its inverse is also a \nkiet, and the minimality of one of them implies the minimality of the other one since the orbits are taken in both directions of time. Also, take into consideration the homeomorphism $h:[0,\ell]\to [0,\ell]$, $h(x)=\ell-x$, which holds $h^2=\Id$, and  define  $T^h:=h\circ T\circ h^{-1}$, the  conjugated map of  $T$ by means of $h$; then, $T^h$ is also an \iet\ (since $h$ is an isometry). \textcolor{black}{In fact, by \cite[Section 1.3]{calisoPeriodosIets1} we  have, among others, the following properties.}

\begin{Theorem}\label{Th:hlamda} Let $T=(\lambda,\pi)$ be a \nkiet. Then:
 \begin{enumerate}
 \item $T^{-1}=(\lambda^i,\pi^i)$ with $\lambda^i_j=\lambda_{|\pi|^{-1}(j)}$ and $\pi^i(j)=\sigma(\pi(|\pi^{-1}|(j))) |\pi|^{-1}(j)$, 
\textcolor{black}{$1\leq j \leq n$.}

\item
$T^h=(\lambda^h,\pi^h)$ with $\lambda^h_j=\lambda_{n+1-j}$ and
$\pi^h(j)=\sigma(\pi(n+1-j)) (n+1-|\pi|(n+1-j))$, $1\leq j \leq n.$

\item $T^h=(\lambda^h,\pi^h)$ is minimal if and only if $T=(\lambda, \pi)$ is minimal.
 \item $T^{-1}=(\lambda^i, \pi^i)$ is minimal if and only if $T=(\lambda, \pi)$ is minimal.
\end{enumerate}
\end{Theorem}
%

\subsection{Minimality and invariant measures}

The concepts of minimality, Keane condition and infinite induction under the Rauzy operator are related for \iets. In this section we  give the connections among them.

For a given  \nkiet\ $T$ with domain $D=\cup_{i=1}^n (a_i,a_{i+1})$, we  write $T(a_i^\+):=\lim_{x\to a_i^+}T(x)$ for $1\leq i\leq n$  and
$T(a_i^\-):=\lim_{x\to a_i^-}T(x)$ for $2\leq i\leq n+1$. We also write $T(a_1^\-)=T(a_1^\+)$ and $T(a_{n+1}^\+)=T(a_{n+1}^\-)$.  Then, we say that $T$ satisfies the \emph{Keane condition} if and only if:
$$T^m(a)\not=a_j \textrm{ for all }m\geq 1,\; 2\leq j\leq n \textrm{ and } a\in\bigcup_{i=1}^{n+1}\{a_i^\+,a_i^\-\}.$$

If a \nkiet\ is minimal then it  satisfies Keane condition, otherwise it would have finite orbits. In fact, we have the following  equivalence (see \cite{MR3170630}):
\begin{theorem}\label{Th:KeaneConditionMinimal}
 $T:D=\cup_{i=1}^n(a_i,a_{i+1})\subset[0,\ell]\to[0,\ell]$ satisfies the Keane condition if and only if it is minimal.
\end{theorem}

We claim that for oriented \iets, if  defined in the whole $[0,\ell)$, minimality does not necesarily implies Keane condition, see \cite[Remark 4.5]{viana}.
With the help of the previous theorem, we have:
\begin{proposition} \label{P:minimalidadimplicainducible}
Let $T=(\lambda,\pi)$ be minimal. Then, $T$ has an infinite orbit under the Rauzy operator; and for any $m\in \mathbb{N}$, $\pi_m\in \irreducibleset$ 
whenever  $R^m(T)=(\lambda^m,\pi_m)$.

Moreover, $|\lambda^m|=\sum_{j=1}^n\lambda^m_j$ goes to 0 as $m$ goes to $+\infty$; and if $v_m$ is the type of $R^m(T)$, then $(v_m)_{m>m_0}$ is not constant for any $m_0\in\N$.

Conversely, let $T=(\lambda,\pi)$ such that  $R^{(m)}(T)=(\lambda^m,\pi_m)$ for all $m>0$ and $\pi_m$ is irreducible. Then $T$ is minimal.
\end{proposition}
\begin{proof}
If  $T=(\lambda,\pi):\bigcup_{i=1}^n(a_i,a_{i+1})\subset[0,\ell]\to[0,\ell]$ is minimal, then $(\lambda,\pi) \in\mathcal{D}$, see (\ref{E:DominioRauzy}), otherwise
$\lambda_n=\lambda_\antin$ and there exist $j\in\{1,2,\dots,n\}$ and  $\X\in\{\+,\-\}$ such that
$T(a_j^{\X})=a_{n}$, which would imply, \textcolor{black}{according to Theorem~\ref{Th:KeaneConditionMinimal}, that $T$ is not minimal.} 

Observe that if $T$ is of type $a$ then $\Disc(R(T))=\Disc (T)$, and if $T$ is of type $b$ then
$\Disc(R(T))=\Disc(T)\menos\{a_{n}\}\cup\{T^{-1}(a_n)\}$; also, the orbits of
$R(T)$ are subsets of those of  $T$. Observe that  $R(T)$ satisfies Keane condition because it holds for $T$. Then, $R(T)=(\lambda^1,\pi_1)$ is minimal and
$\pi_1\in\irreducibleset$.

By recurrence, we obtain that for any $m\in\N$,
$R^m(T)=(\lambda^m,\pi_m)$ is minimal and $\pi_m\in\irreducibleset$.

The second paragraph of the proposition is proved in \cite[Lemma 6 and Proposition 10]{MR3170630} and the  third one  in \cite[Th. A]{MR3170630}.
\end{proof}

Rauzy subgraphs and the cone $\mathcal{M}(T)$ of invariant measures associated to $T$  are closely related.  In some sense, next result is a converse of Proposition~\ref{P:minimalidadimplicainducible}.

\newcommand{\matrizpaso}[1]{M_{v_{#1}}(\pi_{#1})}

\begin{theorem}\label{T-cone}
Let $T=(\lambda^1,\pi_1)$ be an $n$-\iet\ such that: $\pi_1\in S_n^{\sigma,*}$;
$T$ is infinitely inducible;  and
$R^i(\lambda^1,\pi_1)=(\lambda^{i+1},\pi_{i+1})$ for any $i\geq 1$.  Let $\mathcal{G}^T$ be the  Rauzy subgraph of
$T$ associated to $\pi_1$ and
$v\in\{a,b\}^\N$.  Put
\begin{equation}
\mathcal{C}(\mathcal{G}^T):=\bigcap_{i=1}^\infty M_{v_1}(\pi_1)\cdot
M_{v_2}(\pi_2)\cdot \ldots \cdot M_{v_i}(\pi_i) \Lambda^n,
\label{E:cone}
\end{equation}
and let $\mathcal{M}(T)$ be  the cone of invariant measures associated to $T$.
Assume also that, for any $i\in\N$, $\pi_i$ is irreducible. Then:
\begin{enumerate}
\item $\lambda^1\in\mathcal{C}(\mathcal{G}^T)$;
\item if $\gamma\in \mathcal{C}(\mathcal{G}^T)$ and $S=(\gamma,\pi_1)$, \quad
$\mathcal{G}^S=\mathcal{G}^T$;
\item $\mathcal{C}(\mathcal{G}^T)$ and  $\mathcal{M}(T)$ are linearly isomorphic (thus,  if $\mathcal{C}(\mathcal{G}^T)$ is a half line, $T$ is uniquely ergodic);
\item $R^j(T)$ is minimal for any $j\in\N\cup\{0\}$.
\end{enumerate}
\end{theorem}

\begin{proof}
 See \cite[Th. 2]{NoUnicamenteErgodica10}.
\end{proof}

\newcommand{\conografo}[1]{\mathcal{C}(\mathcal{G}^{#1})}

\newcommand{\ensayo}[1]{#1
}

The following Theorem provides a method for constructing minimal \iets\ by means of Rauzy
graphs.

\begin{theorem}\label{conver}\label{T:xxlnoperiodico} Let
$\mathcal{G}^{\pi_1,v}$ be a  the Rauzy-subgraph associated  to
$\pi_1\in S_n^{\sigma,*}$ and $v\in\{a,b\}^\N$. Assume the existence of a   sequence  $(n_k)_k${,} $n_k\in\N$, satisfying
$n_1=1$, $n_{k+1}>n_k$  and  {such that}
$M_{v_{n_k}}(\pi_{n_k})\cdot
M_{v_{n_k+1}}(\pi_{n_k+1})\cdot \ldots \cdot M_{v_{n_{k+1}-1}}(\pi_{v_{n_{k+1}-1}})$ {is positive} for any $k\in\N$. {If} $\mathcal{C}(\mathcal{G}^{\pi_1,v})=\bigcap_{i=1}^\infty M_{v_1}(\pi_1)\cdot
M_{v_2}(\pi_2)\cdot \ldots \cdot M_{v_i}(\pi_i) \Lambda^n$, then:
\begin{enumerate}
\item  $\mathcal{C}(\mathcal{G}^{\pi_1,v})$ is nonempty.
 \item There exists $\lambda^1\in\mathcal{C}(\mathcal{G}^{\textcolor{black}{\pi_1},v})$ such that the associated graph to $T=(\lambda^1,\pi_1)$ is
$\mathcal{G}^{\pi_1,v}$.
\item $R^j(T)$ is minimal for any $j\in\N\cup\{0\}$.
\item $T$ is uniquely ergodic if and only if $R^j(T)$ is  uniquely ergodic for any $j\in\N\cup\{0\}$.
\end{enumerate}
\end{theorem}

\begin{proof}
 See \cite[Theorem 3]{NoUnicamenteErgodica10}.
\end{proof}

To end this section, we introduce \textcolor{black}{a result on unique ergodicity that will be fruitful in the sequel.}

\begin{lemma}\label{C-rutafinita}
 Let $v\in\{a,b\}^L$ for some $L\in\N$ and let   $\pi_1\in S_n^{\sigma,*}$ be such that
 $\mathcal{G}^{\pi_1,v}$ has all its vertices in $S_n^{\sigma,*}$. Then there exists an IET $T=(\lambda,\pi_1)$ such that $\mathcal{G}^T=\mathcal{G}^{\pi_1,v}.$ Moreover, given $\gamma\in\Lambda^n$, there exists an IET $T=(\lambda,\pi_1)$, $\lambda=M_{v_1}(\pi_1)M_{v_2}(\pi_2)\dots M_{v_{L}}(\pi_L)\gamma$, such
 that $R^L(T)=S=(\gamma,\pi_{L+1})$. 
 
 In fact, if  $S$ satisfy the hypothesis of Theorem~\ref{T-cone} (and thus also $T$)  then 
 \begin{equation}
 \mathcal{C}(\mathcal{G}^T)=M_{v_1}(\pi_1)M_{v_2}(\pi_2)\dots M_{v_{L}}(\pi_L)\mathcal{C}(\mathcal{G}^S)
  \label{E:conoTS}
 \end{equation}
 and $S$ is uniquely ergodic if and only if  $T$ is uniquely ergodic. 
  
\end{lemma}

\begin{proof}
The first paragraph is in  \cite[Corollary 19]{linerosolerergodic}.
Equation (\ref{E:conoTS}) is a direct consequence of \cite[Lemma 27-(4)]{NoUnicamenteErgodica10}
and then, by Theorem~\ref{T-cone}-(3) and Equation (\ref{E:conoTS}), 
either both $T$ and $S$ are uniquely ergodic or non uniquely ergodic.
\end{proof}

\subsection{Simplices and projective metric}

\newcommand{\ddp}{\mathrm{d_p}}
%
\textcolor{black}{In this section, we concentrate in clarifying the action of some linear maps on the cone} $\Lambda^n$; \textcolor{black}{we need these properties as an essential tool to prove the unique ergodicity of minimal} $(4,k)$-IETs. 

A set
$K \subset \Lambda^n$ is said to be \emph{ relatively compact} in $\Lambda^n$ if and only if the infimum  of the set of coordinates of the vectors in $K$ is positive, that is,
$$\textcolor{black}{\inf_{u\in K} \{p_j (u): 1\leq j\leq n\}>0,}$$
\textcolor{black}{where $p_j$ is the canonical projection on the $j$-coordinate, $1\leq j\leq n$.}

We are going to introduce the notion of projective pseudometric on the cone $\Lambda^n$  given by Birkhoff \cite{birkhoffmetrica}; we follow here \cite[Sec. 26]{viana} (see also \cite[p. 240]{veech} and \cite[Ch. 4.15]{furstenberglibro}). Let $u,v\in\Lambda^n$, define $|u|=\sum_{i=1}^nu_i$ and write:
$$
a(u,v)\textcolor{black}{:=}\inf\left\{\frac{v_i}{u_i}\right\}_{1\leq i\leq n}\qquad
b(u,v)\textcolor{black}{:=}\sup\left\{\frac{v_i}{u_i}\right\}_{1\leq i\leq n}.$$
Then, it is straightforward to prove:

\begin{lemma}\label{L:a(u,v)b(u,v)} Under the above conditions, 

\begin{equation}\label{Eq:a(u,v)}
v-\lambda u\in\Lambda^n \textrm{ if and only if } \lambda<a(u,v), 
\end{equation}

\begin{equation}\label{Eq:b(u,v)}
\mu u-v\in\Lambda^n \textrm{ if and only if } \mu>b(u,v). 
\end{equation}
with $\lambda,\mu\in\mathbb R$. $\hfill\Box$

\end{lemma}

Now, the \emph{projective pseudometric} in $\Lambda^n$ is:
$$\ddp(u,v)\textcolor{black}{:=}\log\frac{b(u,v)}{a(u,v)}=
\log\sup\left\{\frac{u_i}{v_i}\frac{v_j}{u_j}:1\leq i,j\leq n\right\}.$$

This map has the following properties for any  $u,v,w\in\Lambda^n$:
 $\ddp(u,v)=\ddp(v,u)$;  $\ddp(u,w)\leq \ddp(u,v)+\ddp(v,w);$ $\ddp(u,v)\geq 0$ and  $\ddp(u,v)=0$ if and only if $u=tv$ for some $t>0$. When we restrict $\ddp$ to the simplex $$\Delta^n=\{v\in\Lambda^n:|v|=1\},$$ it becomes a metric.  
 Observe that if   $K\subset\Delta^n $ is relatively compact, 
  then its  diameter, $\diam(K)$, relative
to the projective metric is finite, that is,
$\diam(K)=\sup_{u,v\in K}\{\ddp (u, v)\}<\infty.$
\textcolor{black}{Notice also  that, in general, a relatively compact set $K$ in $\Lambda^n$ does not have necessarily a finite diameter. For instance, consider the set $K=\{(n,1,\ldots,1):n\in\mathbb N\}$, with $\ddp\left((1,1,\ldots,1),(n,1,\ldots,1)\right)=\log n.$}

 From now on, we deal with linear maps $F:\R^n\to\R^n$  satisfying
 \begin{equation}
 F(\Lambda^n)\textcolor{black}{\subseteq}\Lambda^n. \label{eqProyective}
 \end{equation}
Observe that this condition is equivalent to say 
that the associated matrix to $F$ (in the canonical basis) has no zero rows and all the entries are nonnegative.


\textcolor{black}{The following result, whose proof can be consulted in \cite[Section 26, p.83]{viana}, establishes that such $F$ is not expansive with respect to $\ddp$. }
\begin{proposition} \label{Lipschitz1}
Let $F:\R^n\to\R^n$ be a linear map such that $F(\Lambda^n)\textcolor{black}{\subseteq}\Lambda^n$. Then:
$$\ddp(F(u),F(v))\leq \ddp(u,v)\quad \textrm{for all } u,v\in\Lambda^n.$$ 
\end{proposition}

The projectivization of a linear map $F:\R^n\to \R^n$ satisfying (\ref{eqProyective}) is the map
$$f:\Delta^n\to\Delta^n, \textrm{ with } f(u)=\frac{F(u)}{|F(u)|},$$
\textcolor{black}{where $|F(u)|=\sum_{1\leq j\leq n} |F_j(u)|$ if $F=(F_1,\ldots,F_n)$}.
The following proposition can be seen in \cite[Proposition 26.3]{viana}.

\begin{proposition}  \label{Lipschitz2}
Let $F$ be a linear map satisfying  (\ref{eqProyective}), and let $(F_{i,j})_{i,j}$ be the associated matrix to $F$. Assume that there exists $\delta>0$ with 
$\delta<F_{i,j}<\frac1{\delta}$ for any $1\leq i,j\leq n$. Then: 
$$\ddp(F(u),F((v))\leq\frac{1-\delta^2}{1+\delta^2}\ddp(u,v).$$
\end{proposition}

We prove here a generalization of the previous result.

\begin{proposition} \label{Lipschitz3}
Let $F$ be a linear map satisfying  (\ref{eqProyective}), and let $(F_{i,j})_{i,j}$ be the associated matrix to $F$. Assume that
there exist $\delta_1^i$ and $\delta_2^i$ satisfying $0<\delta_1^i<F_{i,j}<\delta_2^i$  for any $1\leq i,j\leq n$. Then:
$$\ddp(F(u),F((v))\leq\frac{\delta-\delta^{-1}}{\delta+\delta^{-1}}\ddp(u,v), $$
with
$$\delta=\sqrt{\sup_{1\leq i\leq n}\left\{\frac{\delta_2^i}{\delta_1^i}\right\}}.$$
\end{proposition}
\begin{proof}
The proof follows the ideas of \cite[Proposition 26.3]{viana}. Let $\Delta_1= \inf_i\frac{\delta_1^i}{\delta_2^i} \textrm{ and }
\Delta_2= \sup_i\frac{\delta_2^i}{\delta_1^i}.$ \textcolor{black}{Notice that $\Delta_1\cdot\Delta_2=1$, with $\Delta_1\leq 1 \leq \Delta_2.$}

For any $z,w\in\Lambda^n$ write $c(z,w)=\frac{|w|}{|z|}$ and observe that
\begin{equation}
 a(F(z),F(w))
=
\inf_i\frac{\sum_jF_{i,j}w_j}{\sum_jF_{i,j}z_j}
>
\inf_i\frac{\delta_1^i\sum_j w_j}{\delta_2^i\sum_jz_j}
=\Delta_1\frac{|w|}{|z|}.
\label{Eq-a}
\end{equation}
Similarly
\begin{equation}
b(F(z),F(w))=\sup_i\frac{\sum_jF_{i,j}w_j}{\sum_jF_{i,j}z_j}
<
\sup_i\frac{\delta_2^i\sum_j w_j}{\delta_1^i\sum_j z_j}=
\Delta_2 \frac{|w|}{|z|}.\label{Eq-b}
\end{equation}

Given $u,v\in\Lambda^n,$ choose arbitrary sequences $(t_n)$ and $(s_n)$ satisfying $0<t_n<a(u,v)\leq b(u,v)<s_n$, thus $z_n=v-t_n u$ and $w_n=s_n u-v$ are, by Lemma~\ref{L:a(u,v)b(u,v)}, in $\Lambda^n$. Now, for these values $z_n, w_n$, by
(\ref{Eq-a}-\ref{Eq-b}) \textcolor{black}{and Lemma~\ref{L:a(u,v)b(u,v)}} we have:
$$F(s_n u-v)-\frac{|s_n u-v|}{|v-t_nu|}\Delta_1F(v-t_n u)\in\Lambda^n
\quad \textrm{ and }\quad
\frac{|s_n u-v|}{|v-t_nu|}\Delta_2 F(v-t_n u)-F(s_n u-v)\in\Lambda^n.
$$
 By taking $T_n:=\frac{|s_n u-v|}{|v-t_nu|}\Delta_1$ and $S_n:=\frac{|s_n u-v|}{|v-t_nu|}\Delta_2$, \textcolor{black}{and using the linearity of $F$}, we rewrite the above expressions as:
 $$(s_n+t_nT_n)F(u)-(1+T_n)F(v)\in\Lambda^n\quad\textrm{ and }\quad
 (1+S_n)F(v)-(s_n+t_nS_n)F(u)\in\Lambda^n.$$
 Then, \textcolor{black}{again by Lemma~\ref{L:a(u,v)b(u,v)}}:
 $$b(F(u),F(v))<\frac{s_n+t_nT_n}{1+T_n}\quad\textrm{ and }\quad
 a(F(u),F(v))>\frac{s_n+t_nS_n}{1+S_n}.$$
 Therefore:
 \begin{eqnarray*}
\ddp(F(u),F(v)) &<&
 \log\left(\frac{s_n+t_nT_n}{1+T_n}\,\frac{1+S_n}{s_n+t_nS_n}\right) \\
 &=& \log\left(\frac{s_n/t_n+T_n}{1+T_n}\,\frac{1+S_n}{s_n/t_n+S_n}\right)
\\
&=&\log\left(\frac{s_n}{t_n}+T_n\right)
-\log(1+T_n)-\log\left(\frac{s_n}{t_n}+S_n\right)
+\log(1+S_n)\\&=&
\int_0^{\log(s_n/t_n)}\left(\frac{e^x}{e^x+T_n}-\frac{e^x}{e^x+S_n}\right)dx
\leq\sup_{x>0}\frac{e^x(S_n-T_n)}{(e^x+T_n)(e^x+S_n)}\log\left(\frac{s_n}{t_n}\right).
 \end{eqnarray*}
\textcolor{black}{It only remains to bound $\sup_{x>0}\frac{e^x(S_n-T_n)}{(e^x+T_n)(e^x+S_n)}$.  Since the real valued map $\varphi(x):=\frac{e^x(S_n-T_n)}{(e^x+T_n)(e^x+S_n)}$ attains its maximum (positive) value at $x_0=\frac{1}{2}\log(S_nT_n)$ and taking into account that   $\frac{S_n}{T_n}=\frac{\Delta_2}{\Delta_1}=\Delta_2^2 >1,$ we find $\varphi(x_0)=\frac{1-\sqrt{\frac{T_n}{S_n}} }  { 1+\sqrt{\frac{T_n}{S_n}} } =  \frac{\sqrt{\Delta_2}-\sqrt{\Delta_2^{-1}}}{\sqrt{\Delta_2}+\sqrt{\Delta_2^{-1}}}=\frac{\delta-\delta^{-1}}{\delta+\delta^{-1}}.$ Thus,  
$\sup_{x>0}\frac{e^x(S_n-T_n)}{(e^x+T_n)(e^x+S_n)}\leq \frac{\delta-\delta^{-1}}{\delta+\delta^{-1}}$ and, consequently, 
$\ddp(F(u),F(v))\leq \frac{\delta-\delta^{-1}}{\delta+\delta^{-1}}\log\left(\frac{s_n}{t_n}\right).$ Finally, since $\delta$ is a fixed number, making $s_n\searrow b(u,v)$ and $t_n\nearrow a(u,v)$, by taking limits when $n$ goes to infinity, we conclude that 
$\ddp(F(u),F(v))\leq \frac{\delta-\delta^{-1}}{\delta+\delta^{-1}} \log\left(\frac{b(u,v)}{a(u,v)}\right)$, i.e., 
\[\ddp(F(u),F(v))\leq \frac{\delta-\delta^{-1}}{\delta+\delta^{-1}} \ddp(u,v).\] 
}
\end{proof}
\CuentasMaxima{
clasetamagno3:[8,[[[-3,-2,-1]],[[-3,-2,1],[-3,1,2],[-2,-3,1]],[[-3,2,-1]],[[-3,2,1],[-2,3,1]],[[-3,-1,-2],[-3,-1,2],[-3,1,-2],[-2,-3,-1],[-2,3,-1],[2,-3,-1],[2,-3,1],[3,-2,1],[3,1,-2]],[[2,3,-1],[3,-2,-1],[3,-1,-2]],[[2,3,1],[3,1,2],[3,2,1]],[[3,-1,2],[3,2,-1]]]][2]
}\CuentasMaxima{
clase1tamagno3:clasetamagno3[1],
clase2tamagno3:clasetamagno3[2],
clase3tamagno3:clasetamagno3[3],
clase4tamagno3:clasetamagno3[4],
clase5tamagno3:clasetamagno3[5],
clase6tamagno3:clasetamagno3[6],
clase7tamagno3:clasetamagno3[7],
clase8tamagno3:clasetamagno3[8]
}
\section{Rauzy classes for $n=2,3$}
In this section, we need the relationship between \iets\ and \cets. Recall that given a minimal  $(n,k)$-\iet\  
$T:D\subset[0,\ell]\to[0,\ell]$, we obtain a minimal  $(n,k)$-\cet\ denoted by $T^c$ by identifying $0\equiv \ell$; however, if $T$ is proper, $T^c$ could  
not be proper. In particular: 
\begin{claim} \label{R-minlinsol} If $T=(\lambda,\pi):D=\bigcup_{i=1}^n I_i\to [0,\ell]$ is a proper $(n,k)$-\iet\ then $T^c$ is a proper $(n-1,k)$-\cet\ (respectively,  a proper  $(n-1,k-1)$-\cet) only if, after calling $\pi(n+1)=\pi(1)$, we have $\pi(i)=n$ and $\pi(i+1)=1$  (respectively, $\pi(i)=-1$ and $\pi(i+1)=-n$) for some $1
  \leq i \leq n$. 
\end{claim}
Also we are going to need  Main Theorem of \cite{gutierrez4b}:
\begin{theorem} \label{T:transigutierrez} Given 
$n\geq k\geq 1,$ there exists a transitive proper  
$(n,k)$-\cet\ if and only if  $n+k\geq 5$.
\end{theorem}

\subsection{The case $n=2$}

Note that the only standard irreducible permutation for $n=2$ is $(2,1);$ then, the set of irreducible signed permutations is
$$S_2^{\sigma,*}=\{(-2,-1),(-2,1),(2,-1),(2,1)\},$$
and a direct application of Claim~\ref{R-minlinsol} and Theorem~\ref{T:transigutierrez} yields:
 \begin{theorem}\label{T:nominimal2}\label{T:transitive2}
 Let $T=(\lambda,\pi)$ be a   $(2,k)$-\iet\ with $k\geq 1$. Then, $T$ is not transitive.
\end{theorem}

\begin{proof}
 \textcolor{black}{The property is obvious for an IET with associated permutation $(-2,-1)$. For the proper cases, } assume that   $T$ is transitive; then,  $T^c$ also is transitive. Observe, moreover, that $T^c$ is a proper transitive $(n',k')$-\cet\ with $n'+k'<5$ by Claim~\ref{R-minlinsol}.
 This contradicts Theorem~\ref{T:transigutierrez} and then   $T^c$ and $T$ are not transitive.
\end{proof}

\subsection{The case $n=3$}

Observe that the standard irreducible permutations  for $n=3$ are:
$
 \imt{clase7tamagno3[1]},\, \imt{clase7tamagno3[2]},\,\imt{clase7tamagno3[3]}.
$
Then, the cardinality of the set $S_3^{\sigma,*}$ of irreducible signed parmutations  is $2^3\cdot 3=24$.

By easy  computations with the help of the software \emph{Maxima}, we obtain eight Rauzy classes that we introduce \textcolor{black}{(recall that in Rauzy classes we only consider irreducible permutations jointly with the corresponding irreducible ones obtained by applying the maps $a$ and $b$)}:

\begin{itemize}
 \item $\mathcal{D}_1=\left\{\tau_{1,1}=\imt{clase1tamagno3[1]}\right\}$.
\item $\mathcal{D}_2=\left\{\tau_{2,1}=\imt{clase2tamagno3[1]},\tau_{2,2}=\imt{clase2tamagno3[2]},\tau_{2,3}=\imt{clase2tamagno3[3]}\right\}$.
\item $\mathcal{D}_3=\left\{\tau_{3,1}=\imt{clase3tamagno3[1]}\right\}$.
\item $\mathcal{D}_4=\left\{\tau_{4,1}=\imt{clase4tamagno3[1]},\tau_{4,2}=\imt{clase4tamagno3[2]}\right\}$.
\item 
$\mathcal{D}_5=\{\tau_{5,1}=\imt{clase5tamagno3[1]},$ $\tau_{5,2}=\imt{clase5tamagno3[2]},$ $\tau_{5,3}=\imt{clase5tamagno3[3]},$ $\tau_{5,4}=\imt{clase5tamagno3[4]},$ $
\tau_{5,5}=\imt{clase5tamagno3[5]},$ $\tau_{5,6}=\imt{clase5tamagno3[6]},$ $\tau_{5,7}=\imt{clase5tamagno3[7]},$ $\tau_{5,8}=\imt{clase5tamagno3[8]},$ $
\tau_{5,9}=\imt{clase5tamagno3[9]}
\}
$.
\item $\mathcal{D}_6=\left\{\tau_{6,1}=\imt{clase6tamagno3[1]},\tau_{6,2}=\imt{clase6tamagno3[2]},\tau_{6,3}=\imt{clase6tamagno3[3]}\right\}$.
\item $\mathcal{D}_7=\left\{\tau_{7,1}=\imt{clase7tamagno3[1]},\tau_{7,2}=\imt{clase7tamagno3[2]},\tau_{7,3}=\imt{clase7tamagno3[3]}\right\}$.
\item $\mathcal{D}_8=\left\{\tau_{8,1}=\imt{clase8tamagno3[1]},\tau_{8,2}=\imt{clase8tamagno3[2]}\right\}$.

\end{itemize}

\begin{remark}\label{R:Class7NoPro}
\textcolor{black}{Class $\mathcal{D}_7$ does not generate IETs with flips, so we leave it out of our study.}  \textcolor{black}{ For the rest of classes, we can also eliminate the non-proper $3$-IETs, since they can be included in the case of $2$-IETs, and by Theorem~\ref{T:transitive2} we know that transitivity is not reached. This is the case for \iets\ with associate permutations: $\tau_{1,1}$ in Class $\mathcal{D}_1$;
$\tau_{2,1}, \tau_{2,2}$ in Class $\mathcal{D}_2$; and $\tau_{6,1}, \tau_{6,2}$ in Class $\mathcal{D}_6$}. 

\end{remark}

Now, we can prove:
\begin{propo}\label{P:notransitamagno3parcial}
 Let $T=(\lambda,\pi)$ be a  proper $(3,k)$-\iet\ , with $k\geq 1$ and $\pi\in
 \{
 \tau_{3,1},$ $ \tau_{4,1},$ $\tau_{4,2},$ $ \tau_{5,7},$ $  \tau_{5,8},$ $
 \tau_{5,9},$ $
 \tau_{8,1}\}$. Then, $T$ is not transitive. 
\end{propo}
\begin{proof}
 If $T$ would be transitive, then $T^c$ would also be a proper transitive $(n',k')$-\cet\ with $n'+k'<5$ by Claim~\ref{R-minlinsol}, a 
 contradiction with Theorem~\ref{T:transigutierrez}. Then, $T^c$ and $T$ are not transitive.
\end{proof}

Also with the help of Lemma~\ref{C-rutafinita}, we can improve the previous result.

\begin{figure}[h]
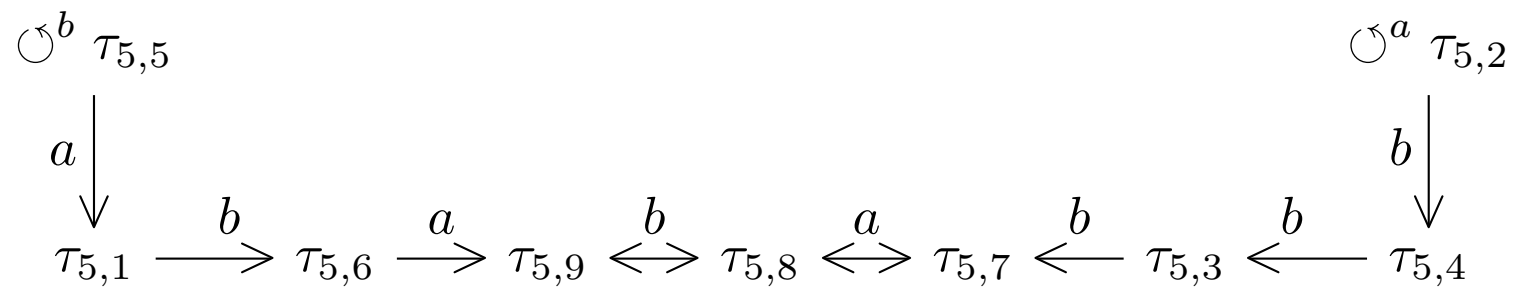

$
\begin{diagram}
\circlearrowleft^b\tau_{5,5}&&&&&& &&&&&&\circlearrowleft^a\tau_{5,2}\\
\dTo^a&&&&&& &&&&&&\dTo^b\\
\tau_{5,1}&\rTo^b&\tau_{5,6}&\rTo^a&\tau_{5,9}&
\rCorresponds^b&\tau_{5,8}&
\rCorresponds^a&\tau_{5,7}&
\lTo^b&\tau_{5,3}&
\lTo^b&\tau_{5,4}&
\end{diagram}
$
\caption{Rauzy graph of $\mathcal{D}_{5}$ \label{clase5tamagno3}}
\end{figure}

\begin{theorem}\label{T:nominimal3}
 Let $T=(\lambda,\pi)$ be a  proper $(3,k)$-\iet\ with $k\geq 1$. 
 Then, $T$ is not minimal.
\end{theorem}
\begin{proof}
\textcolor{black}{The non-minimality when  $\pi=\tau_{1,1}$ is obvious by Remark~\ref{R:Class7NoPro}. Concerning Class $\mathcal{D}_2$, since $a(\tau_{2,1})=\tau_{2,3}$, from Lemma~\ref{C-rutafinita} and Remark~\ref{R:Class7NoPro} we deduce that the corresponding $3$-IETs are not minimal. By a similar reasoning, Class~$\mathcal{D}_6$ does not generate minimal IETs as
$b(\tau_{6,3})=\tau_{6,2}$. The proofs for Class~$\mathcal{D}_3$  and Class~$\mathcal{D}_4$ are  covered by Proposition~\ref{P:notransitamagno3parcial}.  Since $b(\tau_{8,2})=\tau_{8,1}$, we apply jointly Theorem~\ref{T:nominimal3} and Proposition~\ref{P:notransitamagno3parcial} and we obtain that the corresponding \iets\ of Class~$\mathcal{D}_8$ are not minimal.}

\textcolor{black}{Let us focus in   Class~$\mathcal{D}_5$. By Proposition~\ref{P:notransitamagno3parcial}, if
$\pi\in\{\tau_{5,7},\tau_{5,8},\tau_{5,9}\}$  then $T$ is not minimal. }

Finally, assume that $T=(\lambda,\pi)$ is minimal for some 
 $\pi\in\{\tau_{5,1},\tau_{5,2},\tau_{5,3},\tau_{5,4},\tau_{5,5},
 \tau_{5,6}\}$. Then, by Proposition~\ref{P:minimalidadimplicainducible}  $R^k(T)=(\lambda^k,\pi_k)$ is defined for any $k\in\N$ and, since the type of $R^k(T)$ is not eventually constant, it holds $\pi_{k_0}=\tau_{5,8}$ for some $k_0$  (see 
 Figure~\ref{clase5tamagno3}); but this is a contradiction with the previous paragraph. Then $T=(\lambda,\pi)$ is not minimal and 
 the result is proved.
\end{proof}

%

If we take into account oriented $n$-\iets, $2\leq n\leq 3$, then it is clear that almost all   \iets\ are minimal, however they are no 
non uniquely ergodic.
%
%

\begin{theorem}\label{T:Weyl} Let $n\in\{2,3\}$ and let $T:D\subset [0,\ell]\to [0,\ell]$ be a minimal oriented $n$-\iet. 
Then,  $T$ is uniquely ergodic.
\end{theorem}

\begin{proof}
 This is a consequence of  Weyl equidistribution theorem, see \cite{keane1}.
\end{proof}

\section{Rauzy classes for $n=4$}
In \cite{VincentDelecroixCardinality}, it is stated that the cardinality of irreducible standard permutations, $f(n)$, is given by the following recurrence formula:
$$f(n)=n!-\sum_{k=1}^{n-1}f(k)(n-k)!$$
Since $f(1)=f(2)=1$, we can easily compute $f(3)=3$ and $f(4)=13$. These standard irreducible permutations  for $n=4$ are:
\begin{eqnarray*}
 \imt{clase16[1]},\, \imt{clase16[2]},\,\imt{clase16[3]},\,\imt{clase16[4]},\,\imt{clase16[5]},\,\imt{clase16[6]},\,\\
 \imt{clase17[1]},\, \imt{clase17[2]},\, \imt{clase17[3]},\, \imt{clase17[4]},\, \imt{clase17[5]},\,\imt{clase17[6]},\,\imt{clase17[7]}.
\end{eqnarray*}
Then the cardinality of signed irreducible permutations $S_4^{\sigma,*}$ is $2^4\cdot 13=208$.

By easy  computations with the help of the software \emph{Maxima}, we obtain twenty Rauzy classes that we analyze  in the following subsections. We have labelled them by $\mathcal{C}_i$, $1\leq i\leq 20$.

Recall that, in the orientable case, each Rauzy class contains a \emph{standard pair}, that is, a permutation $\pi\in S_n^*$ satisfying $\pi(n)=1$ and $\pi(1)=n$.
We can generalize this notion to signed permutations. We will say that
$\pi\in\permuset$ is a standard pair if $|\pi(n)|=1$ and $|\pi(1)|=n$. It is known that any orientable class contains a standard pair, see \cite[Proposition 6.2]{viana}. This is the case also for the generalized Rauzy classes when $n=4$. We write in Table~\ref{tablestandardpairs} a standard pair of each class and the cardinality of the class.

%
%

\begin{table}[h]
  \begin{tabular}{|c|c|c|||c|c|c|||c|c|c|}
 \hline
$i$&Standard pair in $\mathcal{C}_i$& $\# \mathcal{C}_i$&$i$&Standard pair in $\mathcal{C}_i$& $\# \mathcal{C}_i$&$i$&Standard pair in $\mathcal{C}_i$& $\# \mathcal{C}_i$\\  
 \hline
$1$&$(-4 , -3 , -2 , -1)$& $1$&$2$&$(-4 , -3 , 2 , -1)$& $1$&$3$&$(-4 , -2 , -3 , -1)$& $1$\\
\hline
$4$&$(-4 , -2 , 3 , -1)$& $1$&$5$&$(-4,-3,-2,1)$& $6$&$6$&$(-4 , 2 , -3 , -1)$& $1$\\
\hline
$7$&$(-4 , 2 , 3 , -1)$& $1$&$8$&$(-4,-3,2,1)$& $8$&$9$&$(-4 , 3 , -2 , -1)$& $1$\\
\hline
$10$&$(-4 , 3 , 2 , -1)$& $1$&$11$&$(-4,3,2,1)$& $4$&$12$&$(4,3,-2,1)$& $29$\\
\hline
$13$&$(4,-3,-2,1)$& $36$&$14$&$(-4,3,-2,1)$& $27$&$15$&$(4,-3,-2,-1)$& $6$\\
\hline
$16$&$(4,2,3,1)$& $6$&$17$&$(4,3,2,1)$& $7$&$18$&$(4,3,2,-1)$& $4$\\
\hline
$19$&$(-4,-2,-3,1)$&$59$&$20$ &$(4,3,-2,-1)$ &$8$&&&\\
\hline
\end{tabular}

\caption{Cardinality of Rauzy classes for $n=4$ and standard pairs \label{tablestandardpairs}}

\end{table}

\subsection{Classes with one element}

The following classes has only one permutation and there are no arrows since by applying either $a$ or $b$ to each permutation we obtain a reducible one. Then, according to Proposition~\ref{P:minimalidadimplicainducible}, no mininimal \iet\ can have these permutations as combinatorial data. Thus, we have proved the following proposition.
\begin{propo}\label{clases1element}
 Let $T=(\lambda,\pi)$, with $\pi\in\left(\bigcup_{i=1}^4\mathcal{C}_{i}\right)\cup\mathcal{C}_{6}
 \cup\mathcal{C}_{7}\cup\mathcal{C}_{9}\cup\mathcal{C}_{10}$, where 
\begin{eqnarray*}
\mathcal{C}_1 &=&   \{(-4 , -3 , -2 , -1)\}, \,
\mathcal{C}_2=   \{(-4 , -3 , 2 , -1)\}, \,
\mathcal{C}_3=   \{(-4 , -2 , -3 , -1)\}, \,
\mathcal{C}_4=   \{(-4 , -2 , 3 , -1)\}, \\
\mathcal{C}_6 &=&   \{(-4 , 2 , -3 , -1)\}, \,
\mathcal{C}_7=   \{(-4 , 2 , 3 , -1)\}, \,
\mathcal{C}_9=   \{(-4 , 3 , -2 , -1)\}, \,
\mathcal{C}_{10}=\{(-4 , 3 , 2 , -1)\}.
\end{eqnarray*}
 Then, $T$ is not minimal.
\end{propo}

\subsection{Classes with four elements, $\mathcal{C}_{11}$ and $\mathcal{C}_{18}$ \label{S:clase11clase18}}
\begin{figure}[h]
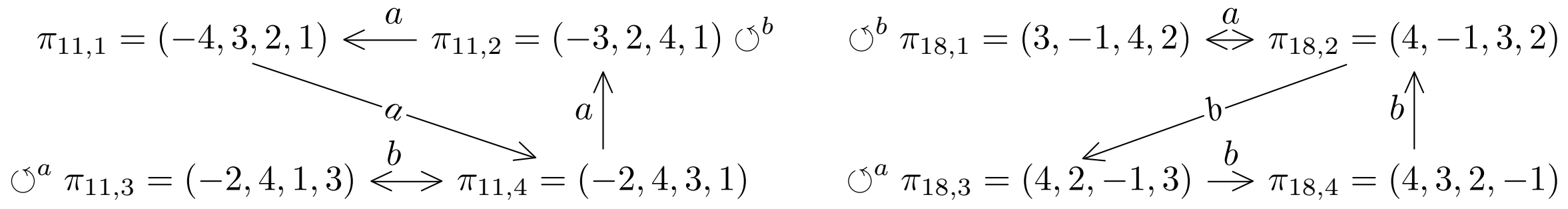

$
\begin{diagram}
\pi_{11,1}=\imt{clase11[1]}&\lTo^a&\pi_{11,2}=\imt{clase11[2]}\circlearrowleft^b\\
 & \rdTo(2,2)~{a} &\uTo^a\\
 \circlearrowleft^a\pi_{11,3}=\imt{clase11[3]}&\rCorresponds^b&\pi_{11,4}=\imt{clase11[4]}\\
\end{diagram}
\qquad
\begin{diagram}
\circlearrowleft^b\pi_{18,1}=\imt{clase18[1]}&\rCorresponds^a&\pi_{18,2}=\imt{clase18[2]}\\
 & \ldTo(2,2)~{b} &\uTo^b\\
 \circlearrowleft^a\pi_{18,3}=\imt{clase18[3]}&\rTo^b&\pi_{18,4}=\imt{clase18[4]}\\
\end{diagram}
$
\caption{$\mathcal{C}_{11}$ and $\mathcal{C}_{18}$ \label{clase11}}
\end{figure}

We introduce the Rauzy graph of    these classes  in Figure~\ref{clase11}.

Observe that, by Claim~\ref{Cl:ExchangeIntervalsBya}, in class $\mathcal{C}_{11}$ the first interval is never shortened, and then the combinatorial data therein do not generate minimal \iets\ by Proposition~\ref{P:minimalidadimplicainducible}. Now, note that $\pi_{11,1}^h=\pi_{18,4}$; \textcolor{black}{then, by Theorem~\ref{Th:hlamda}, $\pi_{18,4}$ cannot be the associated permutation of a minimal $4$-IET.} Finally, since we can connect any vertex of $\mathcal{C}_{18}$ with arrows begining in $\pi_{18,4}$, we obtain that there are not minimal \iets\ with combinatorial data in $\mathcal{C}_{18}$ by applying Lemma~\ref{C-rutafinita}, Proposition~\ref{P:minimalidadimplicainducible} and Theorem~\ref{T-cone}. Therefore, the following proposition is proved.

\begin{propo}\label{clases11Y18element}
 Let $T=(\lambda,\pi)$, with $\pi\in\mathcal{C}_{11}\cup\mathcal{C}_{18}.$ Then, $T$ is not minimal.
\end{propo}

\subsection{Classes 5 and 15}

\begin{figure}[h]
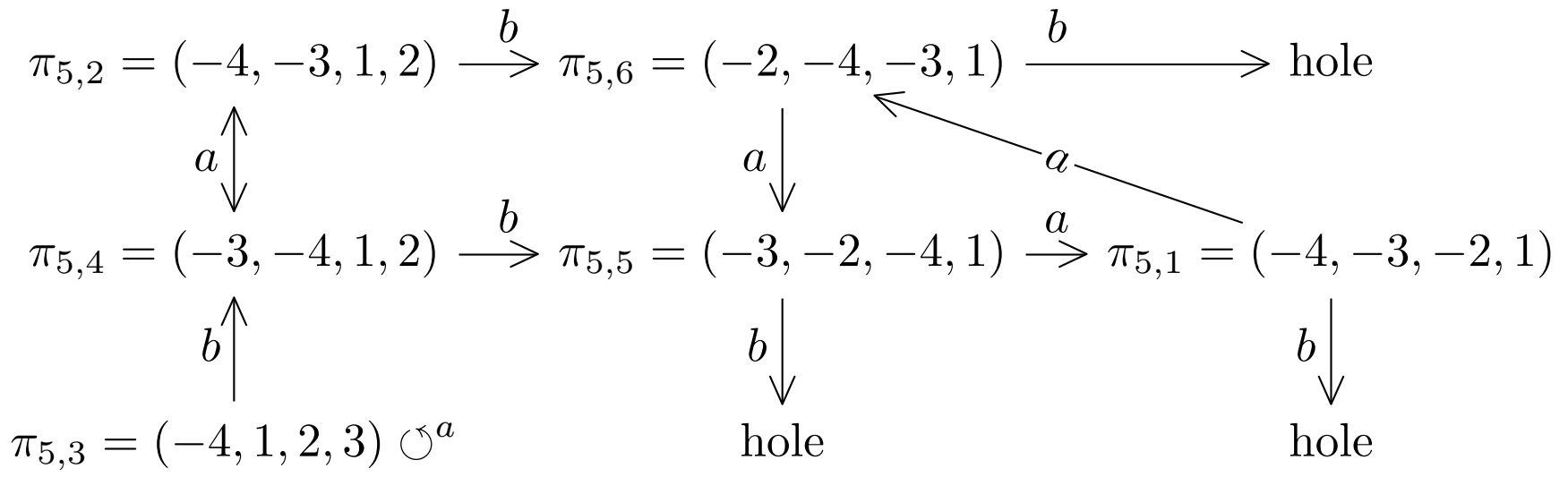

$
\begin{diagram}
\pi_{5,2}=\imt{clase5[2]}&\rTo^b&\pi_{5,6}=\imt{clase5[6]}&\rTo^b&\hole&\\
\dCorresponds^a&&\dTo^a&\luTo~{a}&&\\
\pi_{5,4}=\imt{clase5[4]}&\rTo^b&\pi_{5,5}=\imt{clase5[5]}&\rTo^a&\pi_{5,1}=\imt{clase5[1]}\\
\uTo^b&&\dTo^b&&\dTo^b\\
\pi_{5,3}=\imt{clase5[3]}\circlearrowleft^a&&\hole&&\hole\\
\end{diagram}
$

\caption{Rauzy graph of $\mathcal{C}_5$ \label{clase5}}
\end{figure}

Observe that any compatible path with minimality  must have infinitely many arrows labelled by $a$ and also by $b$ by Proposition~\ref{P:minimalidadimplicainducible}, but it is not possible in the Rauzy graph of Figure~\ref{clase5} \textcolor{black}{since infinite paths are eventually of the form $(aa)^\infty$ or $(aaa)^\infty$}.

Class 15 is the following:

$$\begin{array}{c}
  \mathcal{C}_{15}=\left\{(2 , 3 , 4 , -1),( 3 , 4 , -2 , -1
 ),( 3 , 4 , -1 , -2),( 4 , -3 , -2 , -
 1),( 4 , -2 , -1 , -3),( 4 , -1 , -3
  , -2) \right\}.
 \end{array}
$$

All the elements of Class 15 are conjugated by elements of Class 5 and, \textcolor{black}{by Theorem~\ref{Th:hlamda},} they do  not generate minimal \iets. \textcolor{black}{To be precise:}

$$\begin{array}{lll}
\pi_{5,1}^h=\imt{clase15[4]},&
\pi_{5,2}^h=\imt{clase15[2]}, &
\pi_{5,3}^h=\imt{clase15[1]},\\
\pi_{5,4}^h=\imt{clase15[3]},&
\pi_{5,5}^h=\imt{clase15[6]},&
\pi_{5,6}^h=\imt{clase15[5]}.
 \end{array}
$$

We have finally proved the following result.
\begin{propo}\label{clases5Y15}
 Let $T=(\lambda,\pi)$, with $\pi\in\mathcal{C}_{5}\cup\mathcal{C}_{15}$. Then, $T$ is not minimal.
\end{propo}

\subsection{Class 8} The vertices in this class are the following:
$\mathcal{C}_8=\{(-4 , -3 , 2 , 1),$ $( -4 , 2 , 3 , 1
 ),$ $( -4 , 3 , 1 , 2),$ $( -3 , -2 , 4 , 1
 ),$ $( -3 , 4 , 1 , 2),$ $( -3 , 4 , 2 , 1
 ),$ $( -2 , -4 , 3 , 1),$ $( -2 , 3 , 4 , 1
 )\}.$

\begin{figure}[h]
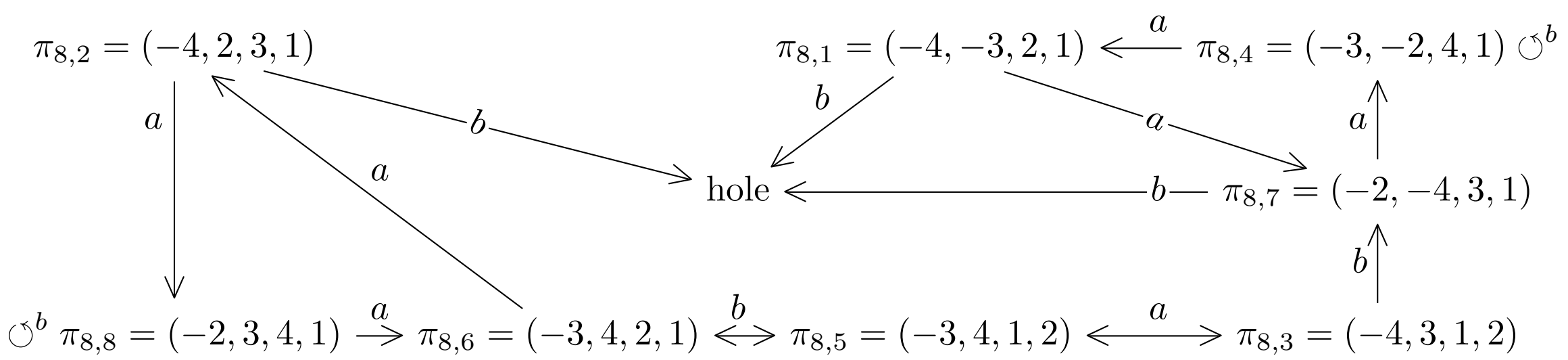

$$
\begin{diagram}
\pi_{8,2}=\imt{clase8[2]}&&&&\pi_{8,1}=\imt{clase8[1]}&&\lTo^a&\pi_{8,4}=\imt{clase8[4]}\circlearrowleft^b&
\\
\dTo^a&\luTo(2,4)^a\rdTo(3,2)~b&&\ldTo(1,2)^b&&      \rdTo(3,2)~{a}         &&\uTo^a&
\\
&&&\hole&&&\lTo~b&\pi_{8,7}=\imt{clase8[7]}&
\\
&&&&&&&\uTo^b&
\\
\circlearrowleft^b\pi_{8,8}=\imt{clase8[8]}&\rTo^a&\pi_{8,6}=\imt{clase8[6]}&\rCorresponds^b&\pi_{8,5}=\imt{clase8[5]}&&\rCorresponds^a&\pi_{8,3}=\imt{clase8[3]}&\\
\end{diagram}
$$
\caption{Rauzy class $\mathcal{C}_{8}$ \label{clase8}}
\end{figure}

In Figure~\ref{clase8}, we draw the Rauzy graph corresponding to this class.

\begin{propo}\label{P:clase8}
 Let $T=(\lambda,\pi)$, with $\pi\in\mathcal{C}_{8}$. Then, $T$ is not minimal.
\end{propo}

\begin{proof}

Observe that in this class we can build infinite paths in the Rauzy graph; however, in any case the first interval is never shortened except when  we are in $\pi_{8,3}$ and we apply $b$, see Claim~\ref{Cl:ExchangeIntervalsBya}. Then, we would arrive to $\pi_{8,7}$ and we never return to $\pi_{8,3}$: consequently, the first interval in an associated \iet\ to  $\pi_{8,7}$ is never the winner, which contradicts the existence of a minimal \iet\ having combinatorial data in this class \textcolor{black}{since by Proposition~\ref{P:minimalidadimplicainducible} the lengths of the vectors $\lambda$ tends to zero.} The proposition is proved.

\end{proof}

\subsection{Class 12} This class has 29 vertices, the following ones:

 \begin{eqnarray*}
\mathcal{C}_{12}=
\{
 \pi_{12,1}=\imt{clase12[1]},\pi_{12,2}=\imt{clase12[2]},\pi_{12,3}=\imt{clase12[3]},\pi_{12,4}=\imt{clase12[4]},\\
 \pi_{12,5}=\imt{clase12[5]},
 \pi_{12,6}=\imt{clase12[6]},
\pi_{12,7}=\imt{clase12[7]},\pi_{12,8}=\imt{clase12[8]},\\\pi_{12,9}=\imt{clase12[9]},\pi_{12,10}=\imt{clase12[10]},
\pi_{12,11}=\imt{clase12[11]},\pi_{12,12}=\imt{clase12[12]},\\
\pi_{12,13}=\imt{clase12[13]},\pi_{12,14}=\imt{clase12[14]},\pi_{12,15}=\imt{clase12[15]},
\pi_{12,16}=\imt{clase12[16]},\\
\pi_{12,17}=\imt{clase12[17]},\pi_{12,18}=\imt{clase12[18]},
\pi_{12,19}=\imt{clase12[19]},\pi_{12,20}=\imt{clase12[20]},\\
\pi_{12,21}=\imt{clase12[21]},\pi_{12,22}=\imt{clase12[22]},\pi_{12,23}=\imt{clase12[23]},\pi_{12,24}=\imt{clase12[24]},\\
\pi_{12,25}=\imt{clase12[25]},
\pi_{12,26}=\imt{clase12[26]},\pi_{12,27}=\imt{clase12[27]},\\
\pi_{12,28}=\imt{clase12[28]}, \pi_{12,29}=\imt{clase12[29]}
\}.
 \end{eqnarray*}

We are going to prove \textcolor{black}{that this class does not provide combinatorial data for minimal IETs.}

\begin{propo}\label{P:clase12}
 Let $T=(\lambda,\pi)$, with $\pi\in\mathcal{C}_{12}$. Then, $T$ is not minimal.
\end{propo}

\begin{proof}

It is important to notice the following relations:

$$\pi_{12,15}^h=\pi_{18,3},\,\,
\pi_{12,18}^h=\pi_{11,3},\,\,
\pi_{12,22}^h=\pi_{18,1},\,\,
\pi_{12,23}^h=\pi_{18,2},\,\,
\pi_{12,27}^h=\pi_{11,2},\,\,
\pi_{12,28}^h=\pi_{11,4}.$$


$\imaximaenvenv{
  }
$

Since the combinatorial data in Classes 18 and 11 are not compatible with minimal \iets\ as we showed in Section~\ref{S:clase11clase18}, \textcolor{black}{by Theorem~\ref{Th:hlamda}} we have that any \iet\ having as combinatorial data the permutations $\pi_{12,15},$ 
$\pi_{12,18},$ 
$\pi_{12,22},$ 
$\pi_{12,23},$ 
$\pi_{12,27}$ or 
$\pi_{12,28}$ is not minimal. Let $N^{(12)}$ be the set  of natural numbers $1\leq j\leq 29$ for which the permutation $\pi_{12,j}$ does not admit minimal \iets. 
 \textcolor{black}{From the above,}  $\{15,18,22,23,27,28\}\subseteq N^{(12)}$.

\begin{figure}[h]
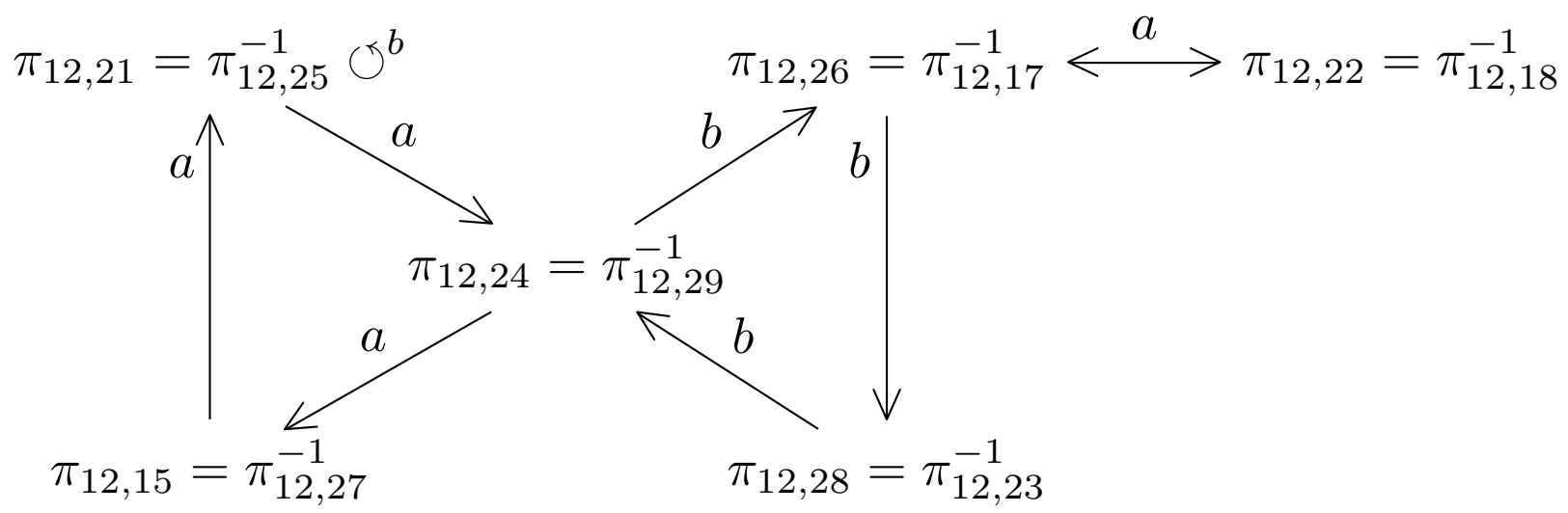

$
\begin{diagram}
&&\pi_{12,21}=\pi_{12,25}^{-1}\circlearrowleft^b&&&&\pi_{12,26}=\pi_{12,17}^{-1} &\rCorresponds[l>=3em]^a& \pi_{12,22}=\pi_{12,18}^{-1}
\\
&&\uTo^a&\rdTo^a&&\ruTo^b
&\dTo^b&&&&&&\\
&&&&\pi_{12,24}=\pi_{12,29}^{-1}
\\
&&&\ldTo^a&&\luTo^b&\\
&&\pi_{12,15}=\pi_{12,27}^{-1}&&&&\pi_{12,28}=\pi_{12,23}^{-1}
\end{diagram}
$

\caption{Subgraph of Rauzy class $\mathcal{C}_{12}$  \label{subclase12}}
\end{figure}

\textrm{
$\begin{maximaenv}
 \end{maximaenv}
$
}

Observe Figure~\ref{subclase12} and recall that the permutations $\pi_{12,22}$, $\pi_{12,15}$ and $\pi_{12,28}$ do not generate minimal \iets.  Then, the other vertices do not generate minimal \iets\ since we can arrive to them from $\pi_{12,15}$ by means of $a$ and $b$, and we can use Lemma~\ref{C-rutafinita},  Proposition~\ref{P:minimalidadimplicainducible} and Theorem~\ref{T-cone} to conclude their non-minimality. Then, bear in mind also Theorem~\ref{Th:hlamda}, we have $\imt{ctoN}\subseteq N^{(12)}.$

\begin{figure}[h]
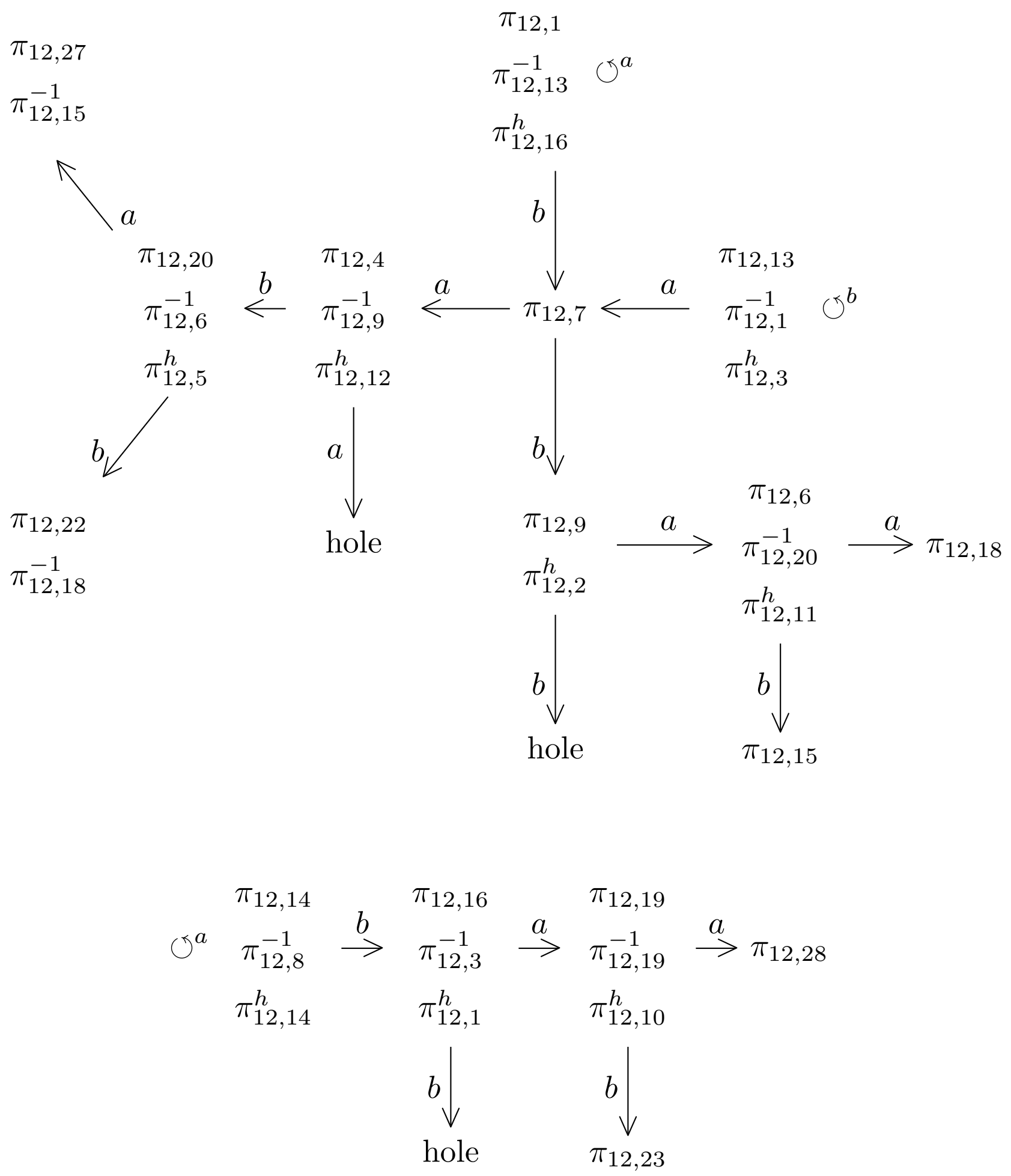

$
\begin{diagram}
\permuclaseinver{12}{27}{15}&&&&&&\permuclaseinverh{12}{1}{13}{16} \circlearrowleft^a &\\
&&&&&&&\\
&\luTo^a&&&&&\dTo^b&\\
&&\permuclaseinverh{12}{20}{6}{5}&\lTo^b&\permuclaseinverh{12}{4}{9}{12}&\lTo^a&\pi_{12,7}&\lTo^a&\permuclaseinverh{12}{13}{1}{3}\circlearrowleft^b\\&&&&\\
&\ldTo^b&&&\dTachado^a&&\dTo^b&\\
\permuclaseinver{12}{22}{18}&&&&\hole&&\permuclaseh{12}{9}{2}&\rTo^a&\permuclaseinverh{12}{6}{20}{11}&\rTo^a&\pi_{12,18}\\
&&&&&&&&&\\
&&&&&&\dTo^b&&\dTo^b&\\
&&&&&&\hole&&\pi_{12,15}\\
\end{diagram}
$

\vspace*{1cm}

$
\begin{diagram}
\circlearrowleft^a\permuclaseinverh{12}{14}{8}{14}
&\rTo^b&\permuclaseinverh{12}{16}{3}{1}&\rTo^a&\permuclaseinverh{12}{19}{19}{10}&\rTo^a&\pi_{12,28}\\
&&&&&&\\
&&\dTo^b&&\dTo^b\\
&&\hole&&\pi_{12,23}
\end{diagram}
$

\caption{Subgraphs of Rauzy class $\mathcal{C}_{12}$ \label{subclase12B}}
\end{figure}

\textrm{
$\begin{maximaenv}
 \end{maximaenv}
$
}

Now, we take into account the Rauzy subgraphs in Figure~ \ref{subclase12B}. Since $\pi_{12,18}$, $\pi_{12,27}$, $\pi_{12,22}$, $\pi_{12,23}$ and $\pi_{12,28}$ do  not admit minimal \iets, again the combination of Proposition~\ref{P:minimalidadimplicainducible}, Theorem~\ref{T-cone} and Lemma~\ref{C-rutafinita}, \textcolor{black}{jointly with Theorem~\ref{Th:hlamda},} gives 
$$\imt{ctoN}\subseteq N^{(12)}.$$

%
Therefore,   $\{1,2,3,\dots,29\}= N^{(12)}$ and there are no minimal \iets\ $T=(\lambda,\pi)$ with $\pi\in\mathcal{C}_{12}$. 
\textcolor{black}{This ends the proof.}

%
%

\end{proof}

\subsection{Class 14}
The list of   $27$ permutations from this class is the following:

\begin{eqnarray*}
    \mathcal{C}_{14}=\{\pi_{14,1}=\imt{clase14[1]},\pi_{14,2}=\imt{clase14[2]},\pi_{14,3}=\imt{clase14[3]},\pi_{14,4}=\imt{clase14[4]},\\\pi_{14,5}=\imt{clase14[5]},\pi_{14,6}=\imt{clase14[6]},
\pi_{14,7}=\imt{clase14[7]},\pi_{14,8}=\imt{clase14[8]},\\
\pi_{14,9}=\imt{clase14[9]},\pi_{14,10}=\imt{clase14[10]},\pi_{14,11}=\imt{clase14[11]},\pi_{14,12}=\imt{clase14[12]},\\
\pi_{14,13}=\imt{clase14[13]},\pi_{14,14}=\imt{clase14[14]},\pi_{14,15}=\imt{clase14[15]},\pi_{14,16}=\imt{clase14[16]},\\
\pi_{14,17}=\imt{clase14[17]},\pi_{14,18}=\imt{clase14[18]},
\pi_{14,19}=\imt{clase14[19]},\pi_{14,20}=\imt{clase14[20]},\\
\pi_{14,21}=\imt{clase14[21]},\pi_{14,22}=\imt{clase14[22]},\pi_{14,23}=\imt{clase14[23]},\pi_{14,24}=\imt{clase14[24]},\\
\pi_{14,25}=\imt{clase14[25]},\pi_{14,26}=\imt{clase14[26]},\pi_{14,27}=\imt{clase14[27]} \}.
 \end{eqnarray*}

 We begin the analysis of this class by studying some parts of the Rauzy graph. In particular, we are going to analyze some parts of the graph that we call \emph{nodes}. A node $\mathcal{N}$ is a subgraph such that the images of the vertices of $\mathcal{N}$, by $a$ and $b$, still remain in $\mathcal{N}$, \textcolor{black}{unless we get a hole, a reducible permutation}.

 Let $N^{(14)}$ be the set  of natural numbers $1\leq j\leq 27$ for which the permutation $\pi_{14,j}$ does not admits minimal \iets.

 \subsubsection{Analysis of nodes}

\begin{figure}[h]
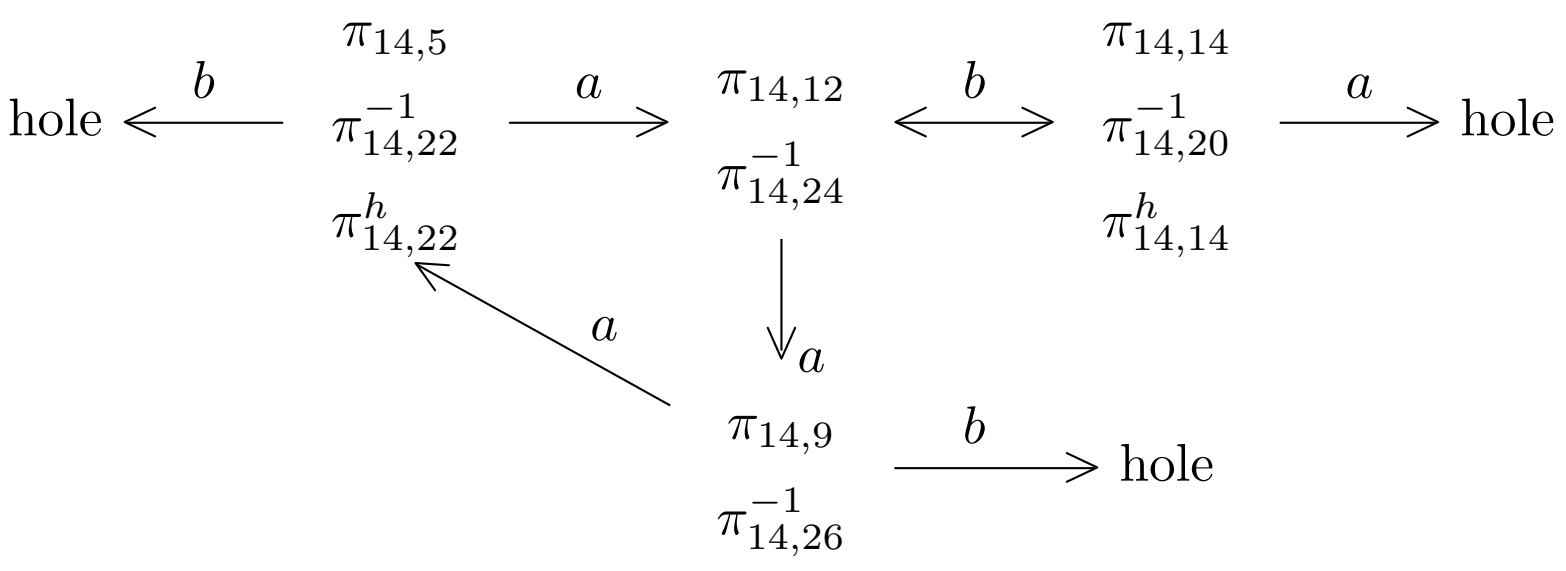

 $
\begin{diagram}
\hole&\lTo[l>=3em]^b&\permuclaseinverh{14}{5}{22}{22}&\rTo[l>=3em]^a&\permuclaseinver{14}{12}{24}&\rCorresponds[l>=3em]^b&\permuclaseinverh{14}{14}{20}{14}&\rTo[l>=3em]^a&\hole
\\
&&&&&&&&&\\
&&&\luTo^a&\dTo_a&\\
&&&&\permuclaseinver{14}{9}{26}&\rTo^b&\hole
\end{diagram}
$
 \caption{Nodes in $\mathcal{C}_{14}$\label{nodes}}
\end{figure}

We introduce a node in Figure~\ref{nodes}.
%

\textrm{
\begin{maximaenv}
\end{maximaenv}
}

Observe that the vertices of any path $v\in\{a,b\}^{\N}$ beginning in one of the vertices of the node of Figure~\ref{nodes} always remain in this node. Moreover, use Claim~\ref{Cl:ExchangeIntervalsBya} to deduce that  the first interval of an \iet\ $T=(\lambda,\pi)$ having combinatorial data in this node is never shortened when passing throw $a$. Moreover, also by Claim~\ref{Cl:ExchangeIntervalsBya}, this first interval is not shortened by $b$ since $\pi_{14,12}=\imt{clase14[12]}$ and $\pi_{14,14}=\imt{clase14[14]}$. Therefore,
\begin{equation}\label{E:clase14-A}
\imt{ctoN}\subseteq N^{(14)}.
\end{equation}

\subsubsection{Other paths}

\begin{figure}[h]
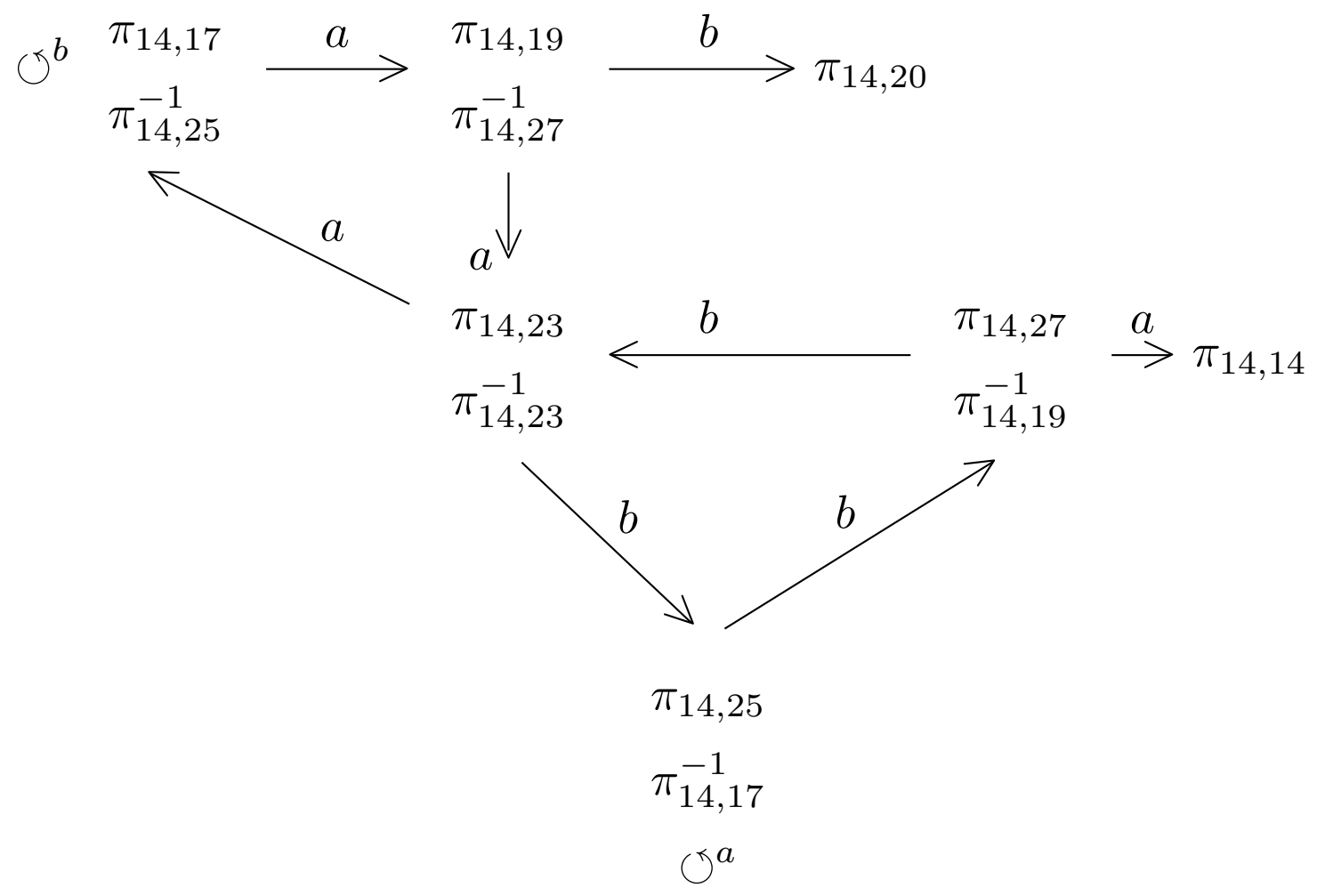


 $
\begin{diagram}
\circlearrowleft^b\permuclaseinver{14}{17}{25}&
\rTo[l>=3em]^a&\permuclaseinver{14}{19}{27}&&\rTo^b&\pi_{14,20}\\
&&&&&&&&&&&&&\\
&\luTo^a&\dTo^a&&
\\
&&\permuclaseinver{14}{23}{23}&&\lTo^b&
&\permuclaseinver{14}{27}{19}&\rTo^a&\pi_{14,14}&\\
&&&&&&&&&
\\
&&&\rdTo^b&
&\ruTo^b&&&&
\\
&&&&&&&&&\\
&&&&
\begin{array}{c}
\permuclaseinver{14}{25}{17}\\
\circlearrowleft^a
\end{array}
&&&&&
\end{diagram}
$
 \caption{Subgraph in $\mathcal{C}_{14}$\label{clase14A}}
\end{figure}

An \iet\ with combinatorial data in the vertices of Figure~\ref{clase14A} can  finish after finitely many times of applying the Rauzy induction in an \iet\ with permtutation $\pi_{14,20}$ or $\pi_{14,14}$,  and then it is not minimal by  Lemma~\ref{C-rutafinita} and because $\{14,20\}\subset N^{(14)}$. The other possibility is to apply to the initial \iet\  infinitely many times  $a$ and $b$. Observe that the permutations to which we apply $b$ are
$\pi_{14,17}=\imt{clase14[17]}$,
$\pi_{14,23}=\imt{clase14[23]}$,
$\pi_{14,27}=\imt{clase14[27]}$,
$\pi_{14,25}=\imt{clase14[25]}$. 
We analyze in Table~\ref{T:binclase14} -for this, remind Claim~\ref{Cl:ExchangeIntervalsBya}- what happens to an \iet\ $T=(\lambda,\pi)$ when we apply $b$ to a permutation $\pi$ given by $\pi_{14,23}$, $\pi_{14,25}$, $\pi_{14,27}$ and  $\pi_{14,17}$.  On the other hand, when we apply the map $a$
 to an \iet\
 $T=(\lambda,\pi)$ we obtain $R(T)=(\lambda',\pi')$, with
\begin{equation}
(\lambda_1',\lambda_2',\lambda_3',\lambda_4')=(\lambda_1,\lambda_2,\lambda_3,\lambda_4-\lambda_{|\pi|^{-1}(4)}). \label{E:applicationofa}
\end{equation}

\begin{table}[h]
\begin{tabular}{|c|c|}
\hline
$\pi$&$\lambda'$\\
\hline
$\pi_{14,23}$&$(\lambda_1-\lambda_4,\lambda_4,\lambda_2,\lambda_3)$\\
\hline
$\pi_{14,25}$&$(\lambda_1-\lambda_4,\lambda_4,\lambda_2,\lambda_3)$\\
\hline
$\pi_{14,27}$&$(\lambda_1-\lambda_4,\lambda_4,\lambda_2,\lambda_3)$\\
\hline
$\pi_{14,17}$&$(\lambda_1,\lambda_2,\lambda_3-\lambda_4,\lambda_4)$\\
\hline
\end{tabular}
\vspace*{.5cm}
 \caption{Relation between $\lambda$ and $\lambda'$ for $T=(\lambda,\pi)$,
 $R(T)=(\lambda',b(\pi))$  and $\lambda=(\lambda_1,\lambda_2,\lambda_3,\lambda_4)$\label{T:binclase14}}
\end{table}

We are going to see that an \iet\ $T=(\lambda,\pi_{14,23})$ is not minimal. Assume the opposite; then, by Proposition~\ref{P:minimalidadimplicainducible} we can apply infinitely many times the Rauzy operator to $T$, i.e., $R^m(T)=(\lambda^m,\pi_m)$ exists for any $m\in\N$ and $|\lambda^m|$ goes to zero as $m$ goes to $+\infty$. Let $v$ be the path followed by $T$, that is, the one  satisfying  $\mathcal{G}^{T}=\mathcal{G}^{\pi_{14,23},v}$. Note  that the path $v$ is going to  travel infinitely many times in the upper or in lower cycle of  Figure~\ref{clase14A} or in both; moreover, infinitely many times $\pi_m$ coincides with $\pi_{14,23}$.

When we follow the upper cycle with the Rauzy induction on $(\lambda,\pi_{14,23})$ and, after applying $k$ times the operator of Rauzy,  come back  to
$T^{(k)}=(\lambda^{k}, \pi_{k})$, we have applied $3r$ times, $r>0$, the operator with the map $a$ and $s$ times,  $s\geq 0$, the operator with the map $b$ and $3r+s=k$.
Taking into account
(\ref{E:applicationofa}) and Table~\ref{T:binclase14}, we deduce that
\begin{equation}\label{E:shorten1}
\lambda^k=(\lambda_1,\lambda_2,\lambda_3^{(s)},\lambda_4^{(3r)}),
\end{equation}
where $\lambda_j^{(t)}$ means that we have shortened  $t$ times the length of $\lambda_j$, with lengths depending on the loser interval in each step of the induction. In any case, it is clear that  the first and second interval have not been shortened.
If we follow the lower cycle with the Rauzy induction on $(\mu,\pi_{14,23})$ and come back -after applying $k$ times the operator of Rauzy- to
$T^{(l)}=(\mu^{l}, \pi_{l})$, we have applied $3r$ times, $r>0$, the operator with the map $b$ and $s$ times,  $s\geq 0$, the operator with the map $a$ and $3r+s=k$.
Taking into account
(\ref{E:applicationofa}) and Table~\ref{T:binclase14}, we deduce now that
\begin{equation}\label{E:shorten2}\lambda^l=(\lambda_1^{(3r)},\lambda_2,\lambda_3^{(s)},\lambda_4),
\end{equation}
which again means that the second and fourth intervals have not been shortened.
From (\ref{E:shorten1})-(\ref{E:shorten2}), we have that  $|\lambda^m|>\lambda_2$ for any $m\in\N$, which contradicts that $|\lambda^m|$ goes to zero as $m$ goes to $\infty$. Then, $T=(\lambda,\pi_{14,23})$ is not minimal. Moreover, with the help of   Lemma~\ref{C-rutafinita}, any \iet\ with associated permutation in $\{\pi_{14,17},
\pi_{14,25},\pi_{14,19},\pi_{14,27}\}$ (see Figure~\ref{clase14A}) is not minimal.

%
\textrm{
$\begin{maximaenv}
ctoN:union({17,25,19,27,23,27,19,25,17},ctoN)
\end{maximaenv}$
}
Therefore we have obtained
\begin{equation}\label{E:nominimal14B}
\imt{ctoN}\subseteq N^{(14)}.
\end{equation}

\begin{figure}[h]
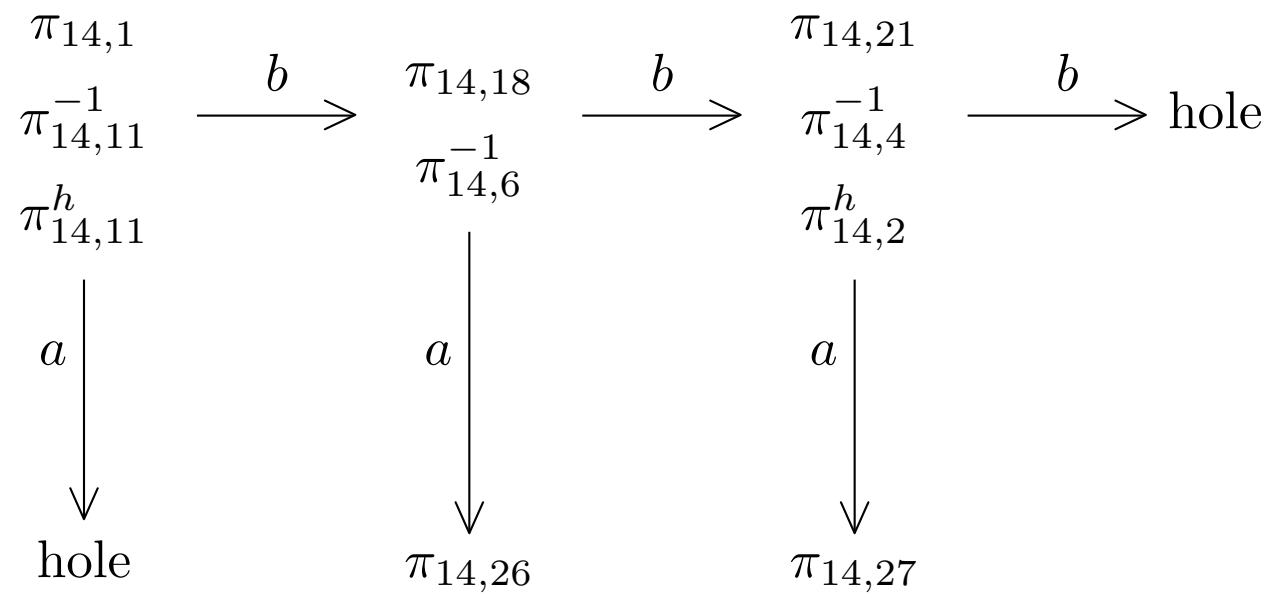

$
 \begin{diagram}
\permuclaseinverh{14}{1}{11}{11}&\rTo[l>=3em]^b&\permuclaseinver{14}{18}{6}&\rTo[l>=3em]^b&\permuclaseinverh{14}{21}{4}{2}&&\rTo[l>=3em]^b&\hole\\
&&&&&&&&\\
\dTo^a&&\dTo^a&&\dTo^a\\
&&&&&&&&\\
\hole&&\pi_{14,26}&&\pi_{14,27}
 \end{diagram}
$
\caption{Subgraph in $\mathcal{C}_{14}$\label{clase14B}}
 \end{figure}

Since $\pi_{14,26}$ and $\pi_{14,27}$ do not generate minimal \iets, the inspection of the graph of Figure~\ref{clase14B} implies that $\{1,11,18,6,4,21,2\}\subset N^{(14)}$.
\textrm{
$\begin{maximaenv}
ctoN:union({1,11,18,6,4,21,2},ctoN)
\end{maximaenv}$
}
Then, $$\imt{ctoN}\subseteq N^{(14)}.$$

\begin{figure}[h]
$
 \begin{diagram}
\permuclaseinver{14}{2}{15}&\rCorresponds[l>=3em]^a&\permuclaseinver{14}{8}{13}\circlearrowleft^b\\
&&&&&\\
\dTo_b&&&&\\
\\
\circlearrowleft^a\permuclaseinver{14}{10}{7}&\rTo^b&\pi_{14,11}
 \end{diagram}
$
\hspace{1cm}
$
 \begin{diagram}
\circlearrowleft^a\permuclaseinver{14}{3}{16}&\rTo[l>=3em]^b&\pi_{14,6}
 \end{diagram}
$
\caption{Subgraphs in $\mathcal{C}_{14}$\label{clase14C}}
 \end{figure}

 Still Figure~\ref{clase14C} and the fact that $11$ and $6$ are in $N^{(14)}$ imply that $\{2,8,15,10,7,13,16,3\}\subset N^{(14)}$ and
\textrm{
$\begin{maximaenv}
ctoN:union({2,8,15,10,7,13,16,3},ctoN)
\end{maximaenv}$
}
\begin{equation}\label{E:minimalidad14}
\imt{ctoN}\subseteq N^{(14)}.
\end{equation}
Equation (\ref{E:minimalidad14}) proves the following proposition.

\begin{propo}\label{P:clase14}
 Let $T=(\lambda,\pi)$, with $\pi\in\mathcal{C}_{14}$. Then, $T$ is not minimal.
\end{propo}

\subsection{Class 20}
The eight permutation of this class are the following:

\begin{eqnarray*}
\mathcal{C}_{20}=\{    \pi_{20,1}=\imt{clase20[1]},\pi_{20,2}=\imt{clase20[2]},\pi_{20,3}=\imt{clase20[3]},\pi_{20,4}=\imt{clase20[4]},\\
    \pi_{20,5}=\imt{clase20[5]},\pi_{20,6}=\imt{clase20[6]},
\pi_{20,7}=\imt{clase20[7]},\pi_{20,8}=\imt{clase20[8]} \}.
 \end{eqnarray*}

 In this class, we have vertices which are conjugated or the inverse to some elements of other classes already studied. In particular:
 \begin{eqnarray*}
 \pi_{20,1}=\pi_{8,5}^{-1}=\imt{clase8[5]}^{-1},\quad
  \pi_{20,2}=\pi_{8,3}^h=\imt{clase8[3]}^h,\quad
  \pi_{20,3}=\pi_{8,4}^{-1}=\imt{clase8[4]}^{-1},\\
  \pi_{20,4}=\pi_{8,8}^{-1}=\imt{clase8[8]}^{-1},\quad
  \pi_{20,5}=\pi_{14,9}^h=\imt{clase14[9]}^h,\quad
  \pi_{20,6}=\pi_{8,2}^h=\imt{clase8[2]}^h,\\
  \pi_{20,7}=\pi_{8,1}^h=\imt{clase8[1]}^h,\quad
  \pi_{20,8}=\pi_{14,19}^h=\imt{clase14[19]}^h
 \end{eqnarray*}

 Therefore, \textcolor{black}{by Theorem~\ref{Th:hlamda} we obtain the non-minimality character in this class.}

\begin{propo}\label{P:clase20}
 Let $T=(\lambda,\pi),$ with $\pi\in\mathcal{C}_{20}$. Then, $T$ is not minimal.
\end{propo}

\subsection{Class 13}

\textcolor{black}{This class contains $36$ elements, namely:}
\begin{eqnarray*}
 \mathcal{C}_{13}=\{  \pi_{13,1}=\imt{clase13[1]},\pi_{13,2}=\imt{clase13[2]},\pi_{13,3}=\imt{clase13[3]},\pi_{13,4}=\imt{clase13[4]},\\
    \pi_{13,5}=\imt{clase13[5]},\pi_{13,6}=\imt{clase13[6]},
\pi_{13,7}=\imt{clase13[7]},\pi_{13,8}=\imt{clase13[8]},\\
\pi_{13,9}=\imt{clase13[9]},\pi_{13,10}=\imt{clase13[10]},\pi_{13,11}=\imt{clase13[11]},\pi_{13,12}=\imt{clase13[12]},\\
\pi_{13,13}=\imt{clase13[13]},\pi_{13,14}=\imt{clase13[14]},\pi_{13,15}=\imt{clase13[15]},\pi_{13,16}=\imt{clase13[16]},\\
\pi_{13,17}=\imt{clase13[17]},\pi_{13,18}=\imt{clase13[18]},
\pi_{13,19}=\imt{clase13[19]},\pi_{13,20}=\imt{clase13[20]},\\
\pi_{13,21}=\imt{clase13[21]},\pi_{13,22}=\imt{clase13[22]},\pi_{13,23}=\imt{clase13[23]},\pi_{13,24}=\imt{clase13[24]},\\
\pi_{13,25}=\imt{clase13[25]},\pi_{13,26}=\imt{clase13[26]},\pi_{13,27}=\imt{clase13[27]},\pi_{13,28}=\imt{clase13[28]},\\
\pi_{13,29}=\imt{clase13[29]},\pi_{13,30}=\imt{clase13[30]},\pi_{13,31}=\imt{clase13[31]},\pi_{13,32}=\imt{clase13[32]},\\
\pi_{13,33}=\imt{clase13[33]},\pi_{13,34}=\imt{clase13[34]},\pi_{13,35}=\imt{clase13[35]},\pi_{13,36}=\imt{clase13[36]}
\}.
 \end{eqnarray*}

Observe that:

 \begin{eqnarray*}
 \pi_{13,11}=\pi_{14,18}^h=\imt{clase14[18]}^h,
  \pi_{13,14}=\pi_{14,7}^h=\imt{clase14[7]}^h,
 \pi_{13,16}=\pi_{14,6}^h=\imt{clase14[6]}^h,\\
  \pi_{13,20}=\pi_{20,3}^h=\imt{clase20[3]}^h,
  \pi_{13,21}=\pi_{14,10}^h=\imt{clase14[10]}^h,
  \pi_{13,22}=\pi_{14,13}^h=\imt{clase14[13]}^h,\\
  \pi_{13,24}=\pi_{20,4}^h=\imt{clase20[4]}^h,
  \pi_{13,26}=\pi_{20,1}^h=\imt{clase20[1]}^h,
  \pi_{13,27}=\pi_{14,24}^h=\imt{clase14[24]}^h,\\
  \pi_{13,28}=\pi_{14,8}^h=\imt{clase14[8]}^h,
  \pi_{13,30}=\pi_{8,5}^h=\imt{clase8[5]}^h,
  \pi_{13,33}=\pi_{14,12}^h=\imt{clase14[12]}^h,\\
  \pi_{13,34}=\pi_{8,4}^h=\imt{clase8[4]}^h,
  \pi_{13,35}=\pi_{8,8}^h=\imt{clase8[8]}^h,
  \pi_{13,36}=\pi_{14,23}^h=\imt{clase14[23]}^h.
 \end{eqnarray*}

 \CuentasMaxima{ctoN:{11,14,16,20,21,22,24,26,27,28,30,33,34,35,36}}

 Let $N^{(13)}$ be the subset of naturals $j$ such that $\pi_{13,j}$ does not generate minimal \iets. \textcolor{black}{From the above 
observation and the previous study of Classes $\mathcal{C}_{14}$, $\mathcal{C}_{20}$, $\mathcal{C}_{8}$, we deduce that } 
 $$\mathcal{U}\textcolor{black}{:=}\imt{ctoN}\subseteq N^{(13)}.$$

 \begin{figure}[h]
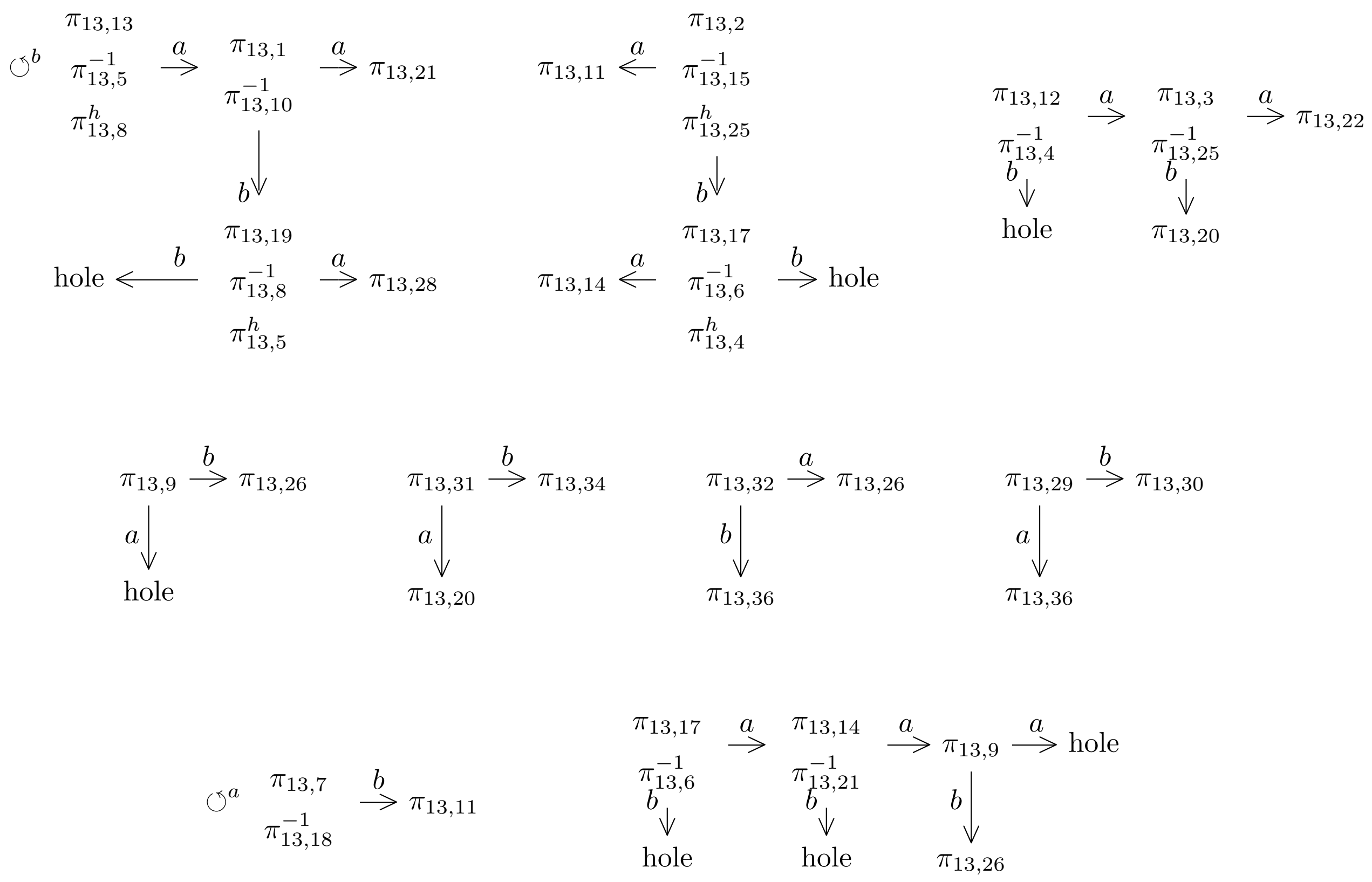

 \begin{tabular}{cp{.5cm}cp{.5cm}c}
  \begin{diagram}
\circlearrowleft^b\permuclaseinverh{13}{13}{5}{8}&\rTo^a& \permuclaseinver{13}{1}{10}&\rTo^a \pi_{13,21} \\
&& &&\\
&& \dTo^b\\
\hole&\lTo^b& \permuclaseinverh{13}{19}{8}{5}&\rTo^a \pi_{13,28}\\
\end{diagram}
&&
\begin{diagram}
\pi_{13,11}&\lTo^a&\permuclaseinverh{13}{2}{15}{25}\\
&\\&&\dTo^b\\
\pi_{13,14}&\lTo^a&\permuclaseinverh{13}{17}{6}{4}&\rTo^b &\hole\\
\end{diagram}
&&
\begin{diagram}
  \permuclaseinver{13}{12}{4}&\rTo^a&\permuclaseinver{13}{3}{25}&\rTo^a&\pi_{13,22}\\
  \dTo^b&&\dTo^b&&\\
  \hole&&\pi_{13,20}&&\\
 \end{diagram}
  \end{tabular}
\vspace*{1cm}

\begin{tabular}{cp{.5cm}cp{.5cm}cp{.5cm}c}
\begin{diagram}
\pi_{13,9}&\rTo^b&\pi_{13,26}\\
  \dTo^a&&\\
  \hole&&\\
 \end{diagram}
 &&
\begin{diagram}
\pi_{13,31}&\rTo^b&\pi_{13,34}\\
  \dTo^a&&\\
  \pi_{13,20}&&\\
 \end{diagram}
 &&
 \begin{diagram}
\pi_{13,32}&\rTo^a&\pi_{13,26}\\
  \dTo^b&&\\
  \pi_{13,36}&&\\
 \end{diagram}
 &&
 \begin{diagram}
 \pi_{13,29}&\rTo^b&\pi_{13,30}\\
  \dTo^a&&\\
  \pi_{13,36}&&\\
 \end{diagram}
 \end{tabular}
\vspace*{1cm}


\begin{tabular}{cp{1cm}c}
\begin{diagram}
\circlearrowleft^a\permuclaseinver{13}{7}{18}&\rTo^b&\pi_{13,11}
 \end{diagram}
 &&
 \begin{diagram}
\permuclaseinver{13}{17}{6}&\rTo^a&
\permuclaseinver{13}{14}{21}&\rTo^a&
\pi_{13,9}&\rTo^a&\hole \\
\dTo^b&&\dTo^b
&&\dTo^b
&&\\
\hole&&
\hole&&
\textcolor{black}{\pi_{13,26}}&&
 \end{diagram}
 \end{tabular}
\vspace*{1cm}
%
  \caption{\textcolor{black}{Some subgraphs for} Class $\mathcal{C}_{13}$\label{clase13}}
 \end{figure}

 \CuentasMaxima{ctoN:union(ctoN,{13,1,19,5,10,8,15,25,4,6,2,17,12,4,3,25,31,9,23,32,7,18,29,6,14,9})}

 \CuentasMaxima{nuevosenN:{13,1,19,5,10,8,15,25,4,6,2,17,12,4,3,25,31,9,23,32,7,18,29,6,14,9}}

Also, \textcolor{black}{repeating the same arguments as in precedent classes}, from the subgraphs represented in Figure~\ref{clase13} we deduce  that $$\mathcal{V}\textcolor{black}{:=}\imt{nuevosenN}\subseteq N^{(13)},$$ and then
$$\mathcal{U}\cup\mathcal{V}=\{n\}_{n=1}^{36}= N^{(13)}.$$
Therefore we have proved:

\begin{propo}\label{P:clase13}
 Let $T=(\lambda,\pi)$, with $\pi\in\mathcal{C}_{13}$. Then, $T$ is not minimal.
\end{propo}

\subsection{Orientable classes, 16 and 17} 
\begin{figure}[htb]
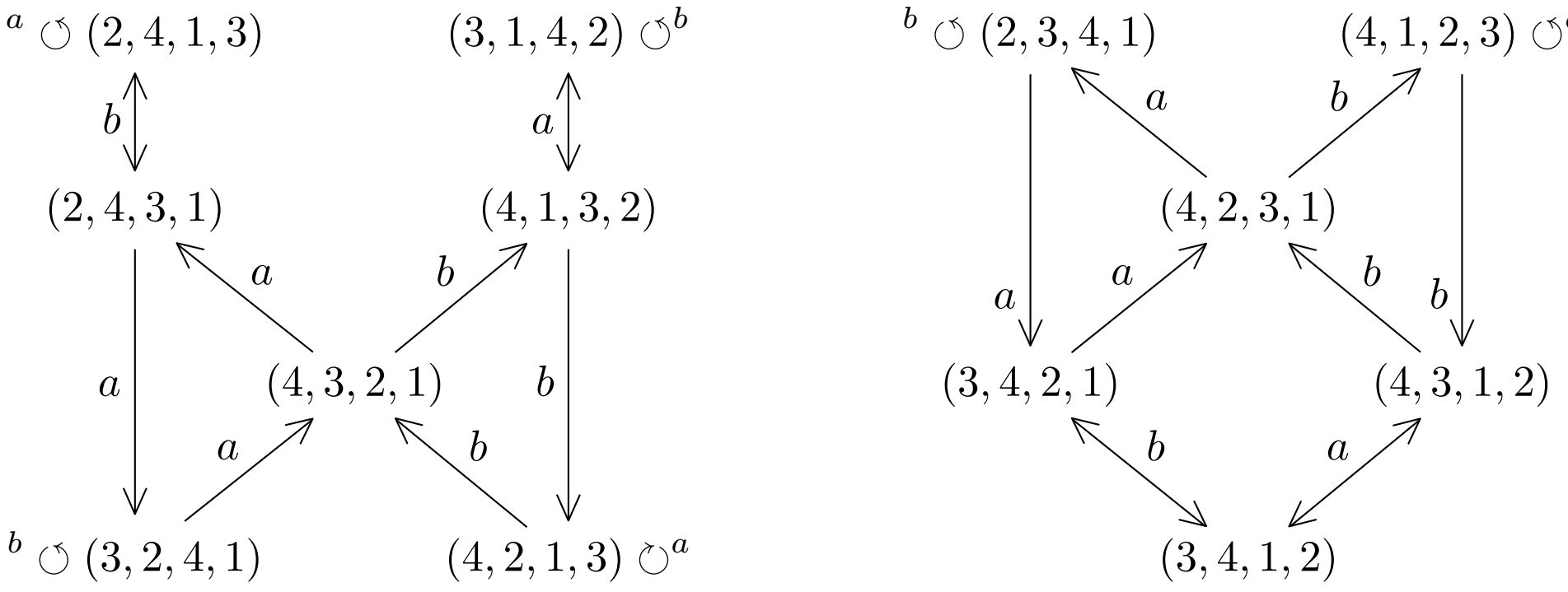

\begin{tabular}{cp{1cm}c}
$\begin{diagram}
^a\circlearrowleft(2,4,1,3)&&&&(3,1,4,2)\circlearrowleft^{b}\\
\dCorresponds^{b}&&&&\dCorresponds^{a}\\
(2,4,3,1)&&&&(4,1,3,2)\\
&\luTo{a}&&\ruTo^b&\\
\dTo^a&&(4,3,2,1)&&\dTo^b\\
&\ruTo^a&&\luTo^b&\\
^b\circlearrowleft(3,2,4,1)&&&&(4,2,1,3)\circlearrowright^a
\end{diagram}$
& &$\begin{diagram}
^b\circlearrowleft(2,3,4,1)&&&&(4,1,2,3)\circlearrowleft^{a}\\
&\luTo^{a}&&\ruTo^{b}&\\
&&(4,2,3,1)&&\\
\dTo^a&\ruTo^a&&\luTo^b&\dTo^b\\
(3,4,2,1)&&&&(4,3,1,2)\\
&\rdCorresponds^{b}&&\ldCorresponds^{a}&\\
&&(3,4,1,2)&&\\
\end{diagram}$
\end{tabular}
\caption{Representation of the two orientable Rauzy classes
\label{clasesorientables}}\end{figure}

There are two orientable classes that we have labelled by $\mathcal{C}_{16}$ and $\mathcal{C}_{17}$. We represent them in Figure~\ref{clasesorientables}.

It is known, see \cite{YoccozCurso}, the existence of minimal non uniquely ergodic \iets\  of the form $T=(\lambda,\pi),$ with
$\pi=(4,3,2,1)$. $T$ has two invariant measures, which is the maximum \textcolor{black}{number} of invariant measures for $n=4$ as proved by Veech in
\cite[Th. 2.12]{veech}.

\subsection{Class 19} This is the \textcolor{black}{largest} 
class and, \textcolor{black}{in contrast to the precedent classes of IETs with flips,} therein there are permutations with associated \iets\ being minimal. In particular, in \cite[Th. 6.1]{gutierrez4b} a minimal uniquely ergodic \iet\ $T=(\lambda,\pi)$, with $\pi=\pi_{19,23}=(-3,4,1,-2)$, was built. A generalization of this minimal \iet\ lead the present authors to build minimal uniquely ergodic \iets\ $S=(\mu,\tau)$, with $\tau=(-3,-4,-5,\dots,-n+1,n,1,-2)$, in \cite{linerosolerergodic}. We introduce the $59$ signed permutations in this class.
 
\begin{eqnarray*}
\mathcal{C}_{19}=\{\pi_{19,1}=\imt{clase19[1]},\pi_{19,2}=\imt{clase19[2]},\pi_{19,3}=\imt{clase19[3]},\pi_{19,4}=\imt{clase19[4]},\\
\pi_{19,5}=\imt{clase19[5]},\pi_{19,6}=\imt{clase19[6]},
\pi_{19,7}=\imt{clase19[7]},\pi_{19,8}=\imt{clase19[8]},\\
\pi_{19,9}=\imt{clase19[9]},\pi_{19,10}=\imt{clase19[10]},\pi_{19,11}=\imt{clase19[11]},\pi_{19,12}=\imt{clase19[12]},\\
\pi_{19,13}=\imt{clase19[13]},\pi_{19,14}=\imt{clase19[14]},\pi_{19,15}=\imt{clase19[15]},\pi_{19,16}=\imt{clase19[16]},\\
\pi_{19,17}=\imt{clase19[17]},\pi_{19,18}=\imt{clase19[18]},
\pi_{19,19}=\imt{clase19[19]},\pi_{19,20}=\imt{clase19[20]},\\
\pi_{19,21}=\imt{clase19[21]},\pi_{19,22}=\imt{clase19[22]},\pi_{19,23}=\imt{clase19[23]},\pi_{19,24}=\imt{clase19[24]},\\
\pi_{19,25}=\imt{clase19[25]},\pi_{19,26}=\imt{clase19[26]},\pi_{19,27}=\imt{clase19[27]},\pi_{19,28}=\imt{clase19[28]},\\
\pi_{19,29}=\imt{clase19[29]},\pi_{19,30}=\imt{clase19[30]},\pi_{19,31}=\imt{clase19[31]},\pi_{19,32}=\imt{clase19[32]},\\
\pi_{19,33}=\imt{clase19[33]},\pi_{19,34}=\imt{clase19[34]},\pi_{19,35}=\imt{clase19[35]},\pi_{19,36}=\imt{clase19[36]},\\
\pi_{19,37}=\imt{clase19[37]},\pi_{19,38}=\imt{clase19[38]},\pi_{19,39}=\imt{clase19[39]},\pi_{19,40}=\imt{clase19[40]},\\
\pi_{19,41}=\imt{clase19[41]},\pi_{19,42}=\imt{clase19[42]},\pi_{19,43}=\imt{clase19[43]},\pi_{19,44}=\imt{clase19[44]},\\
\pi_{19,45}=\imt{clase19[45]},\pi_{19,46}=\imt{clase19[46]},\pi_{19,47}=\imt{clase19[47]},\pi_{19,48}=\imt{clase19[48]},\\
\pi_{19,49}=\imt{clase19[49]},\pi_{19,50}=\imt{clase19[50]},\pi_{19,51}=\imt{clase19[51]},\pi_{19,52}=\imt{clase19[52]},\\
\pi_{19,53}=\imt{clase19[53]},\pi_{19,54}=\imt{clase19[54]},\pi_{19,55}=\imt{clase19[55]},\pi_{19,56}=\imt{clase19[56]},\\
\pi_{19,57}=\imt{clase19[57]},\pi_{19,58}=\imt{clase19[58]},\pi_{19,59}=\imt{clase19[59]}\}.
 \end{eqnarray*}

Permutation  $\pi_{19,50}=\imt{cl4[19][50]}$ is going to  play an important role in  the present section.
We will say that $v\in\{a,b\}^\N$ is a \emph{recurrent path} for a vertex $\pi$ if it contains an unbounded  quantity of operators $a$ and $b$ and there exists an increasing sequence $(m_k)_k \subseteq \N$ such that  $v_{m_k}(v_{m_{k}-1}(\dots v_1(\pi)))=\pi$ for all $k\in\N$; we will also say that \emph{$\pi$ admits a recurrent path} or that \emph{$v$ is a recurrent path for $\pi$}.
The  vertex $\pi$ is \textcolor{black}{said to be} a \emph{repeller} if it does not admit  recurrent paths.

A \emph{vertex $\pi$  is a weak  repeller} if for any recurrent path $v$ for $\pi$  and any $\lambda\in\Lambda^4$ satisfying $\mathcal{G}^{(\lambda,\pi)}=\mathcal{G}^{\pi,v}$ with $R^m((\lambda,\pi))=(\lambda^m,\pi_m)$, it  holds:
$|\lambda^m|$ does not go to zero as $m$ goes to $\infty$.

In order to simplify the study of this class, we introduce some important subgraphs in Figure~\ref{clase19A}.

\begin{figure}[h]
\scalebox{0.9}{
\begin{tabular}{cp{.5cm}c}
\begin{diagram}
&&&&&&\hole&&&\\
\dNoantecesor^b&&\dNoantecesor^a
&&\dNoantecesor^b&\ruTo^a&&&\\
\circlearrowleft^a\pi_{19,11}&\rTo^b&\pi_{19,15}
&\rTo^a&\pi_{19,8}&\rTo^b&\pi_{19,33}\circlearrowleft^a&\rTo^b&\pi_{19,35}\\
&&\dTo^b&&&&&&\\
&&\pi_{19,19}
&\rCorresponds[l>=3em]^a&\pi_{19,7}&\rTo^b&\pi_{19,25}\\
&&\dTo^b
&&\uNoantecesor^b&&\\
&&\hole
&&&&\\
\end{diagram}&&\\

&&\\
\begin{diagram}
&&&&&&\hole&&&\\
\dNoantecesor^a&&\dNoantecesor^b
&&\dNoantecesor^a&\ruTo^b&&&\\
\circlearrowleft^b\pi_{19,21}&\rTo^a&\pi_{19,13}
&\rTo^b&\pi_{19,34}&\rTo^a&\pi_{19,46}\circlearrowleft^b&\rTo^a&\pi_{19,57}\\
&&\dTo^a&&&&&&\\
&&\pi_{19,37}
&\rCorresponds[l>=3em]^b&\pi_{19,36}&\rTo^a&\pi_{19,44}\\
&&\dTo^a
&&\uNoantecesor^a&&\\   
&&\hole
&&&&\\
\end{diagram}
&&
\end{tabular}
}
\caption{Subgraphs with  repellers in class $\mathcal{C}_{19}$\label{clase19A}}

\end{figure}

\CuentasMaxima{R1:{11,15,8,33,19,7,21,13,34,46,37,36}}

\begin{lemma}\label{L:ctoT}
 Let $\mathcal{T}:=\imt{R1}$. If $j\in \mathcal{T}$, $\pi_{19,j}$ is a repeller.
\end{lemma}
\begin{proof}
 It is an easy consequence of the graphs in Figure~\ref{clase19A}.
\end{proof}

\begin{figure}[h]

\begin{diagram}
&&\puntoblanco&&\hole&&\puntoblanco&\rNoantecesor[l>=3em]^b&\pi_{19,3}\circlearrowleft^a\\
&&\dNoantecesor^b&&&\luTo^b&&&\dTo^b&&&&&&\\
\circlearrowleft^b\pi_{19,20}&\rCorresponds^a&\pi_{19,9}
&\rTo^b&\pi_{19,26}\circlearrowleft^a&\rTo^b&\pi_{19,27}&\rTo^a&\pi_{19,14}\\
&&&&&&\uTo^a&\ldTo^a&\dTo^b\\
\puntoblanco\puntoblanco&\rNoantecesor[l>=3em]^a&\pi_{19,47}\circlearrowleft^b&&\puntoblanco&\rNoantecesor[l>=3em]^b&\pi_{19,1}&\rTo^b&&&\hole\\
&&\dTo^a&&&&&&&&\uTo^a\\
\circlearrowleft^a\pi_{19,55}&\rTo^b&\pi_{19,59}&\rTo^a&\pi_{19,32}&\rCorresponds^b&\pi_{19,31}&\rTo^a&\pi_{19,17}&\rTo^a&\pi_{19,6}&\rTo^b&\pi_{19,41}\\
&\luTo^b&\dTo^b&&\dTo^a&&\uNoantecesor^a&&\dTo^b&&\uNoantecesor^b\\
\hole&\lTo^a&\pi_{19,53}&&\hole&&\puntoblanco&&\hole&&\puntoblanco\\
&&\uNoantecesor^a&&&&&&&&&\\
\negro
\end{diagram}

\vspace*{1cm}
\begin{diagram}
&&\puntoblanco&&\hole&&\puntoblanco&\rNoantecesor[l>=3em]^a&\pi_{19,43}\circlearrowleft^b\\
&&\dNoantecesor^a&&&\luTo^a&&&\dTo^a&&&&&&\\
\circlearrowleft^a\pi_{19,38}&\rCorresponds^b&\pi_{19,39}
&\rTo^a&\pi_{19,45}\circlearrowleft^b&\rTo^a&\pi_{19,54}&\rTo^b&\pi_{19,48}\\
&&&&&&\uTo^b&\ldTo^b&\dTo^a\\
\puntoblanco\puntoblanco&\rNoantecesor[l>=3em]^b&\pi_{19,12}\circlearrowleft^a&&\puntoblanco&\rNoantecesor[l>=3em]^a&\pi_{19,51}&\rTo^a&&&\hole\\
&&\dTo^b&&&&&&&&\uTo^b\\
\circlearrowleft^b\pi_{19,28}&\rTo^a&\pi_{19,16}&\rTo^b&\pi_{19,18}&\rCorresponds^a&\pi_{19,5}&\rTo^b&\pi_{19,24}&\rTo^b&\pi_{19,29}&\rTo^a&\pi_{19,23}\\
&\luTo^a&\dTo^a&&\dTo^b&&\uNoantecesor^b&&\dTo^a&&\uNoantecesor^a\\
\hole&\lTo^b&\pi_{19,4}&&\hole&&\puntoblanco&&\hole&&\puntoblanco\\
&&\uNoantecesor^b&&&&&&&&&\\
\negro
\end{diagram}

\caption{Subgraphs with  repellers in class $\mathcal{C}_{19}$\label{clase19B}}
\end{figure}


\CuentasMaxima{R2:{12,47,26,27,14,3,1,32,31,17,6,45,54,48,43,51,18,5,24,29}}


\begin{lemma}\label{L:ctoS}
 Let $$j\in \mathcal{S}:=\imt{R2}.$$  Then, $\pi_{19,j}$ is a repeller.
\end{lemma}
\begin{proof}
From a  simple inspection of Figure~\ref{clase19B}, we have for instance that $\pi_{19,12}$ and $\pi_{19,47}$ are repellers. The same occurs for the remaining permutations by an easy analysis of Figure~\ref{clase19B}.
\end{proof}


\CuentasMaxima{R3:{59,53,55,16,4,28,9,20,38,39}}
\begin{lemma}\label{L:ctoF}
 Let $j\in\mathcal{F}:=\imt{R3}$. Then, $\pi_{19,j}$ is a weak repeller.
\end{lemma}
\begin{proof}
 Recall that $\pi_{19,55}=\imt{clase19[55]}$,
 $\pi_{19,59}=\imt{clase19[59]}$,
 $\pi_{19,53}=\imt{clase19[53]}$,
 $\pi_{19,28}=\imt{clase19[28]}$,
 $\pi_{19,16}=\imt{clase19[16]}$, and
 $\pi_{19,4}=\imt{clase19[4]}$. 
 
 Observe that, according to  
 Figure~\ref{clase19B}, any recurrent path begining in $\pi_{19,55}$ must have alternatively $3r$ times the operator $b$ and $s$ times the operator $a$,  with $r,s$ nonnegative integers. Assume that we have an initial \iet\ $T=(\lambda,\pi_{19,55})$ such that
 $\mathcal{G}^{\pi,v}=\mathcal{G}^T$; then, by Claim~\ref{Cl:ExchangeIntervalsBya},  after applying $3r$ times $b$  and $s$ times $a$ we obtain $R^{3r+s}(T)=(\lambda^{3r+s},\pi_{3r+s})=(\lambda^{3r+s},\pi_{19,55})$, with $\lambda^{3r+s}=(\lambda_1^{(3r)},\lambda_2,\lambda_3,\lambda_4^{(s)})$, where $\lambda_{j}^{(f)}$ is a length obtained from $\lambda_j$ shortening it according to the Rauzy induction; therefore, $|\lambda^m|>\lambda_2+\lambda_3$ and $|\lambda^m|$ does not goes to zero as $m$ goes to $\infty$. Then,  $\pi_{19,55}$ is a weak  repeller and  similar reasonings applies for $\pi_{19,59}$ and $\pi_{19,53}$.

 Note now that, by Figure~\ref{clase19B}, any recurrent path begining in $\pi_{19,28}$ must have alternatively $3r$ times the operator $a$ 
and $s$ times the operator $b$,  with $r,s$ nonnegative integers. Assume that we have an initial \iet\  $T=(\lambda,\pi_{19,28})$ such that
 $\mathcal{G}^{\pi,v}=\mathcal{G}^T$; again by Claim~\ref{Cl:ExchangeIntervalsBya}, after applying $3r$ times $a$  and $s$ times $b$ we obtain $R^{3r+s}(T)=(\lambda^{3r+s},\pi_{3r+s})=(\lambda^{3r+s},\pi_{19,28})$, with $\lambda^{3r+s}=(\lambda_1,\lambda_2,\lambda_3^{(s)},\lambda_4^{(3r)})$; then $|\lambda^m|>\lambda_1+\lambda_2$ and $|\lambda^m|$ does not goes to zero as $m$ goes to $\infty$. Therefore,  $\pi_{19,28}$ is a weak  repeller and  similar reasonings applies for $\pi_{19,4}$ and $\pi_{19,16}$.

\textcolor{black}{For the couple $ \pi_{19,9}=(-4, 1, 3, 2)$, $\pi_{19,20}=(-3, 1, 4, 2)$, the reasoning is similar. Now, by Claim~\ref{Cl:ExchangeIntervalsBya}, it is easily seen again that, for instance, the first interval of the IET is never shortened when applying the Rauzy transformation.} 
\textcolor{black}{The same applies for the couple  $\pi_{19,38}=(2, 4, -1, 3)$, $\pi_{19,39}=(2, 4, 3, -1)$.} 
 \end{proof}
\CuentasMaxima{admisible:{23,44,22,2,56,10,40,58,30,35,52,57,49,25,50,42,41}}


\textcolor{black}{As a consequence of the previous lemmas, we can narrow down on those possible permutations 
that may appear as combinatorial data of a minimal IET.}

\begin{coro}\label{C:previo}
 Let $\textcolor{black}{T}=(\pi,\lambda)$, $\pi\in\mathcal{C}_{19}$, be a minimal \iet\ and write 
$\textcolor{black}{R^j(T)}=(\lambda_j,\pi_{19,\rho(j)})$, $1\leq \rho(j)\leq 59$, for any $j\in\N$. 
 Then, there exists a natural number $N_0$ such that  if $j\geq N_0$ then $\rho(j)\not\in\mathcal{F}\cup\mathcal{T}\cup\mathcal{S}$, or equivalently
$$\rho(j)\in\imt{admisible}.$$ 
\end{coro}
\begin{proof}
 By Lemmas~\ref{L:ctoF}, \ref{L:ctoT} and \ref{L:ctoS} any permutation $\pi=\pi_{19,j}$ with $j\in\mathcal{F}\cup\mathcal{T}\cup\mathcal{S}$ is either a  repeller or a weak  repeller. Assume that $T=(\lambda,\pi)$ is minimal; then, by Proposition~\ref{P:minimalidadimplicainducible}, the Rauzy operator can be applied to $T$ infinitely many times,  $R^m(T)=(\lambda^m,\pi_m)$ for any $m\in\N$, with $\lim_{m\to\infty}|\lambda^m|=0$. Observe that there exits an integer $N_0$ satisfying the properties of the statement, otherwise we would contradict that $\pi$ is either a  repeller or a  weak  repeller.
\end{proof}

\begin{figure}[h]
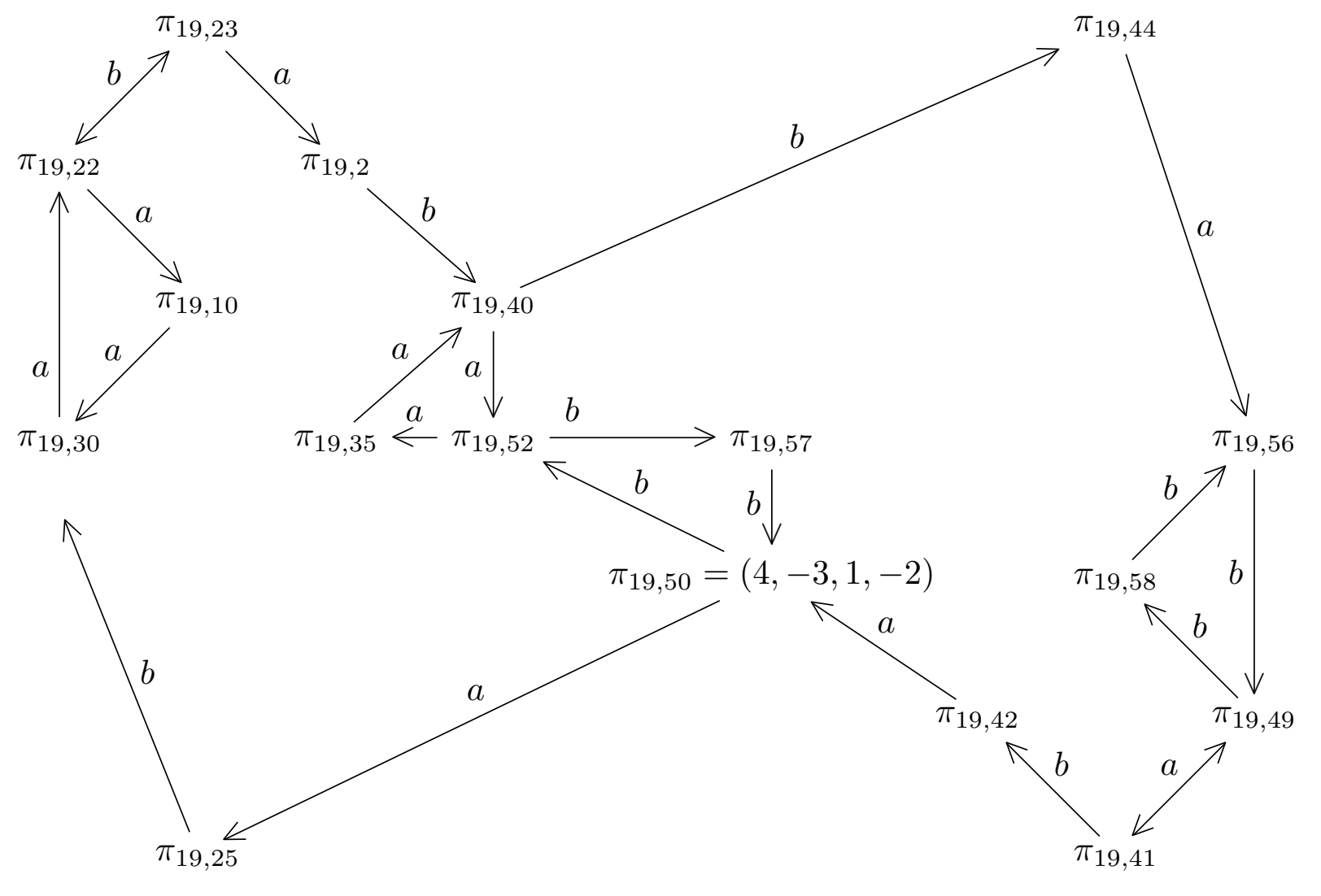


\scalebox{.8}{
\begin{diagram}
  &&\pi_{19,23}&&&&&&&&&\pi_{19,44}&&&\\
&\ruCorresponds^b&&\rdTo^a&&&&&&&\ruTo(5,4)^b&&\rdTo(2,6)^a&&\\
  \pi_{19,22}&&&&\pi_{19,2}&&&&&&&&&&\\
  &\rdTo^a&&&&\rdTo^b&&&&&&&&&\\
  &&\pi_{19,10}&&&&\pi_{19,40}&&&&
  &&&&\\
  \uTo^a&\ldTo^a&&&&\ruTo^a&\dTo^a&&&&&&&&\\
  \pi_{19,30}&&&&\pi_{19,35}&\lTo^a&\pi_{19,52}&\rTo^b&\pi_{19,57}&&&&&\pi_{19,56}&\\
&&&&&&&\luTo^b&\dTo^b&&&&\ruTo^b&&\\
  &\luTo(2,5)^b&&&
  &&&&\pi_{19,50}=\imt{clase19[50]}&&&\pi_{19,58}&&\dTo^b&\\
  &&&&&&&\ldTo(6,4)^a&&\luTo(1,2)^a&&&\luTo^b&\\
  &&&&&&&&&\pi_{19,42}&&&&\pi_{19,49}\\
  &&&&&&&&&&\luTo^b&&\ldCorresponds^a&&&\\
  &&\pi_{19,25}&&&&&&&&&\pi_{19,41}&&&\\
\end{diagram}
}

\caption{Non repeller vertices in $\mathcal{C}_{19}$ \label{nonrepellor}}
\end{figure}


 The objective now is to analyze the Rauzy subgraph of Figure~\ref{nonrepellor}, whose 17 vertices are the vertices of the set  
$\mathcal{C}_{19}\menos\left( \mathcal{S}\cup\mathcal{T}\cup\mathcal{F}\right)$, see Corollary~\ref{C:previo}. For this purpose, 
the following technical results will be useful.

\begin{lemma}\label{L:nociclosmenores}
 Let $T=(\lambda,\pi)$ be a minimal \iet\ with $R^m(T)=(\lambda^m,\pi_m)$. 
 \begin{enumerate}
  \item If $\pi=\pi_{19,30}=\imt{clase19[30]}$, there exists $k$ such that $\pi_k=\pi_{19,2}$.
  \item If $\pi=\pi_{19,56}=\imt{clase19[56]}$, there exists $k$ such that $\pi_k=\pi_{19,42}$.
 \end{enumerate}
\end{lemma}
\begin{proof}
 Assume first that $\pi=\pi_{19,30}$ and that the first statement does not hold, then there exists an increasing sequence of integers $(n_k)_k$ such that $\pi_{n_k}=\pi_{19,30}$ and the  path associated to the Rauzy induction of $T$ must be the concatenation of paths of the form $v_{r,s}=a^{3r+1}b^{2s}a^2$, with $r,s$ nonnegative (see Proposition~\ref{P:minimalidadimplicainducible}). Observe that  $v_{r,s}(\pi)=\pi$ and by the use of Claim~\ref{Cl:ExchangeIntervalsBya} and the identifications $\pi_{19,22}=\imt{clase19[22]}$, $\pi_{19,23}=\imt{clase19[23]}$, $\pi_{19,10}=\imt{clase19[10]}$, $\pi_{19,30}=\imt{clase19[30]}$,  we obtain 
 $\lambda^{3r+2s+3}=(\lambda_1,\lambda_2^{(2s)},\lambda_3,\lambda_4^{(3r+3)})$, recall that $\lambda_j^{(i)}$ means that the length $\lambda_j$ has been shortened  $i$ times according to the Rauzy operator. Therefore $|\lambda^m|>\lambda_1+\lambda_3$ and then we contradict the minimality of $T$ because $\lim_{m\to \infty}|\lambda^m|\not=0$. Therefore there exists an integer $k$  such that  $\pi_k=\pi_{19,2}$.
 
Suppose now that  $\pi=\pi_{19,56}$ and that the second statement does not hold. In this situation, there exists an increasing sequence of integers $(n_k)_k$ such that $\pi_{n_k}=\pi_{19,56}$ and the  path associated to the Rauzy induction of $T$ must be the concatenation of paths of the form $v_{r,s}=b^{3r+1}a^{2s}b^2$, with $r,s$ nonnegative, so  $v_{r,s}(\pi)=\pi$. Taking into account that
 $\pi_{19,56}=\imt{clase19[56]}$, $\pi_{19,49}=\imt{clase19[49]}$, $\pi_{19,58}=\imt{clase19[58]}$, $\pi_{19,41}=\imt{clase19[41]}$ and Claim~\ref{Cl:ExchangeIntervalsBya}, we obtain 
 $\lambda^{3r+2s+3}=(\lambda_1^{(3r+3)},\lambda_2,\textcolor{black}{\lambda_3^{(2s)},\lambda_4})$.  Therefore, 
 $| \lambda^m |>\lambda_2+\textcolor{black}{\lambda_4}$ and then we contradict the minimality of $T$ because $\lim_{m\to \infty}|\lambda^m|\not=0$. Therefore, there exists an integer $k$  such that  $\pi_k=\pi_{19,42}$.
 
\end{proof}

\begin{lemma}\label{L:nociclosdobles}
 Let $T=(\lambda,\pi_{19,50})$ be a minimal \iet\ with $R^m(T)=(\lambda^m,\pi_m)$ then there exists $n$ such that 
 $\pi_n\in\{\pi_{19,44},\pi_{19,25}\}$.
 \end{lemma}
\begin{proof}
Assume the contrary, then there exists a strictly increasing sequence $(m_j)$ such that $\pi_{m_j}=\pi_{19,52}$. Consider now the minimal \iet\  $S=R^{m_1}(T)=(\mu,\pi_{19,52})$; then, by Proposition~\ref{P:minimalidadimplicainducible}, the Rauzy induction will follow a concatenation of paths of the form 
$v_{r,s}=b^{3r}a^{3s}$ or $w_{r,s}=a^{3r}b^{3s}$ for adequate $r>0$ and $s>0$. The use  of these paths maps  
$S=(\mu,\pi_{19,52})$ onto $S^{(3r+3s)}=(\mu^{3r+3s},\pi_{19,52})$ by the Rauzy operator.
Take into account Claim~\ref{Cl:ExchangeIntervalsBya} and  the values  
$\pi_{19,50}=\imt{clase19[50]}$, $\pi_{19,52}=\imt{clase19[52]}$,  
$\pi_{19,57}=\imt{clase19[57]}$ and observe  that 
$\mu^{3r+3s}=(\mu_1^{(3r)},\mu_2,\mu_3,\mu_4^{(3s)})$ when applying 
$v_{r,s}$, and 
$\mu^{3r+3s}=(\mu_1^{(3s)},\mu_2,\mu_3,\mu_4^{(3r)})$ when using  
$w_{r,s}$. In any case, the second and third intervals associated to $S$ are not shortened by the operator $R$, a contradiction with Proposition~\ref{P:minimalidadimplicainducible}.
\end{proof}

Consider $\pi\in\mathcal{C}_{19}$ and  a minimal  \iet\ $T=(\lambda,\pi)$ (infinitely many times inducible by Proposition~\ref{P:minimalidadimplicainducible}), with $R^m(T)=(\lambda^m,\pi_m)$. 
\textcolor{black}{Realize that the Rauzy process cannot end at the infinite concatenation of the loop $\pi_{19,40}, \pi_{19,52}, \pi_{19,35},$ otherwise the type of the induced IETs would be eventually constant, in contradiction with Proposition~\ref{P:minimalidadimplicainducible}. The same happens for the loop 
$\pi_{19,22}, \pi_{19,10}, \pi_{19,30},$ and the cycles $\pi_{19,52}, \pi_{19,57}, \pi_{19,50},$ and $\pi_{19,58}, \pi_{19,56}, \pi_{19,49}.$}
Thus, from   Figure~\ref{nonrepellor}, Corollary~\ref{C:previo} and Lemma~\ref{L:nociclosmenores}, there exists an increasing sequence $(m_k)_k$ such that $\pi_{m_k}=\pi_{19,50}=\imt{clase19[50]}$. Moreover, by Lemma~\ref{L:nociclosdobles}, we can chose  $(m_k)$ such that  the path of operators by the Rauzy induction from $(\lambda^{m_k},\pi_{m_k})$ is a  \textcolor{black}{ loop -cyclic path- in the following types:}

%
%
%
%
%
%
%

\begin{eqnarray}
	v&=&abaa^{3r}bb^{2u}abaa^{3s}b^2,\quad \nonumber\\
	w&=&abaa^{3r}bb^{2u}aba^{3s} babb^{3m}aa^{2q}ba,  \quad \nonumber\\
	z&=&b^{3r+1}a^{3u+2}bab^{3s+1}a^{2q+1}ba,\label{Eq:caminos}
\end{eqnarray}
%
with associated Rauzy matrices $P(v)=M(\mathcal{G}^{\pi_{19,50},v})$, 
$P(w)=M(\mathcal{G}^{\pi_{19,50},w})$, and $P(z)=M(\mathcal{G}^{\pi_{19,50},z})$. 
\textcolor{black}{However, in our analysis, we will compute the Rauzy matrix associated to: $P(v)$; the matrix $P(\widetilde{w})$ of the (non-cyclic) path 	
\begin{equation}\label{Eq:wmodify}
\widetilde{w}:=abaa^{3r}bb^{2u}aba^{3s}bab
\end{equation} 
which connects  $\pi_{19,50}$ with $\pi_{19,49}$; and the matrix $P(\widetilde{z})$ corresponding to the path 
\begin{equation}\label{Eq:zmodify}
\widetilde{z}=zb.
\end{equation}
Notice that, in the Rauzy process applied to a minimal $4$-IET with combinatorial data $\pi_{19,50}$, at least one of type of the paths $v$, $\widetilde{w}$ or $\widetilde{z}$ must appear infinitely many times.} 
In order to build the matrices $P(v), P(\widetilde{w}), P(\widetilde{z})$, we need the following lemma that the reader can easily prove by applying (\ref{E:matrizAB}) and (\ref{E:matrizv}). \textcolor{black}{Realize that, in the following submatrices, we abbreviate the writing; for instance, $M_{aba}(\pi_{19,50})$ is meant the product $M_{a}(\pi_{19,50})\cdot M_{b}(\pi_{19,25})\cdot M_{a}(\pi_{19,30})$. }

\begin{lemma} Let $M_1=M_{aba}(\pi_{19,50})$,
$M_2=M_{a^{3r}}(\pi_{19,22})$, 
$M_3=M_{b}(\pi_{19,22})$,
$M_4=M_{b^{2r}}(\pi_{19,23})$,
$M_5=M_{aba}(\pi_{19,23})$, 
$M_6=M_{a^{3r}}(\pi_{19,52})$, 
$M_7=M_{b^{2}}(\pi_{19,52})$, 
$M_8=M_{ab^2ab}(\pi_{19,23})$, 
$M_9=M_{b^{3s+1}}(\pi_{19,50})$,
$M_{10}=M_{a^{3p+2}}(\pi_{19,52})$,
$M_{11}=M_{bab}(\pi_{19,40})$, 
$M_{12}=M_{b^{3u}}(\pi_{19,49})$, 
$M_{13}=M_{a^{2v+1}}(\pi_{19,49})$, 
$M_{14}=M_{bab}(\pi_{19,41})$. Then the values in Table~\ref{matricesm1a14}  hold.

\begin{table}[htb]
 
\begin{tabular}{ccc}

$
M_1=
\begin{pmatrix}
 1&0&0&0\cr 0&1&1&0\cr 0&0&1&1\cr 1&1&0&0
\end{pmatrix}
$

&

$
M_2=
\begin{pmatrix}
 1&0&0&0\cr 0&1&0&0\cr 0&0&1&0\cr r&r&r&1
\end{pmatrix}$

&

$
M_3=
\begin{pmatrix}
 1&0&0&0\cr 0&1&1&0\cr 0&0&0&1\cr 0&0&1&0
\end{pmatrix}$
\\
&&\\

$
M_4=\begin{pmatrix}
 1&0&0&0\cr 0&1&r&r\cr 0&0&1&0\cr 0&0&0&1
\end{pmatrix}$

&

$
M_5=\begin{pmatrix}
 1&1&0&0\cr 0&0&1&0\cr 0&1&0&1\cr 1&0&1&0
\end{pmatrix}$

&

$
M_6=\begin{pmatrix}
 1&0&0&0\cr 0&1&0&0\cr 0&0&1&0\cr r&r&r&1
\end{pmatrix}$
\\
&&
\\

$
M_7=
\begin{pmatrix}
 1&1&1&0\cr 0&0&0&1\cr 0&1&0&0\cr 0&0&1&0
\end{pmatrix}.$

&
$
  M_8=\begin{pmatrix}
 1&1&1&1\cr 0&1&0&1\cr 0&0&1&0\cr 1&2&0&1
\end{pmatrix}.$

&
 
$
M_9=\begin{pmatrix} 1&s+1&s&s
 \cr 0&0&1&0\cr 0&0&0&1\cr 0&1&0&0\cr
 \end{pmatrix} $

 \\
 &&
 \\

$
M_{10}=
\begin{pmatrix} 1&0&0&0
 \cr 0&1&0&0\cr 0&0&1&0\cr p+1&p&p+1&1\cr
 \end{pmatrix}
 $
& 
$
M_{11}=
\begin{pmatrix}1&1&0&0
 \cr 0&0&1&1\cr 0&1&0&1\cr 0&0&1&0\cr
 \end{pmatrix}
$
 &

$
M_{12}=
\begin{pmatrix} 1&u&u&u
 \cr 0&1&0&0\cr 0&0&1&0\cr 0&0&0&1\cr
 \end{pmatrix}$

 \\
 &&
\\

$
M_{13}=
\begin{pmatrix}
    1&0&0&0
 \cr 0& 1&0&0\cr 0&0&1&0\cr
 v+1&v&0&1\cr
 \end{pmatrix}
$

&

$
M_{14}=
\begin{pmatrix} 1&1&0&0
 \cr 0&0&1&1\cr 0&1&0&1\cr 0&0&1&0\cr
 \end{pmatrix}$

&
\\
&&\\
\end{tabular}

\caption{\label{matricesm1a14} Matrices needed to compute $P(v)$, $P(\widetilde{w})$ and $P(\widetilde{z})$.}
\end{table}

\end{lemma}

\begin{coro} \label{C:PvPw}

The values of $P(v)$, $P(\widetilde{w})$ and $P(\widetilde{z})$  are given respectively in Tables \ref{T:Pv},
\ref{T:Pw} and \ref{T:Pz}.

\begin{table}
 
 \scalebox{1}{
\begin{minipage}{\textwidth}
\begin{flushleft}
$
P(v)=
\left(
\begin{array}{ccc}
1&1&\vdots\\ 
s\,u+u+s+1&2\,s\,u+2\,u+2\,s+3&\vdots\\ 
 r\,s\,u+r\,u+r\,s+s+2\,r+1&2\,r\,s\,u+2\,r\,u+2\,r\,s+2
 \,s+4\,r+2&\vdots\\ 
 s\,u+u+s+1&2\,s\,u+2\,u+2\,s+2&\vdots
 \end{array}
 \right.
$
\end{flushleft}
\begin{flushright}
 $
 \left.
\begin{array}{ccc}
\vdots&1&1\\ 
\vdots&s\,u+2\,u+s+2&s
 \,u+u+s+1\\ 
 \vdots &r\,s\,u+2\,r\,u+r\,s+s+3\,r+2&r\,s\,u+r\,u+r\,s+s+2\,r+1
 \\ 
 \vdots&s\,u+2\,u+s+2&s\,u+u+s+2
 \end{array}
 \right),
$
\end{flushright}
\end{minipage}
 }
 
 \vspace*{.5cm}
 \caption{\label{T:Pv} Matrix $P(v)$}
 \end{table}
 
 
%

 \begin{table}
 
 \scalebox{1}{
\begin{minipage}{\textwidth}

\begin{flushleft}
$P(\widetilde{w})=
\left(\begin{array}{ccc} 
       1&1&\vdots \cr 
       s\,u+u+s+1&2\,s\,u+2\,u+2\,s+3&\vdots\cr 
       r \,s\,u+r\,u+r\,s+s+2\,r+1&2\,r\,s\,u+2\,r\,u+2\,r\,s+2\,s+4\,r+2&\vdots\cr 
 s\,u+ u+s+1&2\,s\,u+2\,u+2\,s+2&\vdots\cr 
      \end{array}\right.$ 
\end{flushleft}
\begin{flushright}
$\left.
\begin{array}{ccc} 
      \vdots&1&1 \cr 
      \vdots&s\,u+u+s+1&2\,s\,u+u+2\,s+2\cr 
      \vdots&r\,
 s\,u+r\,u+r\,s+s+2\,r+1&2\,r\,s\,u+r\,u+2\,r\,s+2\,s+3\,r+1\cr 
 \vdots&s\,u+u+s+2&2\,s\,u+u+2\,s+2\cr 
      \end{array}
      \right)$ 
\end{flushright}


\end{minipage}
 }
 \vspace{.5cm}
 
 \caption{\label{T:Pw} Matrix $P(\widetilde{w})$}
 \end{table}

 \begin{table}
 \scalebox{1}{
\begin{minipage}{\textwidth}
\begin{flushleft}
 $P(\widetilde{z})=
\left(
\begin{array}{ccc}
 A+p\,s\,u+s\,u+
 u+3\,p\,s+4\,s+2&A+2\,p\,
 s\,u+2\,s\,u+2\,u+4\,p\,s+6\,s+3&\vdots
 \cr 
 v+1&v+1&\vdots
  \cr 
  B+p\,u+u+3\,p+2&B+2\,p\,u+2\,u+4\,p+3&\vdots
 \cr 
 v+1&v+2&\vdots
 \cr 
\end{array}
\right.
$
\end{flushleft}
\begin{flushright}
$
 \left.
\begin{array}{ccc}
\vdots&A+2\,p\,s\,u+2\,s\,u+2\,u+4\,p\,s+6\,s+2&
A+2\,p\,s\,u+2\,s\,u+2\,u+3\,p\,s+5\,s+2
 \cr 
 \vdots
 &v+2&v+1
 \cr 
\vdots
&B+2\,p\,u+2\,u+4
 \,p+3&B+2\,p\,u+2\,u+3\,p+3
 \cr \vdots&v+1&v+1
 \cr
 \end{array}
\right)
$
\end{flushright}

\end{minipage}
 }
 with $A=p\,s\,u\,v+s\,u\,v+u\,v+2\,p\,s\,v+3\,s\,v+v$ and
$B=p\,u\,v+u\,v+2\,p\,v+v$. 
\vspace{.5cm}

 \caption{\label{T:Pz} Matrix $P(\widetilde{z})$}
 \end{table}

%
%

\end{coro}

\begin{proof}
 By means of the previous lemma and Figure~\ref{nonrepellor} we obtain:
 \begin{eqnarray*}
  P(v)&=&M_{abaa^{3r}bb^{2u}abaa^{3s}b^2}(\pi_{19,50})\\
  &=&
 M_{aba}(\pi_{19,50})
 M_{a^{3r}}(\pi_{19,22})
 M_{b}(\pi_{19,22})
 M_{b^{2u}}(\pi_{19,23})
 M_{aba}(\pi_{19,23})
 M_{a^{3s}}(\pi_{19,52})
 M_{bb}(\pi_{19,52}),\\ 
 P(\widetilde{w})&=& 
 M_{
 abaa^{3r}bb^{2u}aba^{3s}bab}(\pi_{19,50})\\
 &=&
 M_{aba}(\pi_{19,50})
 M_{a^{3r}}(\pi_{19,22})
 M_{b}(\pi_{19,22})
 M_{b^{2u}}(\pi_{19,23})
 M_{ab}(\pi_{19,23})
 M_{a^{3s}}(\pi_{19,40})
 M_{bab}(\pi_{19,40})
 \end{eqnarray*}
and 
\begin{eqnarray*}
 P(\widetilde{z})&=&
 M_{b^{3s+1}a^{3p+2}babb^{3k}a^{2l+1}bab}(\pi_{19,50})\\
&=&
 M_{b^{3s+1}}(\pi_{19,50})
 M_{a^{3p+2}}(\pi_{19,52})
 M_{bab}(\pi_{19,40})
 M_{b^{3u}}(\pi_{19,49})
 M_{a^{2v+1}}(\pi_{19,49})
 M_{bab}(\pi_{19,41}).
 \end{eqnarray*}
The proof follows by computing these products of matrices.
\end{proof}

 Observe that, in the matrices $P(v), P(\widetilde{w}), P(\widetilde{z})$ of Corollary~\ref{C:PvPw}, we can take, for any $1\leq i\leq 4$, the quantities of Proposition~\ref{Lipschitz3} as follows. Let $P$ be one of the previous matrices and write  
 $\delta_1^i=\min\{P_{i,j}:1\leq j\leq 4\}-\frac1{10}$ and $\delta_2^i=5\delta_1^i$. Then we have 
 $$\delta_1^i< P(x)_{i,j} < \delta_2^i=5\delta_1^i \quad \textrm{ for any }\quad  1\leq i,j\leq 4,\quad x\in\{v,\widetilde{w}, \widetilde{z}\}.$$
 
 \begin{lemma}\label{L:lipschitz}
  Assume the notation in Proposition~\ref{Lipschitz3} and let $x\in\{v,\textcolor{black}{\widetilde{w}, \widetilde{z}}\}$. Then, the linear map defined by  $P(x)$ satisfies (\ref{eqProyective}) and we have $$\delta(x)=\sqrt{\sup_{1\leq i\leq n}\left\{\frac{\delta_2^i}{\delta_1^i}\right\}}=\sqrt{5}.$$
Consequently, if $\lambda,\mu\in\R_+^4$, $$\ddp(P(x)\lambda,P(x)\mu )\leq \frac{\sqrt{5}-\frac1{\sqrt{5}}}{\sqrt{5}+\frac1{\sqrt{5}}}\ddp(\lambda,\mu)=\frac{2}{3}\ddp(\lambda,\mu).$$
 \end{lemma}
\begin{proof}
 $P(x)$ satisfies (\ref{eqProyective}) because it is positive; the other statements are  direct consequences of Proposition~\ref{Lipschitz3}.
\end{proof}

\begin{theorem}\label{T:clase19NoExisteMinimalNoUnicamenteErgodica} Let $\pi\in\mathcal{C}_{19}$ and let $\lambda\in\Lambda^4$ such that 
 $T=(\lambda,\pi)$ is minimal. Then, $T$ is uniquely ergodic. 
\end{theorem}

\begin{proof} 
We begin by claiming that it suffices to prove the result for $\pi=\pi_{19,50}$. Indeed, by means of Corollary~\ref{C:previo} and Figure~\ref{nonrepellor}, for any \textcolor{black}{non repeller vertex}  $\pi'\in\mathcal{C}_{19},$  we can find a finite path $v'\in\{a,b\}^L$, $L\in\N$, such that $v'(\pi')=\pi_{19,50}$ and then the claim follows by applying  Lemma~\ref{C-rutafinita}.

Let $T=(\lambda,\pi_{19,50})$ be minimal. Then, the orbit of $T$ under the Rauzy operator is infinite  by Proposition~\ref{P:minimalidadimplicainducible} and if we write
$$R^j(T)=(\lambda^{j},\pi_{j})\quad\textrm{ for any }\quad j\in\N,$$
then $\pi^{j}$ is irreducible, and by Theorem~\ref{T-cone}:
$$\lambda\in\bigcap_{i\in\N}M_{v_1}(\pi_1)M_{v_2}(\pi_2)\dots M_{v_i}(\pi_i)\Lambda^4.$$

For any pair of natural numbers $j>i$, we write  
\begin{equation}
S(i,j):=M_{v_i}(\pi_i)M_{v_{i+1}}(\pi_{i+1})\dots M_{v_{j}}(\pi_{j}). \label{D:Sij} 
\end{equation}
\textcolor{black}{Since clearly all the entries of the matrices $M_{v_s}(\pi_s)$ are nonnegative  and there are no zero rows for $i\leq s\leq j$, we have that $M_{v_s}(\pi_s)\Lambda^4 \subseteq \Lambda^4,$ and so does $S(i,j)\Lambda^4 \subseteq \Lambda^4$; thus, by Proposition~\ref{Lipschitz1} we deduce that  $\ddp(M_{v_s}(\pi_s)\lambda, M_{v_s}(\pi_s)\mu)\leq \ddp(\lambda,\mu)$ and $\ddp(S(i,j)\lambda,S(i,j)\mu)\leq \ddp(\lambda,\mu)$ for any $\lambda,\mu\in\R_+^4.$}

Note that, by Lemma~\ref{L:nociclosmenores}, the permutation $\pi_{19,50}$ appears infinitely many times in the sequence of permutations $(\pi_j)_{j\in\N}$. 

\textcolor{black}{Taking into account that we can arrive to $\pi_{19,50}$ under the three paths appearing in (\ref{Eq:caminos}), we divide the proof of the non-uniquely minimal ergodicity into the following cases.}

(i) In the first one, we assume that the Rauzy process admits an increasing sequence $(n_k)_{k\geq 1}$ such that:
\begin{itemize}
\item $\pi_{n_1}=\pi_{19,50}$, $R^{n_1}(T)$ is of type $a$ and $S(1,n_1-1)\geq 0$ without null rows.



 \item $\pi_{n_k}=\pi_{19,50}$ and $R^{n_k}=(\lambda^{n_k},\pi_{n_k})$ is of type $a$.  
 \item The induction of $R^{n_k}(T)$ follows a path $z_k\in\{v,\widetilde{w}\}$  given in (\ref{Eq:caminos}), (\ref{Eq:wmodify}),  and $n_{k+1}-n_k$ is bigger or equal  to $\dim (z_k)$ (the number of components of $z_k$)
Moreover, we can assume that the regular matrix
$S(n_k+\dim z_k,n_{k+1})\geq 0.$ 
 \end{itemize}

 Observe that
\begin{eqnarray*}
M_{v_1}(\pi_1)M_{v_2}(\pi_2)\dots M_{v_{n_k}}(\pi_{n_k})&=&S(1,n_k)\\
&=& S(1,n_1-1)
S(n_1,n_1+\dim z_1-1)S(n_1+\dim z_1,n_2)\\
&& S(n_2,n_2+\dim z_2-1)S(n_2+\dim z_2 ,n_3)\\
&& \dots\dots\\
&& S(n_{k-1},n_{k-1}+\dim z_{k-1}-1 ) 
S(n_{k-1}+\dim z_{k-1},n_{k})\\
&=&  
S(1,n_1-1)\\
&& P(z_1)S(n_1+\dim z_1,n_2)\\
&& P(z_2)S(n_2+\dim z_2 ,n_3)\\
&& \dots\dots\\
&&  P(z_{k-1}) 
S(n_{k-1}+\dim z_{k-1},n_{k}).
\end{eqnarray*}
Now, applying Proposition~\ref{Lipschitz1} \textcolor{black}{to the linear maps $S(1,n_1)$ and $S(n_j+\dim z_j, n_{j+1}), \, j\in\{1,\ldots,k-1\}$}, and Lemma~\ref{L:lipschitz} \textcolor{black}{to $P(z_j)$}, we deduce:
$$\ddp(S(1,n_k)\lambda,S(1,n_k)\mu)\leq \left(\frac{2}{3}\right)^{\textcolor{black}{k-1}}\ddp(\lambda,\mu) \quad\textrm{ for any }
\quad \lambda,\mu\in\R_{+}^{4}.$$

Therefore, the diameter of $\bigcap_{i\in\N}M_{v_1}(\pi_1)M_{v_2}(\pi_2)\dots M_{v_i}(\pi_i)\Lambda^n$ with respect to $\ddp$ is zero,
 which implies that if $\lambda,\mu\in\mathcal{C}(\mathcal{G}^T)$ then they are linearly dependent. Finally, $T$ is uniquely ergodic \textcolor{black}{according to Theorem~\ref{T-cone}.}

(ii) In the second case, assume that that the Rauzy induction gives an increasing sequence $(m_k)_{k\geq 1}$ such that now: 
$\pi_{m_1}=\pi_{19,50}$, $R^{m_1}(T)$ is of type $b$ and $S(1,m_1-1)>0$; $\pi_{m_k}=\pi_{19,50}$ and $R^{m_k}=(\lambda^{m_k},\pi_{m_k})$ is of type $b$; and the induction of $R^{m_k}(T)$ follows a path $z_k$ of the form described in~(\ref{Eq:zmodify}),  with  $m_{k+1}-m_k\geq \dim (z_k)$ (the number of components of $z_k$). Repeating the reasoning step by step, we arrive again to $\ddp(S(1,m_k)\lambda,S(1,m_k)\mu)\leq \left(\frac23\right)^{\textcolor{black}{k-1}}\ddp(\lambda,\mu)$ for any $\lambda,\mu\in\R_+^4,$ which implies that $T$ is uniquely ergodic.

\end{proof}

\section{Proof of Main Theorem}

The case $n=1$ is trivial. 
The absence of transitive  $(2,k)$-\iets, $k\geq 1$,  is stated in Theorem~\ref{T:nominimal2}. However in the oriented case we know the existence of minimality, but by Theorem~\ref{T:Weyl} minimality and unique ergodicity are linked. Also, assume the existence of an oriented transitive non uniquely  ergodic $2$-\iet, 
$T=(\lambda, \pi)$; then $\pi=(2,1)$ and the absence of minimality implies, according to Proposition~\ref{P:minimalidadimplicainducible} and having into account that $\pi_m$ is always irreducible, that the Rauzy process associated to $T$ must stop for some $p\geq 0$ arriving to $R^p(T)=T^{(p)}=(\lambda^{p},(2,1)):\bigcup_{i=1}^2 I^{(p)}_i\to \bigcup_{i=1}^2 I^{(p)}_i$  with $\lambda^{p}_{1}=\lambda^{p}_2.$ Although Rauzy operator has stopped, we can take the Poincaré return map,
$S=(\mu,\tau)$, of
$T^{(p)}$ on ${\Dom(T^{(p)})\menos I_2^{(p)}}$ and we   would obtain a  transitive $1$-\iet, a contradiction.

We proceed with $(3,k)$-\iets\ with $k\geq 0$. Observe that there are not minimal $(3,k)$-\iet, $k\geq 1$, by Theorem~\ref{T:nominimal3}, nor minimal non uniquely ergodic oriented $3$-\iets\ by Theorem~\ref{T:Weyl}. 

Assume that $T$ is a transitive non uniquely ergodic $(3,k)$-\iet\ with $k\geq 0$ and with $T$ being not minimal. Then,
\textcolor{black}{according to Proposition~\ref{P:minimalidadimplicainducible}, either $T$ is not infinitely inducible by Rauzy operator $R$, or $R^{m_0}(T)$ has as combinatorial data a non-irreducible permutation $\pi_{m_0}$ for some $m_0\in\mathbb N$}.

In the first case, if the Rauzy process stops at some  $p\in\N\cup\{0\}$, we find that we cannot continue with the process, which implies that
 for $T^{(p)}=(\lambda^{p},\pi_{p})=R^p(T)$ it must happen that $\lambda^{p}_{|\pi_p|^{-1}(3)}=\lambda^{p}_3.$
Thus, if the domain of $T^{(p)}$ is $\bigcup_{i=1}^3 I^{(p)}_i$, then  $T^{(p)}(I^{(p)}_{|\pi_{p}|^{-1}(3)})=I^{(p)}_3$ and it is easy
to see that the Poincaré return map, $S=(\mu,\tau)$, of
$T^{(p)}$ on ${\Dom(T^{(p)})\menos I_3^{(p)}}$ is a transitive non uniquely ergodic $(2,k')$-\iet, $0\leq k'\leq k$,
which preserves the number of independent invariant measures of $T$ but this is is a contradiction since there are not transitive non uniquely ergodic $(2,k')$-\iets\ as we have shown in the first paragraph of this proof.

In the second case, that only can occur with $k\geq 1$, we have that  the transitive \iet\ $R^{m_0}(T)=T^{(m_0)}$ has as combinatorial data a non-irreducible permutation $\pi_{m_0}$ for some $m_0\in\mathbb N$. Then $T^{(m_0)}({J_{m_0}})\subseteq J_{m_0}$ for some union of  intervals $J_{m_0}\subsetneq \Dom(T^{(m_0)})$ which contradicts the transitivity of $T^{(m_0)}$.

Now, let   $T$ be   a minimal $(4,k)$-\iet\ . Then, it is uniquely ergodic by directly applying  
Propositions \ref{clases1element}, \ref{clases11Y18element}, \ref{clases5Y15}, \ref{P:clase8}, \ref{P:clase12}, \ref{P:clase14}, \ref{P:clase20}, \ref{P:clase13}, and Theorem~\ref{T:clase19NoExisteMinimalNoUnicamenteErgodica}. Thus we must finally analyze the case of 
a  transitive non uniquely ergodic $(4,k)$-\iet\ with $k\geq 1$. A simple adaptation of the arguments for $(3,k)$-\iets\ yields the result and the proof is finished.

\section{Conclusions and further investigation}

In this paper, as main result, we have proved that there are no minimal non uniquely ergodic $(4,k)$-IETs, $1\leq k\leq 4$. To do it, we have classified the $4$-IETs into twenty Rauzy classes, eighteen of them consisting of non-oriented IETs. By a direct approach, in which Theorem~\ref{Th:hlamda} on inverse and conjugated IETs plays an important role in order to simplify the cases, for the non-oriented case we have shown that, a priori, the minimality only can occur into 
the class $\mathcal{C}_{19}$. Also, for this candidate class to contain minimal $(4,k)$-IETs, we have detected a subgraph, 
see Figure~\ref{nonrepellor}, consisting of non repeller vertices in $\mathcal{C}_{19}$. Finally, to establish the uniquely ergodic character of these minimal maps, we have used, following the ideas of~\cite{viana}, the so-called projective metric to demonstrate that all the elements of the cone $\Lambda^4$ are directed to a unique direction by the Rauzy matrices, which gets the unique ergodicity of those minimal maps. 

It is worth mentioning that in~\cite{NoUnicamenteErgodica10}, the authors built explicitly minimal non uniquely ergodic $(10,k)$-IETs for any $1\leq k\leq 10$, answering an open question posed in \cite{gutierrez4b}. Moreover, in~\cite{NoUnicamenteErgodica10}, it is proposed the  problem of determining the minimal integer $N$ for which there exists a minimal non uniquely ergodic $(N,k)$-IET with $k\geq 1$. In the present study, we have found that $N>4$. We advance that, in fact, it is possible to obtain minimal non uniquely ergodic $(n,k)$-IETs whenever $n= 6$, so $N\in\{5,6\}$, and we conjecture that $N=6$, see \cite{LineroSolerN=6}. However the ideas introduced in this paper  to prove the absence of minimal non uniquely ergodic $(4,k)$-\iets, $k\geq 1$, are hardly adaptable to the case of $(5,k)$-\iets\ since, in this case, there are $69$ Rauzy classes and one of them contains $530$ permutations while the largest Rauzy class for $n=4$ has $59$ permutations.

\section{Acknowledgements}

This paper has been partially supported by Grant  MTM2017-84079-P (AEI/FEDER, UE), Ministerio de Ciencias, Innovaci\'{on} y Universidades,
Spain.

%
%

\def\cprime{$'$}

\end{document}